\documentclass[11pt,reqno]{amsart}
\usepackage{amsmath,amsfonts,amssymb}
\usepackage[cp1251]{inputenc}
\usepackage{mathrsfs}

\def\wh{\widehat}
\def\wt{\widetilde}
\def\R{\mathbb R}
\def\C{\mathbb C}

\def\u{\mathbf u}
\def\e{\mathbf e}

\def\x{\mathbf x}
\def\v{\mathbf v}

\def\s{\mathbf s}

\def\k{\mathbf k}
\def\w{\mathbf w}

\def\D{\mathbf D}

\def\1{\mathbf 1}

\def\H{\mathfrak H}

\def\NN{\mathfrak N}
\def\bxi{\boldsymbol \xi}

\def\bt{\boldsymbol \theta}

\def\1{\bold 1}

\def\eps{\varepsilon}

\def\Ker{\mathrm{Ker}\,}

\def\ge{\geqslant}

\theoremstyle{theorem}

\numberwithin{equation}{section}

\theoremstyle{plain}
\newtoks\thehProclaim
\newtheorem*{Proclaim}{\the\thehProclaim}

\begin{document}

\title[ Homogenization of high order elliptic operators]
{ Homogenization of high order elliptic operators  with periodic coefficients}

\author{A.~A.~Kukushkin and T.~A.~Suslina}

\address{St.~Petersburg State University\\
Universitetskaya nab. 7/9\\
St.~Petersburg\\
 199034, Russia}

\email{beslave@gmail.com}

\email{t.suslina@spbu.ru}



\thanks{Supported by RFBR (grant no.~14-01-00760) and St.~Petersburg State University (project no.~11.38.263.2014).}

\begin{abstract}
In $L_2(\R^d;\C^n)$, we study a selfadjoint strongly elliptic operator $A_\eps$ of order $2p$
given by the expression $b(\D)^* g(\x/\eps) b(\D)$, $\eps >0$. Here
$g(\x)$ is a bounded and positive definite $(m\times m)$-matrix-valued function in $\R^d$;
it is assumed that $g(\x)$ is periodic with respect to some lattice. Next,
$b(\D)=\sum_{|\alpha|=p}^d b_\alpha \D^\alpha$ is a differential operator of order $p$ with constant coefficients;
$b_\alpha$ are constant $(m\times n)$-matrices. It is assumed that $m\ge n$ and that the symbol $b({\boldsymbol \xi})$
has maximal rank. For the resolvent $(A_\eps - \zeta I)^{-1}$ with $\zeta \in \C \setminus [0,\infty)$,
we obtain approximations in the norm of operators in $L_2(\R^d;\C^n)$ and in the norm of operators acting from
$L_2(\R^d;\C^n)$ to the Sobolev space $H^p(\R^d;\C^n)$, with error estimates depending on $\eps$ and $\zeta$.
\end{abstract}

 \maketitle

\section*{Introduction}

\subsection{Operator-theoretic approach to homogenization theory}

The paper concerns homogenization theory of periodic differential operators (DOs).
It is a wide area of theoretical and applied science.
A broad literature is devoted to homogenization problems; first of all, we mention the books
{[}BeLP{]}, {[}BaPa{]}, {[}ZhKO{]}.

In a series of papers by M.~Sh.~Birman and T.~A.~Suslina {[}BSu1-4{]},
an operator-theoretic approach to homogenization problems was suggested.
By this approach, a wide class of selfadjoint matrix second order operators $\mathcal A$
acting in $L_2(\mathbb{R}^{d};\mathbb{C}^{n})$ was studied.
It was assumed that $\mathcal A$ admits a factorization
$$
\mathcal{A}=b(\mathbf{D})^{*}g(\mathbf{x})b(\mathbf{D}).
\eqno(0.1)
$$
Here $g(\mathbf{x})$ is a bounded and uniformly positive definite
$(m\times m)$-matrix-valued function periodic with respect to some lattice $\Gamma\subset\mathbb{R}^{d}$.
By $\Omega$ we denote the elementary cell of the lattice $\Gamma$.
Next, $b(\mathbf{D})$ is an $(m\times n)$-matrix homogeneous first order DO.
It is assumed that $m\geqslant n$ and that the symbol $b(\boldsymbol{\xi})$
has rank $n$ for any $0 \neq \boldsymbol{\xi}\in \R^d$.
Under the above assumptions, $\mathcal{A}$ is strongly elliptic.
The simplest example of the operator (0.1) is the scalar elliptic operator $-{\rm div}\, g(\x) \nabla$
(the acoustics operator); the operator of elasticity theory also can be represented in the form (0.1).
These and other examples were considered in [BSu1,3,4]
in detail.

Let $\varepsilon>0$ be a small parameter. Denote
$F^{\varepsilon}(\mathbf{x}):=F(\eps^{-1}{\mathbf{x}})$.
Consider the operator
$
\mathcal{A}_{\varepsilon}=b(\mathbf{D})^{*}g^{\varepsilon}(\mathbf{x})b(\mathbf{D})
$
whose coefficients oscillate rapidly as $\varepsilon\to0$.

In {[}BSu1{]}, it was shown that
the resolvent $\left(\mathcal{A}_{\varepsilon}+I\right)^{-1}$ converges in the
operator norm in $L_2(\R^d;\C^n)$ to the resolvent of the \textit{effective operator}
$
\mathcal{A}^{0}=b(\mathbf{D})^{*}g^{0}b(\mathbf{D}),
$
as $\eps \to 0$. Here $g^{0}$ is a constant \emph{effective matrix}.
It was proved that
$$
\left\Vert \left(\mathcal{A}_{\varepsilon}+I\right)^{-1}-\left(\mathcal{A}^{0}+I\right)^{-1}\right\Vert _{L_2(\R^d)\to L_2(\R^d)}\leqslant C\varepsilon.
\eqno(0.2)
$$
In {[}BSu2,3{]}, a more accurate approximation of the resolvent
$\left(\mathcal{A}_{\varepsilon}+I\right)^{-1}$ in the operator norm in $L_2(\mathbb{R}^{d};\mathbb{C}^{n})$
with error $O(\varepsilon^{2})$ was obtained. In {[}BSu4{]}, an approximation of the resolvent
$\left(\mathcal{A}_{\varepsilon}+I\right)^{-1}$ in the norm of operators acting from
$L_2(\mathbb{R}^{d};\mathbb{C}^{n})$ to the Sobolev space $H^1(\mathbb{R}^{d};\mathbb{C}^{n})$ was found:
$$
\left\Vert \left(\mathcal{A}_{\varepsilon}+I\right)^{-1}-\left(\mathcal{A}^{0}+I\right)^{-1}
-\varepsilon\mathcal{K}({\varepsilon})\right\Vert _{L_2(\R^d) \to H^1(\R^d)}\leqslant C \varepsilon.
\eqno(0.3)
$$
Here $\mathcal{K}({\varepsilon})$ is the so called \textit{corrector}.
The operator  $\mathcal{K}({\varepsilon})$ contains rapidly oscillating factors, and so depends on $\eps$; herewith,
$\| \mathcal{K}({\varepsilon})\|_{L_2 \to H^1}=O(\eps^{-1})$.

Estimates (0.2) and (0.3) are order-sharp; the constants are controlled explicitly in terms of the problem data.
Such results are called the \textit{operator error estimates} in homogenization theory.
The method of [BSu1-4] is based on the scaling transformation,
the direct integral expansion for the periodic operator $\mathcal{A}$
(the Floquet-Bloch theory) and the analytic perturbation theory.
It turned out that the resolvent of the operator
$\mathcal{A}_{\varepsilon}$ can be approximated in terms of the threshold characteristics of the operator
$\mathcal{A}$ at the bottom of the spectrum. In this sense, the homogenization procedure
is a \emph{spectral threshold effect}.

We also mention the recent paper [Su], where the analogs of estimates (0.2), (0.3)
for the resolvent $\left(\mathcal{A}_{\varepsilon} - \zeta I\right)^{-1}$ in an arbitrary point
$\zeta \in \C \setminus \R_+$ were obtained;
the error estimates found in [Su] depend on $\eps$ and $\zeta$.

A different approach to operator error estimates (the modified method of the first approximation)
was suggested by V.~V.~Zhikov; by this method, in [Zh] and [ZhPas], the analogs of estimates (0.2) and (0.3)
for the acoustics operator and the elasticity operator were obtained.

A homogenization problem for periodic elliptic DOs of high even order is of separate interest.
In the paper {[}V{]} by N.~A.~Veniaminov, the method suggested in {[}BSu1{]} was developed for such operators.
A homogenization problem for the operator
$$
\mathcal{B}_{\varepsilon}=(\mathbf{D}^{p})^{*}g^{\varepsilon}(\mathbf{x}) \mathbf{D}^{p}
\eqno(0.4)
$$
was studied. Here $g(\mathbf{x})$ is a symmetric uniformly
positive definite and bounded tensor of order $2p$ periodic with respect to the lattice
$\Gamma$. The operator (0.4) with $p=2$ arises in the theory of elastic plates (see [ZhKO]).

The effective operator for $\mathcal{B}_\eps$ is given by
$
\mathcal{B}^{0}=(\mathbf{D}^{p})^{*}g^{0}\mathbf{D}^{p},
$
where $g^{0}$ is an \emph{effective tensor}. In [V], the following analog of estimate (0.2) was proved:
$$
\left\Vert \left(\mathcal{B}_{\varepsilon}+I\right)^{-1}-\left(\mathcal{B}^{0}+I\right)^{-1}\right\Vert _{L_2(\R^d) \to L_2(\R^d)}\leqslant
{C}\varepsilon.
\eqno(0.5)
$$

\subsection{Main results}

We study a more general class of high order elliptic periodic DOs than  (0.4).
Consider the operator
$$
A= A(g) =b(\mathbf{D})^{*}g(\mathbf{x}) b(\mathbf{D}),
\eqno(0.6)
$$
where $g(\mathbf{x})$ is uniformly positive definite and bounded
$\left(m\times m\right)$-matrix-valued function periodic with respect to the lattice $\Gamma$,
and $b(\mathbf{D})$ is an $(m\times n)$-matrix homogeneous DO of order $p$.
The precise definition of the operator (0.6) is given in Subsection~4.1.
We study homogenization problem for the operator $A_{\varepsilon}=A(g^\eps )$.

\textit{Main results of the paper} are approximations of the resolvent
\hbox{$(A_{\varepsilon} - \zeta I)^{-1}$}, where $\zeta \in \C \setminus [0,\infty)$,
in various operator norms with twoparametric error estimates
(depending on $\eps$ and $\zeta$).
Theorem~8.1 gives approximation of the resolvent in the operator norm in $L_2(\R^d;\C^n)$ (the analog of estimate (0.5)):
$$
\left\Vert \left(A_{\varepsilon}-\zeta I\right)^{-1}-\left(A^{0}- \zeta I\right)^{-1}\right\Vert _{L_2(\R^d) \to L_2(\R^d)}\leqslant {C}_1(\zeta) \varepsilon;
\eqno(0.7)
$$
in Theorem~8.2, approximation  of the resolvent in the
"energy" \  norm (i.~e., in the norm of operators acting from
$L_2(\R^d;\C^n)$ to the Sobolev space $H^p(\R^d;\C^n)$) is obtained:
$$
\left\Vert \left(A_{\varepsilon}- \zeta I\right)^{-1}-\left(A^{0}- \zeta I\right)^{-1}-\varepsilon^{p}K(\zeta;\varepsilon)\right\Vert _{L_2(\R^d) \to H^p(\R^d)}\leqslant {C}_2(\zeta)\varepsilon.
\eqno(0.8)
$$
It is shown that the effective operator $A^{0}$ has the same structure as the initial operator:
$A^{0}=b(\mathbf{D})^{*}g^{0}b(\mathbf{D}).$
The corrector $K(\zeta;\eps)$ contains rapidly oscillating factors;
herewith, $\|K(\zeta;\eps)\|_{L_2 \to H^p} = O(\eps^{-p})$.
The dependence of $C_1(\zeta)$ and $C_2(\zeta)$ on the spectral parameter $\zeta$
is searched out.

Besides estimate (0.8), we obtain approximation of the operator $g^\eps b(\D)\left(A_{\varepsilon}- \zeta I\right)^{-1}$
(corresponding to the "flux") in the norm of operators acting from  $L_2(\R^d;\C^n)$ to $L_2(\R^d;\C^m)$.

In the general case, the corrector $K(\zeta;\varepsilon)$
contains an auxiliary smoothing operator. We distinguish a condition that allows to use the
standard corrector (which does not involve the smoothing operator); see Theorem~8.6.

\subsection{The method}

The method is further development of the operator-theoretic approach.

By the scaling transformation, the dependence on the parameter
$\varepsilon$ is transferred from the coefficients of the operator to the resolvent point.
Namely, we have the unitary equivalence:
$$
\begin{aligned}
\left(A_{\varepsilon}- \zeta I\right)^{-1} &\sim\varepsilon^{2p}\left(A- \zeta \varepsilon^{2p}I\right)^{-1},
\\
\left(A^{0}- \zeta I\right)^{-1} &\sim\varepsilon^{2p}\left(A^{0}- \zeta \varepsilon^{2p}I\right)^{-1}.
\end{aligned}
$$
Then estimate (0.7) is reduced to the inequality
$$
\left\Vert \left(A- \zeta \varepsilon^{2p}I\right)^{-1}-\left(A^{0}- \zeta \varepsilon^{2p}I\right)^{-1}\right\Vert _{L_2(\R^d)\to L_2(\R^d)}\leqslant
{C}_1(\zeta) \varepsilon^{1-2p}.
\eqno(0.9)
$$
To prove (0.8), we use (0.7) and the auxiliary inequality
$$
\left\Vert A_{\varepsilon}^{1/2}\bigl(\left(A_{\varepsilon}- \zeta I\right)^{-1}-\left(A^{0}- \zeta I\right)^{-1}-\varepsilon^{p}K({\zeta;\varepsilon})\bigr)\right\Vert _{L_2(\R^d) \to L_2(\R^d)}\leqslant C_3(\zeta)\varepsilon.
$$
By the scaling transformation, the latter is equivalent to
$$
\left\Vert A^{1/2}\bigl(\left(A-\zeta \varepsilon^{2p}I\right)^{-1}-\left(A^{0}- \zeta\varepsilon^{2p}I\right)^{-1}-
\widetilde{K}(\zeta;{\varepsilon}) \bigr)\right\Vert _{L_2 \to L_2}\leqslant C_3(\zeta) \varepsilon^{1-p}.
\eqno(0.10)
$$

Estimates  (0.9) and (0.10) with arbitrary $\zeta \in \C \setminus [0,\infty)$ are deduced from the similar inequalities with $\zeta =-1$
by suitable identities for the resolvents; this trick is borrowed from~[Su].
Therefore, main considerations concern the case where $\zeta =-1$.

The operator $A$ is expanded in the direct integral of the operators $A\left(\mathbf{k}\right)$ acting in
$L_{2}\left(\Omega;\mathbb{C}^{n}\right)$ and depending on the parameter $\mathbf{k}$ (the \textit{quasimomentum}).
The operator $A\left(\mathbf{k}\right)$
is given by the expression $b(\mathbf{D}+\mathbf{k})^{*}g(\mathbf{x})b(\mathbf{D}+\mathbf{k})$
with periodic boundary conditions. As in {[}BSu1{]}, we distinguish the onedimensional parameter $t=\left|\mathbf{k}\right|$,
with respect to which the family $A(\mathbf{k})$ is a \emph{polynomial operator pencil} of order $2p$.
In {[}V{]}, the abstract scheme was developed for such operator pencils.
Using it, we prove estimate  (0.9) (with $\zeta=-1$). To check (0.10), we develop further the abstract scheme for
the polynomial pencils by analogy with {[}BSu4{]}.

\subsection{Plan of the paper}

The paper consists of 8 sections. \S~1--3 are devoted to the abstract scheme.
In \S~1, we describe a factorized operator family $A(t)=X(t)^{*}X(t)$ and introduce the spectral germ.
In \S~2, the results obtained in {[}V{]}
(threshold approximations and the principal order approximation for the resolvent $(A(t)+\eps^{2p}I)^{-1}$) are described.
In \S~3, by further development of the abstract method, we obtain approximation for the resolvent
 $(A(t)+\eps^{2p}I)^{-1}$ with the corrector taken into account.
In \S~4, the class of periodic differential operators $A$ acting in $L_2(\R^d;\C^n)$ is introduced;
 the operator $A$ is expanded in the direct integral of the operators $A(\k)$ acting in $L_2(\Omega;\C^n)$.
In \S~5, the family $A(\k)$ is studied with the help of abstract results, the effective operator is introduced, and the properties
of the effective matrix are described. In \S~6,  we obtain approximation for the resolvent
$(A(\k)+\eps^{2p}I)^{-1}$ by applying abstract theorems.
In \S~7,  using the results of \S~6 and the direct integral expansion for $A$,
we deduce theorems about approximation of the resolvent $(A+\eps^{2p}I)^{-1}$;
next, with the help of appropriate identities for the resolvents, these theorems are carried over to the case of
the resolvent  $(A-\zeta \eps^{2p}I)^{-1}$ at arbitrary point $\zeta \in \C \setminus [0,\infty)$.
We distinguish the condition under which the smoothing operator in the corrector can be removed.
In \S~8, using the scaling transformation, we deduce the main results (approximations of the resolvent
 $(A_\eps - \zeta I)^{-1}$ in various operator norms) from the estimates of \S~7.

\subsection{Notation}

Let $\mathfrak{H}$ and $\mathfrak{G}$ be separable Hilbert spaces.
The symbols $\left\Vert \cdot\right\Vert _{\mathfrak{H}}$ and $\left(\cdot,\cdot\right)_{\mathfrak{H}}$ stand for the norm
and the inner product in $\mathfrak{H}$, respectively; the symbol $\left\Vert \cdot\right\Vert _{\mathfrak{H} \to \mathfrak{G}}$
denotes the norm of a linear continuous operator acting from $\mathfrak{H}$ to $\mathfrak{G}$.
Sometimes we omit the indices. If $G$ is a linear operator acting from $\mathfrak{H}$ to $\mathfrak{G}$,
then $\mathrm{Dom}\, G$ and $\mathrm{Ker}\, G$ denote its domain and kernel, respectively.
If $\NN$ is a subspace of $\H$, then $\NN^\perp$ denotes its orthogonal complement.

The symbols $\langle \cdot ,\cdot \rangle$ and $\vert \cdot \vert$ stand for the inner product and the norm in $\mathbb{C}^n$,
respectively; $\mathbf{1}_n$ is the identity $(n\times n)$-matrix. If $a$ is an $(m\times n)$-matrix,
then $\vert a\vert$ denotes the norm of $a$ viewed as an operator acting from $\mathbb{C}^n$ to $\mathbb{C}^m$.

The $L_{q}$-classes of $\mathbb{C}^{n}$-valued functions in a domain $\mathcal{O}\subset\mathbb{R}^{d}$
are denoted by $L_{q}(\mathcal{O};\mathbb{C}^{n})$, $1\leqslant q\leqslant\infty$.
The Sobolev classes of $\mathbb{C}^{n}$-valued functions (in a domain $\mathcal{O}\subseteq\mathbb{R}^{d}$)
of order $s$ and summability degree $q$ are denoted by $W_{q}^{s}(\mathcal{O};\mathbb{C}^{n})$.
If $q=2$, we use the notation $H^{s}(\mathcal{O};\mathbb{C}^{n})$,
$s\in\mathbb{R}$. If $n=1$, we write $L_{q}({\mathcal O})$,
$W_{q}^{s}(\mathcal{O})$, $H^{s}(\mathcal{O})$, but sometimes we use such abbreviated notation also for
the spaces of vector-valued or matrix-valued functions.

We use the notation $\x = (x_1,\dots,x_d)\in \R^d$, $iD_j = \partial_j = {\partial}/{\partial x_j}$,
$j=1,\dots,d$, $\mathbf{D}=-i{\nabla}= (D_1,\dots,D_d)$.
Next, if $\alpha=(\alpha_1,\dots,\alpha_d) \in {\mathbb Z}_+^d$ is a multiindex and $\k \in \R^d$, then
$|\alpha|= \sum_{j=1}^d \alpha_j$, $\k^\alpha = k_1^{\alpha_1} \cdots k_d^{\alpha_d}$, and
$\D^\alpha = D_1^{\alpha_1} \cdots D_d^{\alpha_d}$.
If $\alpha$ and $\beta$ are two multiindices, we write $\beta \leqslant \alpha$ if
$\beta_j \leqslant \alpha_j$, $j=1,\dots,d$; the binomial coefficients are denoted by
$C_\alpha^\beta= C_{\alpha_1}^{\beta_1}\cdots C_{\alpha_d}^{\beta_d}$.

We denote $\R_+ = [0,\infty)$. By $C$, $B$, $c$, $\mathcal C$, and $\mathfrak C$ (possibly, with indices and marks) we denote
various constants in estimates.

\subsection{} From the talk "Operator error estimates for homogenization of high order elliptic equations"
 given by Svetlana Pastukhova at the XXVI Crimean Autumn Mathematical School (Batiliman, September 17--29, 2015),
 the authors knew that she had obtained close results by the afore mentioned different approach 
 (the modified method of the first approximation).

\section{Abstract scheme. The spectral germ}

\subsection{Polynomial pencils of the form $X(t)^*X(t)$}

Let $\mathfrak{H}$ and $\mathfrak{H}_{*}$ be complex separable Hilbert spaces.
Consider a family (polynomial pencil) of operators
$$
X(t)=\sum_{j=0}^{p}X_{j}t^{j},\quad t\in\mathbb{R},\quad p\in\mathbb{N},\quad p\geqslant2.
$$
(The case where $p=1$ was studied in [BSu1,2,4] in detail.)
The operators $X(t)$ and $X_{j}$ act from $\mathfrak{H}$ to $\mathfrak{H}_{*}$:
$$
X(t),\, X_{j}:\mathfrak{H}\to\mathfrak{H}_{*}.
$$
It is assumed that the operator $X_{0}$ is \emph{densely defined and closed},
and $X_{p}$ is \textit{bounded}. In addition, we impose the following condition on the domains of these operators.

\smallskip
\noindent
\textbf{Condition 1.1.}
$$
\mathrm{Dom}\,X(t)=\mathrm{Dom}\,X_{0}\subset\mathrm{Dom}\,X_{j}\subset\mathrm{Dom}\,X_{p}=\mathfrak{H},\quad j=1,\dots,p-1,\quad t\in\mathbb{R}.
$$

\smallskip
We also assume that the intermediate operators $X_{j}$ with $j=1,\dots,p-1$ are subordinate to $X_{0}$.

\smallskip
\noindent
\textbf{Condition 1.2.}
\textit{For $j=0,\dots,p-1$ and for any $u\in\mathrm{Dom}X_{0}$ we have
$$
\left\Vert X_{j}u\right\Vert _{\mathfrak{H}_{*}}\leqslant\widetilde{C}\left\Vert X_{0}u\right\Vert _{\mathfrak{H}_{*}},
\eqno(1.1)
$$
where $\widetilde{C}$ is some constant \emph{(}of course, $\widetilde{C}\geqslant1$\emph{)}.}

\smallskip
Note that for $j=0$ estimate (1.1) is trivial.
Under the above assumptions, the operator $X(t)$ is \textit{closed} on the domain $\mathrm{Dom}\,X(t)=\mathrm{Dom}\,X_{0}$.

From Condition~1.2 it follows that
$$
\mathrm{Ker}\, X_0 \subset\mathrm{Ker}\, X_{j},\quad j=1,\dots,p-1.
\eqno(1.2)
$$

Our main object is a family of selfadjoint nonnegative operators
$$
A(t)=X(t)^* X(t),\quad t\in\mathbb{R},
\eqno(1.3)
$$
in $\mathfrak{H}$. The operator (1.3) is generated by the closed quadratic form
$$
a(t)[u,u]=\left\Vert X(t)u\right\Vert _{\mathfrak{H}_{*}}^{2},\quad u\in\mathrm{Dom}\,X_{0}.
$$

We denote $A(0)=X_{0}^{*}X_{0}=:A_{0}$, and put
$$
\mathfrak{N}:=\mathrm{Ker}A_{0}=\mathrm{Ker}X_{0},\qquad\mathfrak{N}_{*}:=\mathrm{Ker}A_{0}^{*}=\mathrm{Ker}X_{0}^{*}.
$$
Let $P$ and $P_{*}$ be the orthogonal projections of
$\mathfrak{H}$ onto $\mathfrak{N}$ and of $\mathfrak{H}_{*}$ onto $\mathfrak{N}_{*}$, respectively.

\smallskip
\noindent
\textbf{Condition 1.3.}
\textit{Suppose that the point $\lambda_{0}=0$ is an isolated point in the spectrum of $A_{0}$, and}
$$
n:=\mathrm{dim}\,\mathfrak{N}<\infty,\quad n\leqslant n_{*}:=\mathrm{dim}\,\mathfrak{N}_{*}\leqslant\infty.
$$

\smallskip
\emph{The distance from the point $\lambda_0=0$ to the rest of the spectrum of $A_{0}$ is denoted by $d^{0}$.}
Let ${F}\left(t,s\right)$ denote the spectral projection of $A(t)$ for the interval $\left[0,s\right]$, and let $\mathfrak{F}(t,s):= {F}(t,s)\mathfrak{H}$.
Next, we fix a positive number $\delta$ such that $\delta\leqslant\mathrm{min}\{d^{0}/36,\,1/4\}$, and put
$$
t^{0}=\delta^{1/2}(\widehat{C})^{-1},
\eqno(1.4)
$$
where
$$
\widehat{C}=\mathrm{max}\left\{ \left(p-1\right)\widetilde{C},\left\Vert X_{p}\right\Vert \right\}.
\eqno(1.5)
$$
Here $\widetilde{C}$ is the constant from  (1.1). Note that $t^0\leqslant 1$.

In {[}V, Lemma~3.9{]}, it is deduced from Condition 1.2 that
$$
\left\Vert X_{0}f\right\Vert _{\mathfrak{H}_{*}}
\leqslant
2\left(\left\Vert X(t)f\right\Vert _{\mathfrak{H}_{*}}+\sqrt{\delta}\left\Vert f\right\Vert _{\mathfrak{H}}\right),
\quad f\in\mathrm{Dom}\,X_{0},\quad\left|t\right|\leqslant t^{0}.
\eqno(1.6)
$$
According to {[}V, Proposition~3.10{]}, for $\left|t\right|\leqslant t^{0}$ we have
$$
{F}(t,\delta)= {F}(t,3\delta),\quad\mathrm{rank}\,{F}\left(t,\delta\right)=n.
\eqno(1.7)
$$
It means that for $|t|\leqslant t^0$ the operator $A(t)$ has exactly $n$ eigenvalues (counted with multiplicities)
on the interval $[0,\delta]$, while the interval $(\delta, 3 \delta)$ is free of the spectrum.
For convenience, we denote
$$
{F}(t):= {F}(t,\delta),\quad\mathfrak{F}\left(t\right):=\mathfrak{F}(t,\delta).
$$

\subsection{The operators $Z$, $R$, and the spectral germ $S$}

We put $\mathcal{D}=\mathrm{Dom}\,X_{0}\cap\mathfrak{N}^{\perp}$.
Since the point $\lambda_0 =0$ is isolated in the spectrum of $A_{0}$, we can interpret
$\mathcal{D}$ as the Hilbert space with the inner product $\left(X_{0}\varphi,X_{0}\eta\right)_{\mathfrak{H}_{*}}$,
$\varphi,\eta\in\mathcal{D}$. Let $u\in\mathfrak{H}_{*}$. Consider the equation
$X_0^* (X_0 \psi - u) =0$ for $\psi\in\mathcal{D}$, understood in the weak sense:
$$
\left(X_{0}\psi,X_{0}\zeta\right)_{\mathfrak{H}_{*}}=\left(u,X_{0}\zeta\right)_{\mathfrak{H}_{*}},\quad\forall\zeta\in\mathcal{D}.
\eqno(1.8)
$$
The right-hand side of (1.8) is an antilinear continuous functional of $\zeta\in\mathcal{D}$.
Hence, there exists a unique solution $\psi$ such that
$\left\Vert X_{0}\psi\right\Vert _{\mathfrak{H}_{*}}\leqslant\left\Vert u\right\Vert _{\mathfrak{H}_{*}}$.
Now, let
$$
\omega\in\mathfrak{N},\quad u=-X_{p}\omega;
\eqno(1.9)
$$
in this case we denote the solution of equation (1.8) by $\psi\left(\omega\right)$.
Let $Z:\mathfrak{H}\to\mathfrak{H}$ be a bounded operator defined by
$$
Z\omega=\psi(\omega),\ \omega\in\mathfrak{N};\quad Z\varphi=0,\ \varphi\in\mathfrak{N}^{\perp}.
\eqno(1.10)
$$
In order to estimate the norm of $Z$, we write down (1.8) with $u=-X_{p}\omega$ and $\zeta=\psi\left(\omega\right)$:
$$
\left\Vert X_{0}\psi(\omega)\right\Vert _{\mathfrak{H}_{*}}^{2}=-\left(X_{p}\omega,X_{0}\psi(\omega)\right)_{\mathfrak{H}_{*}}\leqslant\left\Vert X_{0}\psi(\omega)\right\Vert _{\mathfrak{H}_{*}}\left\Vert X_{p}\omega\right\Vert _{\mathfrak{H}_{*}},
$$
whence
$$
\left(A_{0}\psi(\omega),\psi(\omega)\right)_{\mathfrak{H}}\leqslant\left\Vert X_{p}\right\Vert^{2}\left\Vert \omega\right\Vert _{\mathfrak{H}}^{2}.
$$
Recalling that $d^{0}\geqslant36\delta$ and $\psi\left(\omega\right)\in\mathfrak{N}^{\perp}$, we obtain
$$
36\delta\left\Vert \psi\left(\omega\right)\right\Vert _{\mathfrak{H}}^{2}\leqslant\left(A_{0}\psi(\omega),\psi(\omega)\right)_{\mathfrak{H}}\leqslant
\left\Vert X_{p} \right\Vert^{2} \left\Vert \omega\right\Vert _{\mathfrak{H}}^{2}.
$$
Hence,
$$
\left\Vert Z\right\Vert \leqslant (1/6) \delta^{-1/2} \left\Vert X_{p}\right\Vert.
\eqno(1.11)
$$

Now, we put
$$
\omega_{*}:=X_{0}\psi(\omega)+X_{p}\omega\in\mathfrak{N}_{*}
\eqno(1.12)
$$
and define the operator $R$ by the relations
$$
R:\mathfrak{N}\to\mathfrak{N}_{*},\quad R\omega=\omega_{*}.
\eqno(1.13)
$$
The operator $R$ can be also represented as
$$
R=P_{*}X_{p}|_{\mathfrak{N}}.
\eqno(1.14)
$$

By definition, the \textit{spectral germ} of the operator family $A(t)$ at the point
$t=0$ is the selfadjoint operator
$$
S=R^{*}R:\mathfrak{N}\to\mathfrak{N}.
\eqno(1.15)
$$
From (1.14) and (1.15) it follows that
$$
S=PX_{p}^{*}P_{*}X_{p}|_{\mathfrak{N}}.
$$
The germ $S$ is called \textit{nondegenerate} if $\mathrm{Ker}\,S=\left\{ 0\right\} $,
or, equivalently, $\mathrm{rank}\, R=n$.

\subsection{The analytic branches of the eigenvalues and the eigenvectors of $A(t)$}

By the analytic perturbation theory, in {[}V, Subsection 3.3{]}
important properties of the first $n$ eigenvalues and the corresponding eigenvectors of
$A(t)$ for sufficiently small $t$ were found.
Namely, for $\left|t\right|\leqslant t^{0}$ there exist real-analytic functions $\lambda_{j}(t)$ (the branches of the eigenvalues)
 and real-analytic $\mathfrak{H}$-valued functions $\varphi_{j}(t)$ (the branches of the eigenvectors) such that
$$
A(t)\varphi_{j}(t)=\lambda_{j}(t)\varphi_{j}(t),\quad j=1,\dots,n,\quad\left|t\right|\leqslant t^{0},
\eqno(1.16)
$$
and the set $\left\{ \varphi_{j}(t)\right\} _{j=1}^{n}$
forms an orthonormal basis in  $\mathfrak{F}(t)$.
If \hbox{$t_{*}\leqslant t^{0}$} is sufficiently small, then the following power series
expansions are convergent for  $|t|\leqslant t_*$ (see {[}V, Theorem~3.15{]}):
$$
\lambda_{j}(t)  =\gamma_{j}t^{2p}+\dots,\quad\gamma_{j}\geqslant0,\quad j=1,\dots,n,\quad\left|t\right|\leqslant t_{*},
\eqno(1.17)
$$
$$
\varphi_{j}(t)  =\omega_{j}+t\varphi_{j}^{(1)}+t^{2}\varphi_{j}^{(2)}+\dots,\quad j=1,\dots,n,\quad\left|t\right|\leqslant t_{*}.
\eqno(1.18)
$$
Herewith, the set $\left\{ \omega_{j}\right\} _{j=1}^{n}$ is an orthonormal basis in $\mathfrak{N}$.
The numbers $\gamma_j$ and the vectors $\omega_j$ are eigenvalues and eigenvectors of the spectral germ $S$, i.~e.,
$$
S\omega_{j}  =\gamma_{j}\omega_{j},\quad j=1,\dots,n.
$$
We have
\begin{align}
P  &=\sum_{j=1}^{n} \left(\cdot,\omega_{j}\right)_{\mathfrak{H}}\omega_{j},
\tag{1.19}
\\
SP  &=\sum_{j=1}^{n}\gamma_{j}\left(\cdot,\omega_{j}\right)_{\mathfrak{H}}\omega_{j}.
\tag{1.20}
\end{align}
According to [V, Subsection~3.3],  the elements $\varphi_{j}^{\left(i\right)}$ from  (1.18) satisfy
\begin{align}
\varphi_{j}^{\left(i\right)}  &\in\mathfrak{N},\quad j=1,\dots,n,\quad i=1,\dots,p-1;
\tag{1.21}
\\
\varphi_{j}^{\left(p\right)}-\psi\left(\omega_{j}\right)  &\in\mathfrak{N},\quad j=1,\dots,n.
\tag{1.22}
\end{align}

By (1.7) and (1.16),
\begin{align}
{F}(t) &=\sum_{j=1}^{n} \left(\cdot,\varphi_{j}\left(t\right)\right)_{\mathfrak{H}}\varphi_{j}\left(t\right),\quad\left|t\right|\leqslant t^{0},
\tag{1.23}
\\
A(t){F}(t) &=\sum_{j=1}^{n}\lambda_{j}(t)\left(\cdot,\varphi_{j}\left(t\right)\right)_{\mathfrak{H}}\varphi_{j}\left(t\right),\quad\left|t\right|\leqslant t^{0}.
\tag{1.24}
\end{align}

Substituting expansions  (1.17), (1.18) in (1.23), (1.24) and taking (1.19) and (1.20) into account,
we obtain the power series expansions
$$
\begin{aligned}
 {F}(t) &= P+t F_1 + \dots,\quad\left|t\right|\leqslant t_{*},
\\
A(t) {F}(t) &= t^{2p}SP+\dots,\quad\left|t\right|\leqslant t_{*},
\end{aligned}
$$
convergent for $|t|\leqslant t_*$ (where $t_{*}\leqslant t^{0}$ is sufficiently small).
However, we need not these expansions, but approximations for $F(t)$ and $A(t)F(t)$
with one or several first terms (\textit{threshold approximations}) and with error estimates on the whole interval $|t|\leqslant t^0$.

\section{The abstract scheme: threshold approximations}

This section contains main results of the abstract scheme obtained in [V].

\subsection{Auxiliary material}

We need a version of the resolvent identity in the case where
the domains of two operators may not coincide, but
the domains of the corresponding quadratic forms do coincide.
An appropriate version of the resolvent identity was found in [BSu1, Chapter~1, \S~2].

Let $a$ and $b$ be two closed nonnegative quadratic forms in $\mathfrak{H}$ defined on the common domain
$$
\mathfrak{d}:=\mathrm{Dom}\,a=\mathrm{Dom}\,b,
\eqno(2.1)
$$
which is dense in  $\mathfrak{H}$. The operators corresponding to the forms $a$ and $b$ are denoted by $A$ and $B$, respectively.
Consider the sesquilinear form
$$
a_{\gamma}\left[u,v\right]=a\left[u,v\right]+\gamma\left(u,v\right)_{\mathfrak{H}},\quad\gamma>0.
\eqno(2.2)
$$
The corresponding quadratic form is positive definite.
The form $b_{\gamma}$ is introduced in a similar way. The lineal $\mathfrak{d}$ is a Hilbert space
$\mathfrak{d}\left(a_{\gamma}\right)$ with respect to the inner product (2.2).
The norm in $\mathfrak{d}\left(a_{\gamma}\right)$ is denoted by $\left\Vert \cdot\right\Vert _{\mathfrak{d}}$:
$$
\left\Vert u\right\Vert _{\mathfrak{d}}=a_{\gamma}\left[u,u\right]^{1/2},\quad u\in\mathfrak{d}.
\eqno(2.3)
$$
By (2.1), the form $b_{\gamma}$ is continuous in $\mathfrak{d}\left(a_{\gamma}\right)$
and defines a norm equivalent to $\|\cdot\|_{\mathfrak{d}}$. Let $\alpha>0$ be a constant defined by
$$
\alpha^{2}=\sup_{0\neq u\in\mathfrak{d}}\frac{a_{\gamma}\left[u,u\right]}{b_{\gamma}\left[u,u\right]}.
$$
Consider the form $\mathfrak{t}=b-a$. Obviously, it is $a_{\gamma}$-continuous
and generates a selfadjoint operator $T_{\gamma}$ in  $\mathfrak{d}(a_{\gamma})$:
$$
\mathfrak{t}\left[u,v\right]=a_{\gamma}\left[T_{\gamma}u,v\right],\quad u,\, v\in\mathfrak{d}.
$$
We denote
$$
\Omega_{z}(A):=I+\left(z+\gamma\right)R_{z}(A)=\left(A+\gamma I\right)R_{z}(A),
$$
where $R_{z}(A)=(A - zI)^{-1}$ is the resolvent of $A$ at the point
$z\in\rho\left(A\right)$. (Here $\rho(A)$ stands for the resolvent set of $A$.)
Similar notation is introduced for the operator $B$. Then we have
$$
R_{z}(B)-R_{z}(A)=-\Omega_{z}(A)T_{\gamma}R_{z}(B),\quad z\in\rho(A)\cap\rho(B),
\eqno(2.4)
$$
see {[}BSu1, Chapter~1, \S~2{]}. From the definition of the norm (2.3) in
$\mathfrak{d}$ it follows directly that
\begin{align}
\left\Vert u\right\Vert _{\mathfrak{H}} & \leqslant\gamma^{-1/2}\left\Vert u\right\Vert _{\mathfrak{d}},\tag{2.5}\\
\left\Vert R_{z}(A)\right\Vert _{\mathfrak{H} \to \mathfrak{d}} & \leqslant\gamma^{-1/2}\left\Vert \Omega_{z}(A)
\right\Vert _{\mathfrak{H}\to \mathfrak{H}},\tag{2.6}\\
\left\Vert R_{z}(A)\right\Vert _{\mathfrak{d} \to \mathfrak{d}} & \leqslant\gamma^{-1}\left\Vert \Omega_{z}(A)
\right\Vert _{\mathfrak{H}\to \mathfrak{H}},\nonumber \\
\left\Vert \Omega_{z}(A)\right\Vert _{\mathfrak{d}\to \mathfrak{d}} & \leqslant1+\left|z+\gamma\right|\gamma^{-1}
\left\Vert \Omega_{z}(A)\right\Vert _{\mathfrak{H}\to \mathfrak{H}}.\tag{2.7}
\end{align}

\subsection{Estimates for the difference of the resolvents. Threshold approximations}

In this subsection, we formulate some results obtained in [V, Subsection~4.2].

Let $\Gamma_{\delta} \subset \C$ be a contour that envelopes the real interval
$\left[0,\delta\right]$ equidistantly at the distance $\delta$.
Recall that $\delta$ was chosen in Subsection~1.1.
Let $z \in \Gamma_\delta$, and let  $\left|t\right|\leqslant t^{0}$.
The roles of the operators $A$ and $B$ are played by $A_{0}=A(0)$ and $A(t)$.
For brevity, we write $R_{z}(t)$ in place of  $R_{z}(A(t))$ and $\Omega_{z}(t)$ in place of $\Omega_{z}(A(t))$.

We apply the scheme of Subsection 2.1 putting
$$
\gamma=\delta,\quad\mathfrak{d}=\mathrm{Dom}\, X_{0},\quad a\left[u,u\right]=\left\Vert X_{0}u\right\Vert _{\mathfrak{H}_{*}}^{2},\quad b\left[u,u\right]=\left\Vert X(t)u\right\Vert _{\mathfrak{H}_{*}}^{2}.
$$
It is easily seen that $\alpha\leqslant3$. Obviously, $|z| \leqslant 2 \delta$ for $z \in \Gamma_\delta$.
By (1.7),
$\|R_z(t)\|_{\H \to \H} \leqslant \delta^{-1}$ for $z \in \Gamma_\delta$ and $|t| \leqslant t^0$. Hence,
$$
\left\Vert \Omega_{z}(t)\right\Vert _{\mathfrak{H}\to \mathfrak{H}}\leqslant4,\quad z\in\Gamma_{\delta},\quad\left|t\right|\leqslant t^{0}.\eqno(2.8)
$$
Together with  (2.6) this yields
$$
\left\Vert R_{z}(0)\right\Vert _{\mathfrak{H} \to \mathfrak{d}}\leqslant4\delta^{-1/2},\quad z\in\Gamma_{\delta}.
\eqno(2.9)
$$
Similarly,  (2.7) and (2.8) imply that
$$
\left\Vert \Omega_{z}(0)\right\Vert _{\mathfrak{d} \to \mathfrak{d}}\leqslant13.
\eqno(2.10)
$$

Now, consider the form
$$
\begin{aligned}
 & \mathfrak{t}\left[u,u\right]=\left\Vert X(t)u\right\Vert ^{2}_{\H_*} - \left\Vert X_{0}u\right\Vert ^{2}_{\H_*}
\\
 & =2\mathrm{Re}\left(\left(tX_{1}+\dots+t^{p}X_{p}\right)u,X_{0}u\right)_{\H_*}+\left\Vert \left(tX_{1}+\dots+t^{p}X_{p}\right)u\right\Vert ^{2}_{\H_*}.
\end{aligned}
$$
In the space $\mathfrak{d}$ with the metric defined by the form $a_{\delta}$, this quadratic form generates the operator
$T_{\delta}=T_{\delta}(t)$ which can be represented as
$$
T_{\delta}(t)=\sum_{j=1}^{2p}t^{j}T_{\delta}^{(j)},
\eqno(2.11)
$$
where the operators $T_{\delta}^{(j)}$ do not depend on $t$. The norms of $T_{\delta}(t)$ and $T_{\delta}^{(j)}$
were estimated in {[}V, Propositions~4.3, 4.4{]}.

\smallskip
\noindent\textbf{Proposition 2.1.}
\textit{Let $t^0$ be given by} (1.4). \textit{Then}
\begin{align}
\left\Vert T_{\delta}(t)\right\Vert _{\mathfrak{d}\to \mathfrak{d}} &\leqslant  C_{\circ} |t|,\quad |t|\leqslant t^{0},
\tag{2.12}
\\
\left\Vert T_{\delta}^{(j)}\right\Vert _{\mathfrak{d}\to \mathfrak{d}} &\leqslant  \widetilde{B},\quad j=1,\dots,2p,
\tag{2.13}
\end{align}
\textit{where the constants $C_{\circ}$ and $\widetilde{B}$ are given by
\begin{align}
C_{\circ} & =  5\widehat{C}\delta^{-1/2}=5 (t^{0})^{-1},
\tag{2.14}
\\
\widetilde{B} & =  p\widetilde{C}^{2}+\left\Vert X_{p}\right\Vert ^{2}\delta^{-1},
\tag{2.15}
\end{align}
and $\widetilde{C}$ is the constant from \emph{(1.1)} \emph{(}obviously, $\widetilde{B}\geqslant1$\emph{)}.
}

\smallskip

In what follows, we will use the inequality equivalent to (2.12):
$$
\left|\left\Vert X(t)u\right\Vert _{\mathfrak{H}_{*}}^{2}-\left\Vert X_{0}u\right\Vert _{\mathfrak{H}_{*}}^{2}\right|\leqslant\left(\left\Vert X_{0}u\right\Vert _{\mathfrak{H}_{*}}^{2}+\delta\left\Vert u\right\Vert _{\mathfrak{H}}^{2}\right)C_{\circ}\left|t\right|,
\quad u \in \mathfrak{d},\quad |t| \leqslant t^0.
\eqno(2.16)
$$

On the basis of (2.4) in {[}V, (4.16){]} it was proved that
$$
\left\Vert R_{z}(t)-R_{z}(0)\right\Vert _{\mathfrak{H}\to \mathfrak{H}}\leqslant48C_{\circ}\delta^{-1} |t|,\quad\left|t\right|\leqslant t^{0},\quad z\in\Gamma_{\delta}.
\eqno(2.17)
$$

The estimates for the difference of the resolvents allow us to obtain \textit{threshold approximations} for the operators
$F(t)$ and $A(t)F(t)$. The spectral projection $F(t)$ admits the following representation (see {[}K{]})
$$
{F}(t)=-\frac{1}{2\pi i}\oint_{\Gamma_{\delta}}R_{z}(t)\,dz,
$$
where the integration over $\Gamma_{\delta}$ is in positive direction. Then
$$
{F}(t) - P=-\frac{1}{2\pi i}\oint_{\Gamma_{\delta}}\left(R_{z}(t) - R_{z}(0)\right)\,dz.
\eqno(2.18)
$$
Since the length of the contour $\Gamma_\delta$ is equal to $2\delta + 2\pi \delta$, from (2.17) and (2.18) we obtain
$$
\left\Vert {F}(t)-P\right\Vert _{\mathfrak{H}\to \mathfrak{H}} \leqslant
C_1 \left|t\right|,\quad |t|\leqslant t^{0},
\quad C_1 = 48\left(1+\pi^{-1}\right)C_{\circ}.
\eqno(2.19)
$$

Also, in {[}V, (4.25), (4.27){]}
the following approximation for the operator $A(t){F}(t)$ was found:
$$
\left\Vert A(t) {F}(t) - t^{2p}SP\right\Vert _{\mathfrak{H}\to \mathfrak{H}}\leqslant  C_2 \left|t\right|^{2p+1},\quad\left|t\right|\leqslant t^{0},
\eqno(2.20)
$$
where
$$
C_2  =  c(p) (\widetilde{B}^{2p} + C_\circ^{2p+1}),
\eqno(2.21)
$$
$C_\circ$, $\widetilde{B}$ are defined by (2.14), (2.15), and $c(p)$ depends only on $p$.

\subsection{Approximation of the resolvent $(A(t)+\eps^{2p}I)^{-1}$}

In this subsection, we prove theorem about approximation of the resolvent $(A(t)+\eps^{2p}I)^{-1}$.
For this, one more condition is needed.

\smallskip
\noindent
\textbf{Condition 2.2.}
\textit{The eigenvalues $\lambda_{j}(t)$ of $A(t)$ satisfy}
$$
\lambda_{j}(t)\geqslant c_{*}t^{2p},\quad j=1,\dots,n,\quad c_{*}>0,\quad \left|t\right|\leqslant t^{0}.
$$

\smallskip
From Condition~2.2 and relations (1.17), (1.20) it follows that
$$
S\geqslant c_{*}I_{\mathfrak{N}},
\eqno(2.22)
$$
and so the spectral germ $S$ is nondegenerate.

The following statement was obtained in [V, Proposition~4.9]; for completeness, we provide the proof.

\smallskip
\noindent
\textbf{Proposition 2.3.}
\textit{For $\varepsilon>0$ and $|t|\leqslant t^{0}$ we have}
 $$
\varepsilon^{2p-1}\left\Vert \left(A(t)+\varepsilon^{2p}I\right)^{-1} F(t) -\left(t^{2p}SP+\varepsilon^{2p}I\right)^{-1}P\right\Vert _{\mathfrak{H}\to \mathfrak{H}} \leqslant C_3.
\eqno(2.23)
$$
\textit{The constant
$
C_3= c_{*}^{-\frac{1}{2p}}\left(2C_1 +c_{*}^{-1} C_2\right)
$
depends only on $p$, $\delta$, the constant $\widetilde{C}$ from} (1.1), \textit{the norm $\|X_p\|$, and $c_*$}.

\smallskip
\noindent
\textbf{Proof.}
We rely on the identity
$$
\begin{aligned}
&G(t,\eps):=\left(A(t)+\varepsilon^{2p}I\right)^{-1} F(t) -\left(t^{2p}SP+\varepsilon^{2p}I\right)^{-1}P
\\
&=
\left(A(t)+\varepsilon^{2p}I\right)^{-1} F(t)(F(t)-P)
+ (F(t)-P) \left(t^{2p}SP+\varepsilon^{2p}I\right)^{-1}P
\\
&- F(t) \left(A(t)+\varepsilon^{2p}I\right)^{-1} (A(t)F(t) - t^{2p}SP)\left(t^{2p}SP+\varepsilon^{2p}I\right)^{-1}P.
\end{aligned}
\eqno(2.24)
$$
By Condition 2.2,
\begin{align}
\| \left(A(t)+\varepsilon^{2p}I\right)^{-1}\|_{\H \to \H}
&\leqslant (c_* t^{2p} + \eps^{2p})^{-1},\quad |t|\leqslant t^0,
\tag{2.25}
\\
\| \left(t^{2p}SP+\varepsilon^{2p}I\right)^{-1}P\|_{\H \to \H}
&\leqslant (c_* t^{2p} + \eps^{2p})^{-1}.
\tag{2.26}
\end{align}
From (2.19), (2.20), and (2.24)--(2.26) it follows that
$$
\| G(t,\eps)\|_{\H \to \H}
\leqslant 2C_1 |t|(c_* t^{2p} + \eps^{2p})^{-1} + C_2 |t|^{2p+1} (c_* t^{2p} + \eps^{2p})^{-2},\quad |t|\leqslant t^0,
$$
which implies  (2.23). $\square$

\smallskip
By (1.7),  for $\eps>0$ and $|t|\leqslant t^0$ we have
$$
\begin{aligned}
&\eps^{2p-1} \| \left(A(t)+\varepsilon^{2p}I\right)^{-1} F(t)^\perp\|_{\H \to \H}
\\
&\leqslant \eps^{2p-1} \| \left(A(t)+\varepsilon^{2p}I\right)^{-1+1/2p}\|
\| \left(A(t)+\varepsilon^{2p}I\right)^{-1/2p} F(t)^\perp\|
\leqslant (3\delta)^{-1/2p}.
\end{aligned}
$$
Combining this with Proposition 2.3, we arrive at the following result.

\smallskip
\noindent
\textbf{Theorem 2.4.}
\textit{Let $A(t)$ be the operator family} (1.3) \textit{satisfying the assumptions of Subsection} 1.1 \textit{and Condition}~2.2.
\textit{Let $P$ be the orthogonal projection of $\H$ onto the subspace $\NN$, and let $S$ be the spectral germ of the family
 $A(t)$ at $t=0$. Then}
 $$
 \begin{aligned}
\varepsilon^{2p-1}\left\Vert \left(A(t)+\varepsilon^{2p}I\right)^{-1}-\left(t^{2p}SP+\varepsilon^{2p}I\right)^{-1}P\right\Vert _{\mathfrak{H}\to \mathfrak{H}}  \leqslant C_{A},
\\
\varepsilon>0,\quad\left|t\right|\leqslant t^{0}.
\end{aligned}
\eqno(2.27)
$$
\textit{Here $\delta$ is chosen in Subsection} 1.1, \textit{and $t^0$ is defined by} (1.4).
\textit{The constant $C_{A}$ is given by
$$
C_{A}= C_3+ \left(3\delta\right)^{-1/2p}
= c_{*}^{-\frac{1}{2p}}\left(2C_1 +c_{*}^{-1} C_2\right)
+ \left(3\delta\right)^{-1/2p}
\eqno(2.28)
$$
and depends only on $p$, $\delta$, the constant $\widetilde{C}$ from} (1.1), \textit{the norm $\|X_p\|$, and $c_*$}.

\smallskip
\noindent
\textbf{Remark 2.5.} It is possible to write down a combersome explicit expression for the constant $C_A$ using
relations (1.5), (2.14), (2.15), (2.19), (2.21), and (2.28). For further application to differential operators,
it is important to know how this constant depends on the problem data.
After possible overstating, $C_A$ can be considered as a polynomial of the variables
$\widetilde{C}$, $\|X_p\|$, $\delta^{-1/2p}$, and $c_*^{-1/2p}$ with positive coefficients depending only on $p$.

\section{The abstract scheme: approximation of the resolvent $\left(A(t)+\varepsilon^{2p}I\right)^{-1}$ with the corrector
taken into account}

In the present section, we obtain approximation of the  resolvent $\left(A(t)+\varepsilon^{2p}I\right)^{-1}$
with the corrector taken into account. Our goal is to prove the following theorem.

\smallskip
\noindent
\textbf{Theorem 3.1.}
\textit{Let $A(t)$ be the operator family} (1.3) \textit{satisfying the assumptions of Subsection} 1.1 \textit{and Condition}~2.2.
\textit{Let $P$ be the orthogonal projection of $\H$ onto the subspace $\NN$. Suppose that $Z$ is the operator} (1.10) \textit{and
$S$ is the spectral germ of $A(t)$ at the point $t=0$. Then}
$$
\begin{aligned}
\varepsilon^{p-1} &\left\Vert
A(t)^{1/2}\left(\left(A(t)+\varepsilon^{2p}I\right)^{-1}-\left(I+t^{p}Z\right)\left(t^{2p}SP+\varepsilon^{2p}I\right)^{-1}P\right)\right\Vert _{\mathfrak{H}\to \mathfrak{H}}
\\
&\leqslant \check{C}_{A}, \quad \eps>0, \quad |t|\leqslant t^0.
\end{aligned}
\eqno(3.1)
$$
\textit{Here $t^0$ is defined by} (1.4). \textit{The constant $\check{C}_{A}$ depends only on
$p$, $\delta$, the constant $\widetilde{C}$ from} (1.1), \textit{the norm $\|X_p\|$, and $c_*$}.

\subsection{Proof of Theorem 3.1: step 1}

Denote
\begin{align}
\mathfrak{A}_{\varepsilon}(t) & :=A(t)^{1/2}\left(A(t)+\varepsilon^{2p}I\right)^{-1},
\tag{3.2}
\\
\Xi\left(t,\varepsilon\right) & :=\left(t^{2p}SP+\varepsilon^{2p}I\right)^{-1}P.
\tag{3.3}
\end{align}
We need to estimate the operator
$$
\Upsilon(t,\eps):= \mathfrak{A}_{\varepsilon}(t) - A(t)^{1/2}\left(I+t^{p}Z\right)\Xi(t,\varepsilon).
\eqno(3.4)
$$
Note that the operator (3.3) satisfies (2.26).

The operator (3.4) can be represented as the sum of four terms
$$
\Upsilon(t,\eps) = J_1(t,\eps)+J_2(t,\eps)+ J_3(t,\eps) +J_4(t,\eps),
\eqno(3.5)
$$
where
\begin{align}
J_1(t,\eps)&:= \mathfrak{A}_{\varepsilon}(t) F(t)^\perp,
\tag{3.6}
\\
J_2(t,\eps)&:= \mathfrak{A}_{\varepsilon}(t) F(t) (F(t)-P),
\tag{3.7}
\\
J_3(t,\eps)&:= F(t)\mathfrak{A}_{\varepsilon}(t) P - F(t) A(t)^{1/2} \Xi(t,\eps),
\tag{3.8}
\\
J_4(t,\eps)&:=  {A}(t)^{1/2}  (F(t)-P) \Xi(t,\eps)  - t^p {A}(t)^{1/2}  Z \Xi(t,\eps).
\tag{3.9}
\end{align}

To estimate the operator (3.6), we apply the Young inequality in the form
$$
\left(\lambda +\varepsilon^{2p}\right)^{-1}\leqslant \lambda^{-1/2-1/2p} \varepsilon^{1-p},
\quad \lambda>0,\quad \eps>0.
\eqno(3.10)
$$
By (1.7), (3.2), (3.6), and (3.10),
$$
 \left\Vert J_1 (t,\eps)\right\Vert _{\mathfrak{H}\to \mathfrak{H}} \leqslant \sup_{\lambda\geqslant3\delta} \lambda^{1/2} \left(\lambda+\varepsilon^{2p}\right)^{-1}
\leqslant \left(3\delta\right)^{-1/2p} \varepsilon^{1-p},\quad |t|\leqslant t^{0}.
\eqno(3.11)
$$
Condition 2.2 and (3.10) imply that
$$
\begin{aligned}
&\left\Vert \mathfrak{A}_{\varepsilon}(t) {F}(t) \right\Vert _{\mathfrak{H}\to \mathfrak{H}}
=\sup_{1\leqslant l\leqslant n}\sqrt{\lambda_{l}(t)}\left(\lambda_{l}(t)+\varepsilon^{2p}\right)^{-1}
\\
&\leqslant
\eps^{1-p} \sup_{1\leqslant l\leqslant n} \lambda_{l}(t)^{-1/2p}
\leqslant
c_*^{-1/2p} |t|^{-1}\eps^{1-p},
\quad 0< \left|t\right|\leqslant t^{0}.
\end{aligned}
\eqno(3.12)
$$
Together with (2.19) this yields the following estimate for the operator (3.7):
$$
\|J_2(t,\eps)\|_{\mathfrak{H}\to \mathfrak{H}} \leqslant
C_1 c_*^{-1/2p} \eps^{1-p}, \quad |t| \leqslant t^0.
\eqno(3.13)
$$

Next, to handle the operator (3.8), we apply the following analog of the resolvent identity:
$$
\begin{aligned}
 & {F}(t) \left(A\left(t\right)+\varepsilon^{2p}I\right)^{-1}P -  {F}(t)\Xi(t,\varepsilon)
 \\
&= -{F}(t)\left(A(t)+\varepsilon^{2p}I\right)^{-1}\left(A(t) {F}(t)-t^{2p}SP\right)\Xi(t,\varepsilon),
\end{aligned}
$$
and multiply both sides by $A(t)^{1/2}$. Then the operator (3.8) can be written as
$$
J_3(t,\eps) = -{F}(t) \mathfrak{A}_\eps(t)  \left(A(t) {F}(t)-t^{2p}SP\right)\Xi(t,\varepsilon).
$$
Combining this with (2.20), (2.26), and (3.12), we obtain
$$
\begin{aligned}
\|J_3(t,\eps)\|_{\mathfrak{H}\to \mathfrak{H}} &\leqslant
C_2 c_*^{-1/2p} (c_* t^{2p} + \eps^{2p})^{-1} |t|^{2p+1} |t|^{-1} \eps^{1-p}
\\
&\leqslant C_2 c_*^{-1/2p-1} \eps^{1-p}, \quad |t| \leqslant t^0.
\end{aligned}
\eqno(3.14)
$$

As a result, relations (3.5), (3.11), (3.13), and (3.14) imply that
$$
\|\Upsilon(t,\eps)\|_{\mathfrak{H}\to \mathfrak{H}} \leqslant
C_4 \eps^{1-p}+
\|J_4(t,\eps)\|_{\mathfrak{H}\to \mathfrak{H}},\quad |t| \leqslant t^0,
\eqno(3.15)
$$
where
$$
C_4 = (3\delta)^{-1/2p} + C_1 c_*^{-1/2p} + C_2 c_*^{-1/2p-1}.
\eqno(3.16)
$$
Thus, the proof of estimate (3.1) is reduced to estimation of the operator~(3.9).

\subsection{The iteration procedure}

Now, we rewrite the resolvent identity (2.4) as
$$
R_{z}(t)=R_{z}(0)-\Omega_{z}(0)T_{\delta}(t)R_{z}(t).
\eqno(3.17)
$$
We will iterate it using representation (2.11) for $T_\delta(t)$.
After $p$ iterations we obtain
$$
R_{z}(t)-R_{z}(0)=t\Psi_{1}(z)+\dots+t^{p}\Psi_{p}(z)+\Psi_{*}(t,z).
\eqno(3.18)
$$
Together with (2.18) this implies
$$
{F}(t)- P= t {F}_{1}+\dots+t^{p} {F}_{p}+{F}_{*}(t),
\eqno(3.19)
$$
where
\begin{align}
{F}_{i} & =-\frac{1}{2\pi i}\oint_{\Gamma_{\delta}}\Psi_{i}(z)\,dz,\quad i=1,\dots,p,
\tag{3.20}
\\
{F}_{*}\left(t\right) & =-\frac{1}{2\pi i}\oint_{\Gamma_{\delta}}\Psi_{*}(t,z)\,dz.
\tag{3.21}
\end{align}
By (3.19), the operator (3.9) can be represented as
$$
J_4 (t,\eps)= J_4^{(1)} (t,\eps)+ J_4^{(2)} (t,\eps) + J_4^{(3)} (t,\eps),
\eqno(3.22)
$$
where
\begin{align}
J_4^{(1)} (t,\eps)&:= \sum_{i=1}^{p-1} t^i A(t)^{1/2} F_i \Xi(t,\eps),
\tag{3.23}
\\
J_4^{(2)} (t,\eps)&:=  t^p A(t)^{1/2}( F_p -Z) \Xi(t,\eps),
\tag{3.24}
\\
J_4^{(3)} (t,\eps)&:= A(t)^{1/2} F_*(t) \Xi(t,\eps).
\tag{3.25}
\end{align}

\subsection{Estimates for the operators $F_i$}

Now we find expressions for the operators ${F}_{i}$ in terms of the coefficients of the expansions (1.18)
for the eigenvectors $\varphi_{j}(t)$. By (1.23),
\begin{align*}
 & {F}(t)=\sum_{j=1}^{n}\left(\cdot,\varphi_{j}(t)\right)_{\mathfrak{H}}\varphi_{j}(t)
\\
 &=\sum_{j=1}^{n}\left(\cdot, \omega_{j}+t\varphi_{j}^{(1)}+\dots+t^{p}\varphi_{j}^{(p)}\right)_{\mathfrak{H}}\left(\omega_{j}
+t\varphi_{j}^{(1)}+\dots+t^{p}\varphi_{j}^{(p)}\right)+O\left(t^{p+1}\right)
\\
 & =\sum_{j=1}^{n}\left(\cdot,\omega_{j}\right)_{\mathfrak{H}}\omega_{j}+t\sum_{j=1}^{n}\left\{ \left(\cdot,\omega_{j}\right)_{\mathfrak{H}}\varphi_{j}^{(1)}+\left(\cdot,\varphi_{j}^{(1)}\right)_{\mathfrak{H}}\omega_{j}\right\} +\dots
\\
 & +t^{p}\sum_{j=1}^{n}\left\{ \left(\cdot,\omega_{j}\right)_{\mathfrak{H}}\varphi_{j}^{(p)}+
\left(\cdot,\varphi^{(1)}_{j}\right)_{\mathfrak{H}}\varphi_{j}^{(p-1)}+ \dots+\left(\cdot,\varphi_{j}^{(p)}\right)_{\mathfrak{H}}\omega_{j}\right\} +{F}_{*}\left(t\right).
\end{align*}
According to (3.19),
$$
{F}_{i}  = \sum_{j=1}^{n}\sum_{k=0}^{i}\left(\cdot,\varphi_{j}^{(k)}\right)\varphi_{j}^{(i-k)},\quad i=0,\dots,p,
\eqno(3.26)
$$
here, for convenience, it is assumed that $\varphi_{j}^{(0)}=\omega_{j}$, $j=1,\dots,n$, and $F_0=P$.

By (1.21), we have $\varphi_{j}^{\left(l\right)}\in\mathfrak{N}$,
$l=0,\dots,p-1$, $j=1,\dots,n$. Hence, the operators
${F}_{i}$ with $i=1,\dots,p-1$ take $\mathfrak{N}$ to $\mathfrak{N}$ and $\mathfrak{N}^\perp$ to $\{0\}$.
This fact will help us to estimate the norm of the operator (3.23) by  $O(\varepsilon^{-p+1})$.

To estimate the operators $F_i$, we use the invariant representations (3.20) and estimate the integrands
 uniformly in $z$. For this, we first find the invariant representations
for the operators $\Psi_i(z)$. We iterate the identity (3.17) using (2.11).
The symbol  ``$\sim$'' will be used instead of  ``$=$'' if the terms of order $t^k$ with $k> p$ are dropped. We have:
\begin{align*}
R_{z}&(t)  =R_{z}(0)-\Omega_{z}(0)T_{\delta}(t)\left(R_{z}(0)-\Omega_{z}(0)T_{\delta}(t)R_{z}(t)\right)
\\
 & \sim R_{z}(0)-\Omega_{z}(0)\sum_{i_{1}=1}^{p}t^{i_{1}}T_{\delta}^{\left(i_{1}\right)}R_{z}(0)
+\left(\Omega_{z}(0) T_{\delta}(t)\right)^{2}R_{z}(t)
\\
 & \sim R_{z}(0)-\Omega_{z}(0)\sum_{i_{1}=1}^{p}t^{i_{1}}T_{\delta}^{\left(i_{1}\right)}R_{z}(0)
\\
 & +\Omega_{z}(0) \sum_{i_{1}=1}^{p-1}t^{i_{1}}T_{\delta}^{\left(i_{1}\right)}\Omega_{z}(0) \sum_{i_{2}=1}^{p-1}t^{i_{2}}T_{\delta}^{\left(i_{2}\right)}R_{z}(0)
 -\left(\Omega_{z}(0) T_{\delta}(t)\right)^{3}R_{z}(t)
\\
 & \sim R_{z}(0) -\Omega_{z}(0) \sum_{i_{1}=1}^{p}t^{i_{1}}T_{\delta}^{\left(i_{1}\right)}R_{z}(0)
\\
 &+\Omega_{z}(0) \sum_{i_{1}=1}^{p-1}t^{i_{1}}T_{\delta}^{\left(i_{1}\right)}\Omega_{z}(0) \sum_{i_{2}=1}^{p-1}t^{i_{2}}T_{\delta}^{\left(i_{2}\right)}R_{z}(0)
\\
 & -\Omega_{z}(0) \sum_{i_{1}=1}^{p-2}t^{i_{1}}T_{\delta}^{\left(i_{1}\right)}\Omega_{z}(0) \sum_{i_{2}=1}^{p-2}t^{i_{2}}T_{\delta}^{\left(i_{2}\right)}\Omega_{z}(0) \sum_{i_{3}=1}^{p-2}t^{i_{3}}T_{\delta}^{\left(i_{3}\right)}R_{z}(0)
\\
 & +
\left(\Omega_{z}(0) T_{\delta}(t)\right)^{4}R_{z}(t).
\end{align*}
We continue this iteration procedure until the last term becomes
$\left(\Omega_{z}(0)T_{\delta}(t)\right)^{p+1}R_{z}(t).$
The final expression takes the form
$$
\begin{aligned}
R_{z}(t) & \sim R_{z}(0)-\Omega_{z}\left(0\right)\sum_{i_{1}=1}^{p}t^{i_{1}}T_{\delta}^{\left(i_{1}\right)}R_{z}(0)
\\
 & +\Omega_{z}(0) \sum_{i_{1}=1}^{p-1}t^{i_{1}}T_{\delta}^{\left(i_{1}\right)}\Omega_{z}(0) \sum_{i_{2}=1}^{p-1}t^{i_{2}}T_{\delta}^{\left(i_{2}\right)}R_{z}(0) +\dots
\\
 & +\left(-1\right)^{k}\Omega_{z}(0)\sum_{i_{1}=1}^{p+1-k}t^{i_{1}}T_{\delta}^{\left(i_{1}\right)}\cdots\Omega_{z}(0) \sum_{i_{k}=1}^{p+1-k}t^{i_{k}}T_{\delta}^{\left(i_{k}\right)}R_{z}(0)+\dots
\\
 & +\left(-1\right)^{p}t^{p}\left(\Omega_{z}(0)T_{\delta}^{\left(1\right)}\right)^{p}R_{z}(0).
\end{aligned}
$$
Let us extract $\Psi_{i}(z)$. For simplicity, we use the following notation.
Let $\gamma^{k}=\left(\gamma_{1}^{k},\dots,\gamma_{k}^{k}\right)$ be a multiindex of length  $k$
such that $\gamma_{i}^{k}\geqslant1$, $i=1,\dots,k$. Denote
$$
\left(\Omega_{z}(0) T_{\delta}^{(\cdot)}\right)^{\gamma^{k}}=\Omega_{z}(0) T_{\delta}^{\left(\gamma_{1}^{k}\right)}\cdots\Omega_{z}(0)T_{\delta}^{\left(\gamma_{k}^{k}\right)}.
$$
Then $\Psi_{i}(z)$ is given by
$$
\Psi_{i}(z)=\sum_{k=1}^{i}\left(-1\right)^{k}\sum_{\left|\gamma^{k}\right|=i}\left(\Omega_{z}(0) T_{\delta}^{(\cdot)}\right)^{\gamma^{k}}R_{z}(0),
\quad i=1,\dots,p.
\eqno(3.27)
$$
Now, relations  (2.5), (2.9), (2.10), and (2.13) imply the following estimates for the operators (3.27):
\begin{align*}
 & \left\Vert \Psi_{i}(z)\right\Vert _{\mathfrak{H}\to \mathfrak{H}}
\leqslant\delta^{-1/2}\left\Vert \Psi_{i}(z) \right\Vert _{\mathfrak{H}\to \mathfrak{d}}\leqslant
\\
 & \leqslant\delta^{-1/2}\sum_{k=1}^{i}\sum_{\left|\gamma^{k}\right|=i}\left\Vert \left(\Omega_{z}(0) T_{\delta}^{(\cdot)}\right)^{\gamma^{k}}R_{z}(0) \right\Vert _{\mathfrak{H}\to \mathfrak{d}} \leqslant
\\
 & \leqslant4\delta^{-1}\sum_{k=1}^{i}\sum_{\left|\gamma^{k}\right|=i}\left(13\widetilde{B}\right)^{k},\quad i=1,\dots,p.
\end{align*}
Since $\widetilde{B}\geqslant1$, then $\left(13\widetilde{B}\right)^{k}$
does not exceed $\left(13\widetilde{B}\right)^{i}$ if $k=1,\dots,i$, whence
$$
\left\Vert \Psi_{i}(z)\right\Vert _{\mathfrak{H}\to \mathfrak{H}} \leqslant4\delta^{-1}\left(13\widetilde{B}\right)^{i}\biggl(\sum_{k=1}^{i}\sum_{\left|\gamma^{k}\right|=i}1\biggr), \quad i=1,\dots,p.
$$
Using (3.20) and noting that the length of the contour $\Gamma_{\delta}$ is equal to $2\pi\delta+2\delta$, we obtain
$$
\left\Vert {F}_{i}\right\Vert _{\mathfrak{H}\to \mathfrak{H}}
\leqslant4\left(13\widetilde{B}\right)^{i}\biggl(\sum_{k=1}^{i}\sum_{\left|\gamma^{k}\right|=i}1\biggr)\left(\pi^{-1}+1\right):=C^{\left(i\right)},
\quad i=1,\dots,p.
\eqno(3.28)
$$

\subsection{Estimate for the operator $J_4^{(1)}(t,\eps)$.}

Let $u\in\mathfrak{H}$, and let $v=\Xi(t,\varepsilon)u\in\mathfrak{N}$.
We estimate the norms
$$
\left\Vert A(t)^{1/2} {F}_{i}\Xi\left(t,\varepsilon\right)u\right\Vert _{\mathfrak{H}}=
\left\Vert X(t) {F}_{i}v\right\Vert _{\mathfrak{H}_{*}},\quad i=1,\dots,p-1.
\eqno(3.29)
$$
As has already been mentioned, the operators ${F}_{i}$, $i=1,\dots,p-1$, take $\mathfrak{N}$ to $\mathfrak{N}$.
Taking (1.2) into account, we simplify the expression under the norm sign
in the right-hand side of (3.29):
\begin{align*}
 & X(t){F}_{i}v=\left(X_{0}+tX_1 + \dots+t^{p-1}X_{p-1}+t^{p}X_{p}\right){F}_{i}v=\\
 & =t^{p}X_{p}{F}_{i}v,\quad v\in\mathfrak{N},\quad i=1,\dots,p-1.
\end{align*}
Consequently, by (2.26) and (3.28),
$$
\left\Vert X(t) {F}_{i}v\right\Vert _{\mathfrak{H}_{*}}
\leqslant C^{(i)} \left|t\right|^{p}\left\Vert X_{p}\right\Vert \left\Vert v\right\Vert _{\mathfrak{H}}
\leqslant C^{(i)} \left|t\right|^{p}\left\Vert X_{p}\right\Vert (c_* t^{2p} + \eps^{2p})^{-1}
\left\Vert u\right\Vert _{\mathfrak{H}}.
\eqno(3.30)
$$
By the Young inequality (3.10),
$$
(c_* t^{2p} + \eps^{2p})^{-1} \leqslant c_*^{-1/2-1/2p} |t|^{-1-p} \eps^{1-p}.
\eqno(3.31)
$$
Combining this with  (3.29), (3.30) and noting that $t^0\leqslant 1$, we obtain the following estimate for the operator (3.23):
$$
\| J_4^{(1)}(t,\eps)\|_{\mathfrak{H}\to \mathfrak{H}}
\leqslant C_5 \eps^{1-p}, \quad |t|\leqslant t^0.
\eqno(3.32)
$$
where
$$
C_5 = c_*^{-1/2-1/2p} \|X_p\| \biggl(\sum_{i=1}^{p-1} C^{(i)} \biggr).
\eqno(3.33)
$$

\subsection{Estimate for the operator $J_4^{(2)}(t,\eps)$}

Now, we consider the operator $F_p$. According to (3.26),
$$
{F}_{p}=\sum_{j=1}^{n}\sum_{i=0}^{p}\left(\cdot,\varphi_{j}^{\left(i\right)}\right)\varphi_{j}^{\left(p-i\right)}.
$$
Denote
$
\widetilde{\omega}_{j}:=\varphi_{j}^{\left(p\right)}-Z\omega_{j}.
$
By (1.10) and (1.22), $\widetilde{\omega}_{j}\in\mathfrak{N}$.
The operator ${F}_{p}$ can be represented as
$$
{F}_{p}=\check{{F}}_{p}+\widetilde{F}_{p},
\eqno(3.34)
$$
where
\begin{align}
\check{F}_{p}&=\sum_{j=1}^{n}\left( \left(\cdot,\omega_{j}\right)Z\omega_{j}+\left(\cdot,Z\omega_{j}\right)\omega_{j}\right),
\tag{3.35}
\\
\widetilde{F}_{p}&=\sum_{j=1}^{n}\sum_{i=1}^{p-1}\left(\cdot,\varphi_{j}^{\left(i\right)}\right)\varphi_{j}^{\left(p-i\right)}+\sum_{j=1}^{n}\left(\left(\cdot,\omega_{j}\right)\widetilde{\omega}_{j}+\left(\cdot,\widetilde{\omega}_{j}\right)\omega_{j}\right).
\tag{3.36}
\end{align}
Relation $\widetilde{\omega}_{j}\in\mathfrak{N}$ together with  (1.21) show that the operator
(3.36) takes $\mathfrak N$ to $\mathfrak N$ and ${\mathfrak N}^\perp$ to $\{0\}$.

Using (3.35) and (1.19), we see that
$
\check{F}_{p}=ZP+PZ^{*}.
$
From the definition (1.10) of $Z$ it follows that
$PZ=0$, whence $Z^*P=0$. Consequently, $(\check{F}_{p}-Z)P=0$.
Together with (3.34) this implies
$$
({F}_{p} - Z)P= \widetilde{F}_{p}P.
\eqno(3.37)
$$
Thus, the operator (3.24) can be written as
$$
J_4^{(2)}(t,\eps) = t^p A(t)^{1/2} \widetilde{F}_{p} \Xi(t,\eps).
\eqno(3.38)
$$

By (1.11), (3.28), and (3.37),
$$
\| \widetilde{F}_{p}P \|_{\mathfrak{H}\to \mathfrak{H}}
\leqslant  (1/6) \delta^{-1/2} \|X_p\| + C^{(p)}.
\eqno(3.39)
$$
Again, let $u\in\mathfrak{H}$ and let $v=\Xi\left(t,\varepsilon\right)u\in\mathfrak{N}$.
The fact that $\widetilde{F}_p$ takes $\mathfrak{N}$ to $\mathfrak{N}$ together with (1.2)
allow us to simplify the expression for $X(t)\widetilde{F}_{p}v$:
$$
 X(t)\widetilde{F}_{p}v=\left(X_{0}+tX_1+ \dots+t^{p-1}X_{p-1}+t^{p}X_{p}\right) \widetilde{F}_{p}v=
 t^{p}X_{p}\widetilde{F}_{p}v.
$$
Hence, using  (2.26), (3.38), and  (3.39), we obtain
\begin{align*}
 & \left\Vert J_4^{(2)}(t,\eps) u\right\Vert _{\mathfrak{H}}
 =\left|t\right|^{p} \left\Vert X(t) \widetilde{F}_{p}v \right\Vert _{\mathfrak{H}_{*}} = t^{2p}
\left\Vert X_{p}\widetilde{F}_{p}v\right\Vert _{\mathfrak{H}_{*}}
\\
 & \leqslant t^{2p} \left\Vert X_{p}\right\Vert \left(C^{\left(p\right)}+(1/6) \delta^{-1/2} \left\Vert X_{p}\right\Vert \right)
 (c_* t^{2p} + \eps^{2p})^{-1}
\left\Vert u\right\Vert _{\mathfrak{H}}.
\end{align*}
Combining this with (3.31) and noting that $t^0\leqslant 1$, we arrive at
$$
\left\Vert J_4^{(2)}(t,\eps) u\right\Vert _{\mathfrak{H}}
\leqslant C_6 \eps^{1-p},\quad |t| \leqslant t^0,
\eqno(3.40)
$$
where
$$
C_6 = \left\Vert X_{p}\right\Vert \left(C^{\left(p\right)}+(1/6) \delta^{-1/2} \left\Vert X_{p}\right\Vert \right) c_*^{-1/2-1/2p}.
\eqno(3.41)
$$

\subsection{Estimate for the operator $J_4^{(3)}(t,\eps)$}

To estimate the operator (3.25), we use representation (3.21).
First, we estimate the operator $A(t)^{1/2} \Psi_*(t,z)$ uniformly with respect to $z\in \Gamma_\delta$.
For this, we need some auxiliary statements.

\smallskip
\noindent
\textbf{Lemma 3.2.} \textit{We have}
$$
\left\Vert A(t)^{1/2}\right\Vert _{\mathfrak{d} \to \mathfrak{H}}\leqslant \sqrt{6},\quad\left|t\right|\leqslant t^{0}.
$$

\smallskip\noindent
\textbf{Proof.}
Let $u\in\mathfrak{d}$. Using the definition of the operator $A(t)^{1/2}$, we have
$\left\Vert A(t)^{1/2}u\right\Vert _{\mathfrak{H}} = \left\Vert X(t)u\right\Vert _{\mathfrak{H}_{*}}$.
By (2.16),
$$
\left\Vert X (t)u\right\Vert _{\mathfrak{H}_{*}}^{2} \leqslant
\left\Vert X_{0}u\right\Vert _{\mathfrak{H}_{*}}^{2} + \left(\left\Vert X_{0}u\right\Vert _{\mathfrak{H}_{*}}^{2}+\delta\left\Vert u\right\Vert _{\mathfrak{H}}^{2}\right) C_\circ \left|t\right|, \quad |t| \leqslant t^0.
$$
Recalling that $\|u\|^2_{\mathfrak{d}} = \left\Vert X_{0}u\right\Vert _{\mathfrak{H}_{*}}^{2}+\delta\left\Vert u\right\Vert _{\mathfrak{H}}^{2}$
and taking (2.14) into accout, we arrive at
$$
\left\Vert A (t)^{1/2}u\right\Vert _{\mathfrak{H}}^{2} \leqslant (1+ C_\circ t^0) \|u\|^2_{\mathfrak{d}}= 6 \|u\|^2_{\mathfrak{d}},
\quad u\in \mathfrak{d},\quad |t|\leqslant t^0.
$$
$\square$

\smallskip
\noindent
\textbf{Lemma 3.3.} \textit{We have}
$$
\left\Vert R_z(t)\right\Vert _{\mathfrak{H} \to \mathfrak{d}}\leqslant (2\sqrt{3}+3)\delta^{-1/2},\quad z \in \Gamma_\delta,\quad \left|t\right|\leqslant t^{0}.
\eqno(3.42)
$$

\smallskip\noindent
\textbf{Proof.}
Let $u\in\mathfrak{H}$. By the definition of the norm in $\mathfrak{d}$, we have
$$
\left\Vert R_{z}(t)u\right\Vert _{\mathfrak{d}} \leqslant  \delta^{1/2} \left\Vert R_{z}(t)u\right\Vert _{\mathfrak{H}} +\left\Vert
X_{0}R_{z}(t)u\right\Vert _{\mathfrak{H}_{*}}.
$$
Together with (1.6) this implies
$$
\begin{aligned}
 & \left\Vert R_{z}(t)u\right\Vert _{\mathfrak{d}}\leqslant2\left\Vert X(t)R_{z}(t)u\right\Vert _{\mathfrak{H}_{*}}+3\delta^{1/2}\left\Vert R_{z}(t) u\right\Vert _{\mathfrak{H}}\leqslant
\\
 & \leqslant2\left\Vert X(t)R_{z}(t)u\right\Vert _{\mathfrak{H}_{*}}+3\delta^{-1/2}\left\Vert u\right\Vert _{\mathfrak{H}},\quad\left|t\right|\leqslant t^{0},\quad z\in\Gamma_{\delta}.
\end{aligned}
\eqno(3.43)
$$
Since $\left|z\right|\leqslant2\delta$ for $z\in\Gamma_{\delta}$, then
$$
\begin{aligned}
 & \left\Vert X(t)R_{z}(t)u\right\Vert _{\mathfrak{H}_{*}}^{2}
=\left(A(t)R_{z}(t)u,R_{z}(t)u\right)_{\mathfrak{H}}
\\
 & =\left(u,R_{z}(t)u\right)_{\mathfrak{H}}+z\left\Vert R_{z}(t)u\right\Vert _{\mathfrak{H}}^{2}
\\
 & \leqslant\left\Vert R_{z}(t)\right\Vert _{\mathfrak{H}\to \mathfrak{H}}
\left\Vert u\right\Vert _{\mathfrak{H}}^{2}+2\delta\left\Vert R_{z}(t)\right\Vert _{\mathfrak{H}\to \mathfrak{H}}^{2}
\left\Vert u\right\Vert _{\mathfrak{H}}^{2}\leqslant 3\delta^{-1}\left\Vert u\right\Vert _{\mathfrak{H}}^{2}.
\end{aligned}
\eqno(3.44)
$$
Combining (3.43) and (3.44), we arrive at the required inequality (3.42). $\square$

\smallskip
Now, we estimate the operator $\Psi_{*}\left(t,z\right)$ in the $(\mathfrak{H}\to \mathfrak{d})$-norm.
For this, we apply the iteration procedure once more in order to extract the term $\Psi_{*}\left(t,z\right)$ in the expansion (3.18).
As before, we iterate (3.17) using (2.11). Our goal is to estimate the norms of all
operators at $t^k$ with  $k> p$. Therefore, we will use the sign $\sim$ if some terms of order $t^k$ with  $k\leqslant p$ are dropped.
The first iteration is
\begin{align*}
 & R_{z}(t) =R_{z}(0) -\Omega_{z}(0) T_{\delta}(t)\left(R_{z}(0)-\Omega_{z}(0) T_{\delta}(t)R_{z}(t)\right)
\\
 & =R_{z}(0)-\Omega_{z}(0)\left(tT_{\delta}^{\left(1\right)}+\dots+t^{p}T_{\delta}^{\left(p\right)}\right)R_{z}(0)
\\
 & -\Omega_{z}(0) \left(t^{p+1}T_{\delta}^{\left(p+1\right)}+\dots+t^{2p}T_{\delta}^{\left(2p\right)}\right)R_{z}(0) +\left(\Omega_{z}(0) T_{\delta}(t)\right)^{2} R_{z}(t).
\end{align*}
Hence,
$$
R_{z}(t) \sim -\Omega_{z}(0) \left(t^{p+1}T_{\delta}^{\left(p+1\right)}+\dots+t^{2p}T_{\delta}^{\left(2p\right)}\right)R_{z}(0) +\left(\Omega_{z}(0) T_{\delta}(t)\right)^{2}R_{z}(t).
\eqno(3.45)
$$
Denote the first term in the right-hand side of (3.45) by ${\mathcal I}_1(t,z)$
and estimate it with the help of (2.9), (2.10), (2.13), and the inequality $t^0 \leqslant 1$:
$$
\left\Vert {\mathcal I}_1(t,z) \right\Vert_{\mathfrak{H} \to \mathfrak{d}}\leqslant C_{(1)} |t |^{p+1},
\quad |t| \leqslant t^0, \quad C_{(1)}= 52 \delta^{-1/2} p \widetilde{B}.
$$
Now we write down the unrecorded terms for the second iteration leaving only the terms
with $t^k$, where $k>p$:
\begin{align*}
 & \left(\Omega_{z}(0) T_{\delta}(t)\right)^{2} R_{z}(t) =\left(\Omega_{z}(0) T_{\delta}(t)\right)^{2} \left(R_{z}(0) -\Omega_{z}(0) T_{\delta}(t) R_{z}(t)\right)
\\
 & =t^{2}\Omega_{z}(0) \left(T_{\delta}^{\left(1\right)}\Omega_{z}(0)T_{\delta}^{\left(1\right)} \right)R_{z}(0)
\\
 & +t^{3}\Omega_{z}(0) \left(T_{\delta}^{\left(1\right)}\Omega_{z}(0)T_{\delta}^{\left(2\right)} +T_{\delta}^{\left(2\right)}\Omega_{z}(0) T_{\delta}^{\left(1\right)}\right)R_{z}(0)+\dots
\\
 & +t^{p}\Omega_{z}(0) \left(T_{\delta}^{\left(1\right)}\Omega_{z}(0) T_{\delta}^{\left(p-1\right)}+\dots+T_{\delta}^{\left(p-1\right)}\Omega_{z}(0) T_{\delta}^{\left(1\right)}\right)R_{z}(0)
\\
 & +t^{p+1}\Omega_{z}\left(0\right)\left(T_{\delta}^{\left(1\right)}\Omega_{z}(0) T_{\delta}^{\left(p\right)}+\dots+T_{\delta}^{\left(p\right)}\Omega_{z}(0) T_{\delta}^{\left(1\right)}\right)R_{z}(0)+\dots
\\
 & +t^{4p}\Omega_{z}(0) \left(T_{\delta}^{\left(2p\right)}\Omega_{z}(0) T_{\delta}^{\left(2p\right)}\right)R_{z}(0) -\left(\Omega_{z}(0) T_{\delta}(t) \right)^{3}R_{z}(t)
\\
 & \sim t^{p+1}\Omega_{z}(0) \left(T_{\delta}^{\left(1\right)}\Omega_{z}(0) T_{\delta}^{\left(p\right)}+\dots+T_{\delta}^{\left(p\right)}\Omega_{z}(0) T_{\delta}^{\left(1\right)}\right)R_{z}(0)+\dots
\\
 & +t^{4p}\Omega_{z}(0) \left(T_{\delta}^{\left(2p\right)}\Omega_{z}(0) T_{\delta}^{\left(2p\right)}\right)R_{z}(0)
-\left(\Omega_{z}(0) T_{\delta}(t)\right)^{3}R_{z}(t)
\\
& =:
{\mathcal I}_2(t,z) - \left(\Omega_{z}(0) T_{\delta}(t)\right)^{3}R_{z}(t).
\end{align*}
Again, we estimate the extracted term ${\mathcal I}_2(t,z)$ with the help of (2.9), (2.10), (2.13), and the inequality $t^0 \leqslant 1$:
$$
\left\Vert {\mathcal I}_2(t,z) \right\Vert_{\mathfrak{H} \to \mathfrak{d}}\leqslant C_{(2)} |t |^{p+1},
\quad |t| \leqslant t^0,
$$
where
$$
C_{(2)}= 4 \cdot 13^2 \delta^{-1/2} c_p^{(2)} \widetilde{B}^2,
\quad c_p^{(2)} = p + (p+1)+\dots + 2p+ (2p-1)+\dots +1.
$$
Now the unrecorded term is  $- \left(\Omega_{z}(0) T_{\delta}(t)\right)^{3}R_{z}(t)$.

We continue this iteration procedure until the unrecorded term becomes
$(-1)^{p+1}\left(\Omega_{z}(0) T_{\delta}(t)\right)^{p+1}R_{z}(t)=: {\mathcal I}^0(t,z)$.
All the terms ${\mathcal I}_j(t,z)$, $j=1,\dots,p,$ extracted in the consequent iterations are estimated as follows:
$$
\left\Vert {\mathcal I}_j(t,z) \right\Vert_{\mathfrak{H} \to \mathfrak{d}}\leqslant C_{(j)} |t |^{p+1},
\quad |t| \leqslant t^0,\quad j=1,\dots,p,
\eqno(3.46)
$$
where
$$
C_{(j)}= 4 \cdot 13^j \delta^{-1/2} c_p^{(j)} \widetilde{B}^j,
\eqno(3.47)
$$
and $c_p^{(j)}$ depends only on $p$ and $j$.
Finally, the term  ${\mathcal I}^0(t,z)$ is estimated with the help of (2.10), (2.12), and Lemma~3.3:
$$
\left\Vert {\mathcal I}^0(t,z) \right\Vert_{\mathfrak{H} \to \mathfrak{d}}
\leqslant C_{(p+1)} |t |^{p+1},
\quad |t| \leqslant t^0,
\eqno(3.48)
$$
where
$$
C_{(p+1)} = (2\sqrt{3}+3)\cdot 13^{p+1} \delta^{-1/2}C_\circ^{p+1}.
\eqno(3.49)
$$
Clearly, $\Psi_*(t,z) = {\mathcal I}_1(t,z)+\dots+ {\mathcal I}_p(t,z) + {\mathcal I}^0(t,z)$.
As a result, relations (3.46) and (3.48) imply that
$$
\left\Vert \Psi_*(t,z) \right\Vert_{\mathfrak{H} \to \mathfrak{d}}
\leqslant C_{7} |t |^{p+1},
\quad z\in \Gamma_\delta,\quad  |t| \leqslant t^0,
\eqno(3.50)
$$
where
$$
C_7 = \sum_{j=1}^{p+1} C_{(j)}.
\eqno(3.51)
$$

Now,  Lemma~3.2 and estimate (3.50) yield
$$
\left\Vert A(t)^{1/2}\Psi_*(t,z) \right\Vert_{\mathfrak{H} \to \mathfrak{H}}
\leqslant \sqrt{6} C_{7} |t |^{p+1},
\quad z\in \Gamma_\delta,\quad  |t| \leqslant t^0.
$$
Combining this with (3.21) and the estimate for the length of the contour $\Gamma_\delta$
by $2\pi\delta+2\delta \leqslant 2\pi+2$, we obtain
$$
\left\Vert A(t)^{1/2} F_*(t) \right\Vert_{\mathfrak{H} \to \mathfrak{H}}
\leqslant (1+\pi^{-1}) \sqrt{6} C_{7}  |t |^{p+1}, \quad  |t| \leqslant t^0.
\eqno(3.52)
$$

As a result, relations (2.26) and (3.52) imply the following estimate for the operator (3.25):
$$
\left\Vert J_4^{(3)}(t, \eps) \right\Vert_{\mathfrak{H} \to \mathfrak{H}}
\leqslant
(1+\pi^{-1})\sqrt{6} C_{7} |t |^{p+1} (c_* t^{2p} + \eps^{2p})^{-1},
\quad |t| \leqslant t^0.
$$
Using (3.31), we deduce
$$
\left\Vert J_4^{(3)}(t, \eps) \right\Vert_{\mathfrak{H} \to \mathfrak{H}}
\leqslant
C_8 \eps^{1-p}, \quad |t| \leqslant t^0,
\eqno(3.53)
$$
where
$$
C_8 = (1+\pi^{-1}) \sqrt{6} C_{7} c_*^{-1/2-1/2p}.
\eqno(3.54)
$$

\subsection{Completion of the proof of Theorem 3.1}
Relations (3.4), (3.15), (3.22), (3.32), (3.40), and (3.53) imply the required estimate (3.1) with the constant
$$
\check{C}_A= C_4 + C_5 + C_6 + C_8.
\eqno(3.55)
$$

\smallskip
\noindent
\textbf{Remark 3.4.} It is possible to write down a combersome explicit expression for the constant $\check{C}_A$
using relations (1.5), (2.14), (2.15), (2.19), (2.21), (3.16), (3.28), (3.33), (3.41), (3.47), (3.49), (3.51), (3.54), and (3.55).
For further application to differential operators, it is important to know how this constant
depends on the problem data. After possible overstating, $\check{C}_A$ can be considered as
a polynomial of the variables $\widetilde{C}$, $\|X_p\|$, $\delta^{-1/2p}$, and $c_*^{-1/2p}$
with positive coefficients depending only on $p$.

\section{Periodic differential operators in $\R^d$.\\ Direct integral expansion}

\subsection{Factorized operators of order $2p$ in $L_2(\R^d;\C^n)$}

In $L_{2}(\mathbb{R}^{d};\mathbb{C}^{n})$, we consider a differential operator
given formally by the expression
$$
A=b(\mathbf{D})^{*}g(\mathbf{x})b(\mathbf{D}).
\eqno(4.1)
$$
Here $g(\mathbf{x})$ is uniformly positive definite and bounded $(m\times m)$-matrix-valued function
(in general, $g(\mathbf{x})$ is a Hermitian matrix with complex entries):
$$
\begin{aligned}
g,\, g^{-1} & \in L_{\infty}(\mathbb{R}^{d}),
\\
g\left(\mathbf{x}\right) & \geqslant c \mathbf{1}_{m},\quad c >0,\quad \mathbf{x}\in\mathbb{R}^{d}.
\end{aligned}
\eqno(4.2)
$$
The operator $b\left(\mathbf{D}\right)$ is given by
$$
b\left(\mathbf{D}\right)=\sum_{\left|\alpha\right|=p}b_{\alpha}\mathbf{D}^{\alpha},
\eqno(4.3)
$$
where $b_{\alpha}$ are constant $\left(m\times n\right)$-matrices, in general, with complex entries.
It is assumed that $m \geqslant n$ and that the symbol
$b({\boldsymbol{\xi}})= \sum_{\left|\alpha\right|=p}b_{\alpha} {\boldsymbol \xi}^{\alpha}$ satisfies
$$
\mathrm{rank}\, b(\boldsymbol{\xi}) =n,\quad0\neq\boldsymbol{\xi}\in\mathbb{R}^{d}.
$$
This condition is equivalent to the following estimates
$$
\begin{aligned}
 & \alpha_{0}\mathbf{1}_{n}\leqslant b({\boldsymbol{\theta}})^{*}b(\boldsymbol{\theta})\leqslant\alpha_{1}\mathbf{1}_{n},\quad\boldsymbol{\theta}\in\mathbb{S}^{d-1},
\\
 & 0<\alpha_{0}\leqslant\alpha_{1}<\infty.
\end{aligned}
\eqno(4.4)
$$
Without loss of generality, we assume that
$$
\left|b_{\alpha}\right|\leqslant\alpha_{1}^{1/2},\quad\left|\alpha\right|=p.
\eqno(4.5)
$$

The precise definition of the operator  $A$ is given in terms of the quadratic form.
By (4.2), the matrix $g$ can be written as
$$
g(\mathbf{x})=h(\mathbf{x})^{*}h(\mathbf{x});
$$
where $h,\, h^{-1}\in L_{\infty}$. For instance, one can put $h=g^{1/2}$.

Consider the operator $X$ acting from $L_{2}(\mathbb{R}^{d};\mathbb{C}^{n})$
to $L_{2}(\mathbb{R}^{d};\mathbb{C}^{m})$ and defined by
$$
(X\mathbf{u})(\mathbf{x})=h(\mathbf{x})b(\mathbf{D})\mathbf{u}(\mathbf{x}),\quad\mathbf{u}\in\mathrm{Dom}\, X=H^{p}(\mathbb{R}^{d};\mathbb{C}^{n}),
$$
and the quadratic form
 $$
 a\left[\mathbf{u},\mathbf{u}\right]=\left\Vert X\mathbf{u}\right\Vert _{L_{2}(\mathbb{R}^{d})}^{2}=
 \intop_{\R^d} \langle g(\x)b(\D)\u(\x), b(\D)\u(\x)\rangle \,d\x,\quad \u \in H^p(\R^d;\C^n).
 \eqno(4.6)
$$
Using the Fourier transformation and relations (4.2), (4.4), it is easy to check that
$$
 c_0 \intop_{\R^d} |\D^p \u|^2\,d\x \leqslant a\left[\mathbf{u},\mathbf{u}\right] \leqslant c_1 \intop_{\R^d} |\D^p \u|^2\,d\x,
 \quad  \u \in H^p(\R^d;\C^n),
 \eqno(4.7)
$$
where the notation $|\D^p \u|^2:= \sum_{|\alpha|=p} |\D^\alpha \u|^2$ is used.
Indeed, by the Parceval identity, we have
$$
 \|g^{-1}\|^{-1}_{L_\infty} \intop_{\R^d} |b(\boldsymbol{\xi}) \wh{\u}(\boldsymbol{\xi})|^2\, d\boldsymbol{\xi} \leqslant
a\left[\mathbf{u},\mathbf{u}\right] \leqslant \|g\|_{L_\infty} \intop_{\R^d} |b(\boldsymbol{\xi}) \wh{\u}(\boldsymbol{\xi})|^2 \, d\boldsymbol{\xi},
$$
where $\wh{\u}(\bxi)$ is the Fourier image of the function $\u(\x)$.
Together with  (4.4) this yields
$$
 \alpha_0 \|g^{-1}\|^{-1}_{L_\infty}  \intop_{\R^d} |\boldsymbol{\xi}|^{2p}
 |\wh{\u}(\boldsymbol{\xi})|^2\, d\boldsymbol{\xi} \leqslant
a\left[\mathbf{u},\mathbf{u}\right] \leqslant
\alpha_1 \|g\|_{L_\infty} \intop_{\R^d} |\boldsymbol{\xi}|^{2p} |\wh{\u}(\boldsymbol{\xi})|^2 \, d\boldsymbol{\xi}.
\eqno(4.8)
$$
Using the elementary inequalities
$$
c_p' \sum_{|\alpha|=p} |\bxi^\alpha|^{2} \leqslant |\bxi|^{2p} \leqslant c_p'' \sum_{|\alpha|=p} |\bxi^\alpha|^{2},
\eqno(4.9)
$$
with the constants $c_p'$, $c_p''$ depending only on $d$ and $p$, we arrive at
the required inequalities  (4.7) with
$$
c_0 =  c_p' \alpha_0 \|g^{-1}\|^{-1}_{L_\infty},\quad
c_1=   c_p'' \alpha_1 \|g\|_{L_\infty}.
\eqno(4.10)
$$

Hence, the form (4.6) is closed and nonnegative.
The selfadjoint operator in $L_2(\R^d;\C^n)$ corresponding to this form is denoted by $A$.

\subsection{Lattices in $\mathbb{R}^{d}$}

In what follows, it is assumed that the matrix-valued functions $g$ and $h$ are \textit{periodic} with respect to
some lattice $\Gamma\subset\mathbb{R}^{d}$:
$$
g(\mathbf{x}+\mathbf{n})=g(\mathbf{x}),\qquad h(\mathbf{x}+\mathbf{n})=h(\mathbf{x}),
\quad\x \in \R^d,\quad \mathbf{n}\in\Gamma.
$$

Let $\mathbf{n}_{1},\dots,\mathbf{n}_{d}$ be a basis in $\mathbb{R}^{d}$ generating the lattice $\Gamma$:
$$
\Gamma=\left\{ \mathbf{n}\in\mathbb{R}^{d}:\,\mathbf{n}=\sum_{i=1}^{d}l_{i}\mathbf{n}_{i},\, l_{i}\in\mathbb{Z}\right\} ,
$$
and let $\Omega$ be the elementary cell of the lattice $\Gamma$:
$$
\Omega=\left\{ \mathbf{x}\in\mathbb{R}^{d}:\,\mathbf{x}=\sum_{i=1}^{d}t_{i}\mathbf{n}_{i},\,0<t_{i}<1\right\} .
$$

The basis $\mathbf{s}^{1},\dots,\mathbf{s}^{d}$ in $\mathbb{R}^{d}$ dual to the basis
$\mathbf{n}_{1},\dots,\mathbf{n}_{d}$ is defined by the relations
$\left\langle \mathbf{s}^{i},\mathbf{n}_{j}\right\rangle _{\mathbb{R}^{d}}=2\pi\delta_{j}^{i}$.
This basis generates the \textit{lattice $\widetilde{\Gamma}$ dual to $\Gamma$}:
$$
\widetilde{\Gamma}=\left\{ \mathbf{s}\in\mathbb{R}^{d}:\,\mathbf{s}=\sum_{i=1}^{d}q_{i}\mathbf{s}^{i},\, q_{i}\in\mathbb{Z}\right\} .
$$
Instead of the cell of the dual lattice, it is convenient to consider the \emph{central Brillouin zone}
$$
\widetilde{\Omega}=\left\{ \mathbf{k}\in\mathbb{R}^{d}\,:\,\left|\mathbf{k}\right|<\left|\mathbf{k}-\mathbf{s}\right|,\,0\neq\mathbf{s}\in\widetilde{\Gamma}\right\} ,
$$
which is a fundamental domain of the lattice $\widetilde{\Gamma}$.
We use the notation  $|\Omega|=\mathrm{meas}\,\Omega$,
$|\widetilde{\Omega}|=\mathrm{meas}\,\widetilde{\Omega}$ and note that
$|\Omega| \, |\widetilde{\Omega}|=\left(2\pi\right)^{d}$. Let $r_{0}$ be the maximal radius of the ball
containing in $\mathrm{clos}\,\widetilde{\Omega}$. Then
$$
2r_{0}=\min\left|\mathbf{s}\right|,\quad0\neq\mathbf{s}\in\widetilde{\Gamma}.
\eqno(4.11)
$$
Denote
$$
\mathcal{B}\left(r\right)=\left\{ \mathbf{k}\in\mathbb{R}^{d}\,:\,\left|\mathbf{k}\right|\leqslant r\right\} ,\quad r>0.
$$

With the lattice $\Gamma$, a discrete Fourier transformation $\{\widehat{v}_{\mathbf{s}} \}_{\mathbf{s} \in \wt{\Gamma}}
\mapsto v$ is associated:
$$
v(\mathbf{x})=\left|\Omega\right|^{-1/2}\sum_{\mathbf{s}\in\widetilde{\Gamma}}\widehat{v}_{\mathbf{s}}\exp\left(i\left\langle \mathbf{s},\mathbf{x}\right\rangle \right),\quad\mathbf{x}\in\Omega,
$$
which is a unitary mapping of $l_{2}(\widetilde{\Gamma})$ onto $L_{2}(\Omega)$:
$$
\intop_{\Omega}\left|v\left(\mathbf{x}\right)\right|^{2}d\mathbf{x}=\sum_{\mathbf{s}\in\widetilde{\Gamma}}\left|\widehat{v}_{\mathbf{s}}\right|^{2}.
$$

By $\widetilde{W}^{s}_q(\Omega;\mathbb{C}^{n})$ we denote the subspace of all functions in
$W_q^{s}(\Omega;\mathbb{C}^{n})$ whose $\Gamma$-periodic extension to $\mathbb{R}^{d}$
belongs to $W_{q,\mathrm{loc}}^{s}(\mathbb{R}^{d};\mathbb{C}^{n})$.
If  $q=2$, we use the notation $\widetilde{H}^{s}(\Omega;\mathbb{C}^{n})$.

\subsection{The operators $X(\k)$ and $A(\k)$ in $L_2(\Omega;\C^n)$}

Let $\mathbf{k}\in\mathbb{R}^{d}$.
We consider the operator $X\left(\mathbf{k}\right):\, L_{2}\left(\Omega;\mathbb{C}^{n}\right)\to L_{2}\left(\Omega;\mathbb{C}^{m}\right)$
defined on the domain
$$
\mathrm{Dom}\,X\left(\mathbf{k}\right)=\widetilde{H}^{p}(\Omega;\mathbb{C}^{n})
$$
by the relation
$$
(X(\mathbf{k})\u)(\x)= h(\mathbf{x})b(\mathbf{D}+\mathbf{k})\u(\x).
\eqno(4.12)
$$
Consider the quadratic form $a\left(\mathbf{k}\right)$ given by
$$
\begin{aligned}
a(\mathbf{k})\left[\mathbf{u},\mathbf{u}\right]=\left\Vert X(\mathbf{k})\mathbf{u}\right\Vert _{L_{2}(\Omega)}^{2}= \intop_\Omega \langle g(\x) b(\D+\k) \u, b(\D+\k)\u \rangle\,d\x,
\\ \mathbf{u}
\in \widetilde{H}^{p}(\Omega;\mathbb{C}^{n}).
\end{aligned}
\eqno(4.13)
$$
Using the discrete Fourier transformation and relations (4.2), (4.4), it is easy to check that
$$
\begin{aligned}
\alpha_0 \|g^{-1}\|^{-1}_{L_\infty} a_*(\k)\left[\mathbf{u},\mathbf{u}\right]
 \leqslant a(\k)\left[\mathbf{u},\mathbf{u}\right] \leqslant
\alpha_1 \|g\|_{L_\infty} a_*(\k)\left[\mathbf{u},\mathbf{u}\right],
\\
\u \in \widetilde{H}^p(\Omega;\C^n),
\end{aligned}
\eqno(4.14)
$$
for any $\k\in \R^d$, where
$$
a_*(\k)\left[\mathbf{u},\mathbf{u}\right] := \sum_{{\mathbf s} \in \wt{\Gamma}} |{\mathbf s}+\k|^{2p} |\wh{\u}_{\mathbf s}|^2,
\quad  \u \in \widetilde{H}^p(\Omega;\C^n).
\eqno(4.15)
$$
Together with (4.9) this yields
$$
 \begin{aligned}
 c_0 \intop_{\Omega} |(\D+\k)^p \u|^2\,d\x \leqslant a(\k)\left[\mathbf{u},\mathbf{u}\right] \leqslant c_1 \intop_{\Omega}
 |(\D+\k)^p \u|^2\,d\x,
 \\
   \u \in \widetilde{H}^p(\Omega;\C^n),
\end{aligned}
$$
where the constants $c_0$, $c_1$ are defined by (4.10).
Hence, the operator $X(\k)$ is closed, and the form (4.13) is closed and nonnegative.
The selfadjoint operator in $L_2(\Omega;\C^n)$ corresponding to the form $a(\mathbf{k})$ is denoted by $A(\mathbf{k})$.
Formally, we have
$$
A(\k) = b(\D+\k)^* g(\x) b(\D+\k).
$$

\subsection{Direct integral expansion for the operator $A$}

First, the Gelfand transformation $\mathcal{U}$ is defined on the Schwartz class $\mathcal{S}(\mathbb{R}^{d};\mathbb{C}^{n})$
by the relation
\begin{align*}
 & \widetilde{\mathbf{v}}(\mathbf{k},\mathbf{x})=\left(\mathcal{U}\mathbf{v}\right)(\mathbf{k},\mathbf{x})= |\widetilde{\Omega}|^{-1/2}\sum_{\mathbf{n}\in\Gamma}\mathrm{exp}(-i \langle \mathbf{k},\mathbf{x}+\mathbf{n}\rangle)
\mathbf{v}(\mathbf{x}+\mathbf{n}),\\
 & \mathbf{v}\in\mathcal{S}(\mathbb{R}^{d};\mathbb{C}^{n}),\quad
\mathbf{x}\in\mathbb{R}^{d},\quad\mathbf{k}\in\mathbb{R}^{d}.
\end{align*}
Then
$$
\intop_{\widetilde{\Omega}}\intop_{\Omega}\left|\widetilde{\mathbf{v}}\left(\mathbf{k},\mathbf{x}\right)\right|^{2}d\mathbf{x}\, d\mathbf{k}=\intop_{\mathbb{R}^{d}}\left|\mathbf{v}(\mathbf{x})\right|^{2}d\mathbf{x},\quad\widetilde{\mathbf{v}}=\mathcal{U}\mathbf{v},
$$
and $\mathcal{U}$ is extended by continuity up to unitary mapping
$$
\mathcal{U}\,:\, L_{2}(\mathbb{R}^{d};\mathbb{C}^{n})\to\intop_{\widetilde{\Omega}}\oplus L_{2}(\Omega;\mathbb{C}^{n})d\mathbf{k}=:\mathcal{K}.
\eqno(4.16)
$$

The relation $\v \in H^p(\R^d;\C^n)$ is equivalent to the fact that  $\wt{\v}(\k,\cdot)\in \wt{H}^p(\Omega;\C^n)$
for a.~e. $\k \in \wt{\Omega}$ and
$$
\intop_{\wt{\Omega}}\intop_\Omega \left(|(\D+ \k)^p \wt{\v}(\k,\x)|^2 + |\wt{\v}(\k,\x)|^2\right) \,d\x\,d\k < \infty.
$$
Under the Gelfand transformation $\mathcal U$, the operator of multiplication by a bounded $\Gamma$-periodic
function in $L_2(\R^d;\C^n)$ turns into the operator of multiplication by the same function on the fibers
of the direct integral $\mathcal K$ from (4.16). On these fibers, the action of the operator  $b(\D)$ on $\v \in H^p(\R^d;\C^n)$ turns into
the action of the operator $b(\D+\k)$ on $\wt{\v}(\k,\cdot) \in \wt{H}^p(\Omega;\C^n)$.

It follows that the form (4.6) can be written as
$$
a\left[\mathbf{v},\mathbf{v}\right]=\intop_{\widetilde{\Omega}}
a(\mathbf{k}) \left[\widetilde{\mathbf{v}}(\mathbf{k},\cdot), \widetilde{\mathbf{v}}(\mathbf{k},\cdot) \right]\,d\mathbf{k},\quad \v \in H^p(\R^d;\C^n).
$$
Thus,  the operator $A$ is unitarily equivalent (with the affinity $\mathcal U$)
to the direct integral of the operators $A(\mathbf{k})$:
$$
\mathcal{U}A\mathcal{U}^{-1}=\intop_{\widetilde{\Omega}}\oplus A(\mathbf{k}) d\mathbf{k}.
\eqno(4.17)
$$

\section{Differential operators $A(\k)$ on the cell $\Omega$.\\ Application of the abstract scheme}

\subsection{Investigation of the operators $X(\mathbf{k})$ and $A(\mathbf{k})$}

Our goal is to check that the abstract scheme can be applied to the operators $A(\k)$.
As in [BSu1], we put $t:=\left|\mathbf{k}\right|$ and $\boldsymbol{\theta}:=\mathbf{k}/t$.
The operators $X(\mathbf{k})=: X(t,\boldsymbol{\theta})$ and $A(\mathbf{k})=: A(t,\boldsymbol{\theta})$
depend on the onedimensional parameter $t$ and the additional parameter
$\boldsymbol{\theta} \in {\mathbb S}^{d-1}$ (the latter  was absent in the abstract scheme).
We will try to make our constructions and estimates uniform with respect to $\boldsymbol{\theta}$.

By (4.3) and (4.12),
$$
\begin{aligned}
&X(\mathbf{k})=h \sum_{|\alpha|=p} b_{\alpha} \left( \mathbf{D}+\mathbf{k} \right)^{\alpha}=
h \sum_{|\alpha|=p} b_{\alpha} \sum_{\beta\leqslant\alpha} C_{\alpha}^{\beta} {\k}^{\alpha - \beta}\mathbf{D}^{\beta}
\\
&=h \sum_{|\alpha|=p} b_{\alpha}\sum_{\beta\leqslant\alpha} C_{\alpha}^{\beta} t^{|\alpha-\beta|} \boldsymbol{\theta}^{\alpha-\beta} \mathbf{D}^{\beta}.
\end{aligned}
$$
Hence, the operator $X(\mathbf{k})$ can be written as
$$
X(\mathbf{k})=X(t,\boldsymbol{\theta})= X_0 + \sum_{j=1}^p t^j X_j(\boldsymbol{\theta}).
\eqno(5.1)
$$
Here the operator
$$
X_{0}=h \sum_{|\alpha|=p} b_{\alpha}\mathbf{D}^{\alpha}= h b(\mathbf{D})
\eqno(5.2)
$$
is closed on the domain
$$
\mathrm{Dom}\, X_{0}=\widetilde{H}^{p}(\Omega;\mathbb{C}^{n}),
\eqno(5.3)
$$
the ``intermediate'' operators $X_{j}(\boldsymbol{\theta})$, $j=1,\dots,p-1$, are given by
$$
X_{j}(\boldsymbol{\theta})=h \sum_{\left|\alpha\right|=p} b_{\alpha}
\sum_{ \beta\leqslant\alpha,\,\left|\beta\right|=p-j } C_{\alpha}^{\beta} \boldsymbol{\theta}^{\alpha-\beta}\mathbf{D}^{\beta}
\eqno(5.4)
$$
on the domains
$$
\mathrm{Dom}\,X_{j}(\boldsymbol{\theta})=\widetilde{H}^{p-j}(\Omega;\mathbb{C}^{n}),
\eqno(5.5)
$$
and the operator
$$
X_{p}(\boldsymbol{\theta})=h \sum_{|\alpha|=p} b_{\alpha}\boldsymbol{\theta}^{\alpha}= h b(\boldsymbol{\theta})
$$
 is bounded from $L_{2}(\Omega;\mathbb{C}^{n})$ to $L_{2}(\Omega;\mathbb{C}^{m})$.

From (5.3) and (5.5) it follows that Condition 1.1 is satisfied:
$$
\mathrm{Dom}\,X_{0} \subset \mathrm{Dom}\,X_{j}\left(\boldsymbol{\theta}\right) \subset
\mathrm{Dom}\,X_{p}\left(\boldsymbol{\theta}\right)=L_{2}(\Omega;\mathbb{C}^{n}),\quad j=1,\dots,p-1.
$$
By  (4.4), we have
$$
\left\Vert X_{p}(\boldsymbol{\theta}) \right\Vert
\leqslant  \alpha_{1}^{1/2} \left\Vert g \right\Vert^{1/2} _{L_{\infty}}.
\eqno(5.6)
$$

Now we consider the kernels of the operators $X_{j}(\boldsymbol{\theta})$.

\smallskip
\noindent
\textbf{Proposition 5.1.}
\textit{The kernel of $X_0$ consists of constant vector-valued functions}:
$$
\mathfrak{N}:=\mathrm{Ker}X_{0}=\left\{ \mathbf{u}\in L_{2}(\Omega;\mathbb{C}^{n}):\ \mathbf{u}(\mathbf{x}) =\mathbf{c}\in\mathbb{C}^{n}\right\}.
\eqno(5.7)
$$
\textit{For  $j=1,\dots,p-1$ we have}
$$
\mathfrak{N} \subset\mathrm{Ker}\, X_{j}(\boldsymbol{\theta}), \quad \bt \in {\mathbb S}^{d-1}.
\eqno(5.8)
$$

\smallskip
\noindent
\textbf{Proof.}
Let $\mathbf{u}\in\mathfrak{N}$. It means that $\mathbf{u}\in \wt{H}^p(\Omega;\C^n)$ and
$b(\mathbf{D}) \mathbf{u}=0.$ Using the Parvceval identity for the discrete Fourier transformation,
we write this condition as
$$
0=\intop_{\Omega}\left|b(\mathbf{D})\mathbf{u}(\mathbf{x})\right|^{2}d\mathbf{x}
=\sum_{\mathbf{s}\in\widetilde{\Gamma}}\left|b(\mathbf{s})\widehat{\mathbf{u}}_{\mathbf{s}}\right|^{2}=\sum_{\mathbf{s}\in\widetilde{\Gamma}}\left\langle b(\mathbf{s})^{*}b(\mathbf{s})\widehat{\mathbf{u}}_{\mathbf{s}},\widehat{\mathbf{u}}_{\mathbf{s}}\right\rangle _{\mathbb{C}^{n}}.
\eqno(5.9)
$$
By (4.4), relation (5.9) is equivalent to
$$
|\mathbf{s}|^{2p} \left|\widehat{\mathbf{u}}_{\mathbf{s}}\right|^2 =0,\quad\mathbf{s}\in\widetilde{\Gamma}.
$$
In its turn, this is equivalent to the fact that all the Fourier coefficients $\widehat{\mathbf{u}}_{\mathbf{s}}$
besides $\widehat{\mathbf{u}}_{\mathbf{0}}$ are equal to zero, i/~e.,
$\mathbf{u}(\mathbf{x}) = \mathbf{c}\in\mathbb{C}^{n}$.

By (5.4), relations (5.8) are obvious.
$\square$
\smallskip

The orthogonal projection of  $L_2(\Omega;\C^n)$ onto  $\mathfrak{N}$
acts as averaging over the cell:
$$
P\mathbf{u}=\left|\Omega\right|^{-1}\intop_{\Omega}\mathbf{u}(\mathbf{x})\, d\mathbf{x},\quad\mathbf{u}\in L_2(\Omega;\C^n).
$$

Let $n_* = {\rm Ker}\,X_0^*$. Relation $m\geqslant n$ ensures that $n_* \geqslant n$.
Moreover, since
$$
\NN_* = {\rm Ker}\,X_0^* = \{{\mathbf q} \in L_2(\Omega;\C^m):\ h^*{\mathbf q} \in \wt{H}^p(\Omega;\C^m):\
b(\D)^* h^* {\mathbf q} =0\},
$$
then the following alternative takes place: either $n_*=\infty$ (if $m>n$), or $n_*=n$ (if $m=n$).

Now we check that Condition 1.2 is satisfied.

\smallskip\noindent\textbf{Proposition 5.2.}
\textit{For $j=1,\dots,p-1$ we have
$$
\left\Vert X_{j}(\boldsymbol{\theta}) \mathbf{u}\right\Vert _{L_{2}(\Omega)}
\leqslant\widetilde{C}_{j}\left\Vert X_{0}\mathbf{u}\right\Vert _{L_{2}(\Omega)},
\quad \mathbf{u}\in\widetilde{H}^{p}(\Omega;\mathbb{C}^{n}).
\eqno(5.10)
$$
Here the constants $\widetilde{C}_{j}$ do not depend on the parameter $\boldsymbol{\theta}\in {\mathbb S}^{d-1}$
and depend only on $d$, $p$, $j$, $\left\Vert g\right\Vert _{L_{\infty}}$, $\left\Vert g^{-1}\right\Vert _{L_{\infty}}$,
$\alpha_{0}$, $\alpha_{1}$, and $r_{0}$.}
\smallskip

\noindent\textbf{Proof.}
By (4.5) and (5.4),
$$
\left\Vert X_{j}(\boldsymbol{\theta})\mathbf{u}\right\Vert _{L_{2}\left(\Omega\right)}\leqslant
\alpha_{1}^{1/2} \left\Vert g\right\Vert _{L_{\infty}}^{1/2}
\sum_{|\alpha|\leqslant p} \sum_{ \beta\leqslant\alpha,\,\left|\beta\right|=p-j } C_{\alpha}^{\beta} \|\mathbf{D}^{\beta} \u\|_{L_{2}(\Omega)}.
\eqno(5.11)
$$
We expand a function $\mathbf{u}\in\widetilde{H}^{p}(\Omega;\mathbb{C}^{n})$ in the Fourier series
$$
\mathbf{u}(\x)=\left|\Omega\right|^{-1/2}\sum_{\mathbf{s}\in\widetilde{\Gamma}}
\widehat{\mathbf{u}}_{\mathbf{s}} e^{i\left\langle \mathbf{x},\mathbf{s}\right\rangle }.
\eqno(5.12)
$$
By (4.11), for $j=1,\dots, p-1$ we have
$$
|\mathbf{s}^{\beta}|^{2}\leqslant |\mathbf{s}|^{2\left|\beta\right|}\leqslant
\left(2r_{0}\right)^{-2j} |\mathbf{s}|^{2p},\quad\mathbf{0}\neq\mathbf{s}\in\widetilde{\Gamma},\quad\left|\beta\right|=p-j.
\eqno(5.13)
$$
From the Parceval identity for the Fourier series and from (5.13) we deduce
$$
\|\mathbf{D}^{\beta} \u\|_{L_{2}(\Omega)}^2 =
\sum_{\mathbf{0} \neq\mathbf{s}\in\widetilde{\Gamma}}
|\mathbf{s}^{\beta} \widehat{\mathbf{u}}_{\mathbf{s}} |^{2}
\leqslant (2r_0)^{-2j}
\sum_{\mathbf{0} \neq\mathbf{s}\in\widetilde{\Gamma}} |\mathbf{s}|^{2p}| \widehat{\mathbf{u}}_{\mathbf{s}} |^{2},
\quad |\beta| = p-j.
\eqno(5.14)
$$
Next, from the definition  (5.2), (5.3) of $X_{0}$, expansion (5.12), and the lower estimate (4.4) it follows that
$$
\left\Vert X_{0}\mathbf{u}\right\Vert _{L_{2}\left(\Omega\right)}^{2}\geqslant
\left\Vert g^{-1}\right\Vert _{L_{\infty}}^{-1} \alpha_{0} \sum_{\mathbf{s}\in\widetilde{\Gamma}}
\left|\mathbf{s}\right|^{2p}\left|\widehat{\mathbf{u}}_{\mathbf{s}}\right|^{2},
\quad \mathbf{u}\in\widetilde{H}^{p}(\Omega;\mathbb{C}^{n}).
\eqno(5.15)
$$
As a result, combining (5.11), (5.14), and (5.15), we arrive at the required inequality (5.10) with the constant
$$
\widetilde{C}_{j}=  \alpha_1^{1/2} \alpha_0^{-1/2}
\|g\|_{L_\infty}^{1/2} \|g^{-1}\|_{L_\infty}^{1/2} (2r_0)^{-j} \biggl( \sum_{|\alpha| \leqslant p}
\sum_{ \beta\leqslant\alpha,\,\left|\beta\right|=p-j } C_{\alpha}^{\beta}\biggr).
\eqno(5.16)
$$
$\square$

\smallskip

From the compactness of the embedding of $\mathrm{Dom}\, a(0)=\widetilde{H}^{p}(\Omega;\mathbb{C}^{n})$
into $L_{2}(\Omega;\mathbb{C}^{n})$ it follows that the spectrum of the operator
$A(0)$ is discrete. The point $\lambda_{0}=0$ is an isolated eigenvalue of the operator $A(0)=X_0^*X_0$
of multiplicity $n$; the corresponding eigenspace $\NN$ consists of constant vector-valued functions (see (5.7)).
We estimate the distance $d^0$ from the point $\lambda_{0}=0$ to the rest of the spectrum of $A(0)$.
From (4.11) and (5.15) it follows that
$$
\begin{aligned}
a(0)[\mathbf{u}, \mathbf{u}]\geqslant \alpha_{0}\left\Vert g^{-1}\right\Vert _{L_{\infty}}^{-1}\left(2r_{0}\right)^{2p}
\left\Vert \mathbf{u}\right\Vert _{L_{2}(\Omega)}^{2},
\\
\mathbf{u}\in\widetilde{H}^{p}(\Omega;\mathbb{C}^{n}),\quad\intop_{\Omega}\mathbf{u}(\mathbf{x})\,d\mathbf{x}=0.
\end{aligned}
$$
Thus,
$$
d^{0}\geqslant\alpha_{0}\left\Vert g^{-1}\right\Vert _{L_{\infty}}^{-1}\left(2r_{0}\right)^{2p}.
\eqno(5.17)
$$

According to the abstract scheme, we should fix a positive number $\delta$ such that
$\delta\leqslant\min\left\{ d^{0}/36,\,1/4\right\}$.
Taking (5.17) into account, we put
$$
\delta=\min\left\{ \alpha_{0}\left\Vert g^{-1}\right\Vert _{L_{\infty}}^{-1}\left(2r_{0}\right)^{2p}/36,\,1/4\right\}.
\eqno(5.18)
$$
Inequalities (5.10) allow us to choose $\widetilde{C}$ (see (1.1)) as follows
$$
\widetilde{C}=\max\left\{ 1,\,\widetilde{C}_{1},\,\dots,\,\widetilde{C}_{p-1}\right\},
\eqno(5.19)
$$
where the constants $\widetilde{C}_{j}$ are defined by  (5.16).

Now the constant $\widehat{C}=\max\left\{ \left(p-1\right)\widetilde{C},\,\left\Vert X_{p}\left(\boldsymbol{\theta}\right)\right\Vert \right\} $
(see (1.5)) depends on  $\boldsymbol{\theta}$.
Using (5.6), we take the following overstated value
$$
\widehat{C}=\max\left\{ \left(p-1\right)\widetilde{C},\, \alpha_{1}^{1/2} \left\Vert g\right\Vert _{L_{\infty}}^{1/2} \right\},
$$
which does not depend on $\boldsymbol{\theta}$.
According to (1.4), we put
$$
t^{0}=\frac{\delta^{1/2}}{\wh{C}}= \frac{\delta^{1/2}}{\max\left\{ \left(p-1\right)\widetilde{C},\,
\alpha_{1}^{1/2} \left\Vert g\right\Vert _{L_{\infty}}^{1/2} \right\} }.
\eqno(5.20)
$$

\subsection{Incorporation of the operators $A(t,\bt)$ in the framework of the abstract scheme}

We apply the abstract scheme putting
$$
\mathfrak{H}=L_{2}(\Omega;\mathbb{C}^{n}),\quad\mathfrak{H}_{*}=L_{2}(\Omega;\mathbb{C}^{m}).
$$
The role of the polynomial pencil $X(t)$ is played by
$X(t,\boldsymbol{\theta}):= X(\k)=X(t\boldsymbol{\theta})$ (see (5.1)); this pencil depends also on
$\boldsymbol{\theta}$. In Subsection~5.1 it was checked that Conditions 1.1 and 1.2 are satisfied.
The role of  $A(t)$ is played by the operator $A(t,\bt):= A\left(\mathbf{k}\right)=A\left(t\boldsymbol{\theta}\right)$.
According to the definition of $A(\k)$ (see Subsection~4.3), we have
$$
A(t,\bt)=X(t,\bt)^{*} X(t,\bt).
$$
In Subsection 5.1 it was checked that Condition 1.3 is also satisfied.
The kernel $\NN = \Ker A(0)= \Ker X_0$ is described in (5.7).

It remains to check that Condition 2.2 is satisfied. Denote by $E_{j}(\mathbf{k})$, $j\in\mathbb{N}$,
$\mathbf{k}\in\widetilde{\Omega}$, the consecutive eigenvalues of the operator $A(\mathbf{k})$
counted with multiplicities.

We rely on the twosided estimates (4.14) for the form $a(\mathbf{k})$ in terms of the auxiliary form (4.15).
The selfadjoint operator in $\H$ corresponding to the form (4.15) is denoted by $A_*(\k)$,
and its consecutive eigenvalues are denoted by $E_{j}^{0}\left(\mathbf{k}\right)$, $j\in\mathbb{N}$.
For a different way of enumeration, the eigenvalues of $A_*(\k)$ can be found from (4.15): they coincide with
the numbers $|\mathbf{s}+\mathbf{k}|^{2p}$, $\mathbf{s} \in \wt{\Gamma}$, and each such
eigenvalue is of multiplicity $n$ (they are enumerated by the index $\mathbf{s} \in \wt{\Gamma}$).
Then it is easily seen that
\begin{align}
E_{l}^{0}\left(\mathbf{k}\right) & =\left|\mathbf{k}\right|^{2p},\quad l=1,\dots,n,\quad\mathbf{k}\in\mathrm{clos}\,\widetilde{\Omega},
\tag{5.21}
\\
E_{1}^{0}\left(\mathbf{k}\right) & \geqslant r^{2p},\quad\mathbf{k}\in\mathrm{clos}\,\widetilde{\Omega}\setminus\mathcal{B}\left(r\right),\quad0<r\leqslant r_{0}.
\tag{5.22}
\end{align}
Combining the lower estimate (4.14),  (5.21), and (5.22), we obtain
\begin{align}
E_{l}\left(\mathbf{k}\right) & \geqslant\alpha_{0}\left\Vert g^{-1}\right\Vert _{L_{\infty}}^{-1}\left|\mathbf{k}\right|^{2p},\quad l=1,\dots,n,\quad\mathbf{k}\in\mathrm{clos}\,\widetilde{\Omega},
\nonumber
\\
E_{1}\left(\mathbf{k}\right) & \geqslant\alpha_{0}\left\Vert g^{-1}\right\Vert _{L_{\infty}}^{-1}r^{2p},\quad\mathbf{k}\in\mathrm{clos}\,\widetilde{\Omega}\setminus\mathcal{B}\left(r\right),\quad0<r\leqslant r_{0},
\tag{5.23}
\end{align}
which implies that Condition 2.2 is satisfied with the constant
$$
c_{*}=\alpha_{0}\left\Vert g^{-1}\right\Vert _{L_{\infty}}^{-1}.
\eqno(5.24)
$$

\subsection{Construction of the spectral germ}

The operators $Z$, $R$, and $S$ introduced in Subsection~1.2
now depend on  $\boldsymbol{\theta}$. To construct them, we introduce the auxiliary operator $\Lambda$. Let
$$
\mathfrak{M}=\left\{ \mathbf{w}\in L_2(\Omega;\C^m): \ \w(\x)=\mathbf{C}\in\mathbb{C}^{m}\right\}
$$
be the subspace of constant vector-valued functions in $\H_*=L_2(\Omega;\C^m)$.

By definition, the operator $\Lambda:\mathfrak{M}\to\mathfrak{H}$ takes a vector
$\mathbf{C}\in\mathfrak{M}$ to the \emph{weak $\Gamma$-periodic solution}
$\mathbf{v}_{\mathbf{C}}\in\widetilde{H}^{p}(\Omega;\mathbb{C}^{n})$ of the problem
$$
b(\mathbf{D})^{*}g(\mathbf{x})\left(b(\mathbf{D})\mathbf{v}_{\mathbf{C}}(\mathbf{x})+\mathbf{C}\right)=0,\quad\intop_{\Omega}\mathbf{v}_{\mathbf{C}}(\x)\,d\mathbf{x}=0.
\eqno(5.25)
$$
Problem (5.25) with $\mathbf{C}=b(\boldsymbol{\theta})\mathbf{c}$, $\mathbf{c}\in \C^n$,
is realization of the problem (1.8), (1.9) (now $\omega = \mathbf{c}\in \NN$).
Let $\mathbf{e}_{1},\dots,\mathbf{e}_{m}$ be the standard orthonormal basis in
$\mathbb{C}^{m}$, and let $\mathbf{v}_{j}=\mathbf{v}_{\mathbf{e}_{j}}$.
In the standard basis in $\mathbb{C}^{n}$, the vector-valued functions $\mathbf{v}_{j}(\x)$
can be written as the columns of length  $n$. Let $\Lambda\left(\mathbf{x}\right)$ be the $\left(n\times m\right)$-matrix with the columns
$\mathbf{v}_{1}(\mathbf{x}),\dots,\mathbf{v}_{m}(\mathbf{x})$.
Then $\Lambda$ is the operator of multiplication by the matrix-valued function $\Lambda(\mathbf{x})$.
Note that the matrix-valued function $\Lambda(\mathbf{x})$ is a $\Gamma$-periodic solution of the problem
$$
b(\mathbf{D})^{*}g(\mathbf{x})\left(b(\mathbf{D})\Lambda(\mathbf{x})+\mathbf{1}_{m}\right)=0,\qquad\intop_{\Omega}\Lambda(\mathbf{x})\,d\mathbf{x}=0.
\eqno(5.26)
$$
According to (1.10), we obtain
\begin{align*}
(Z(\boldsymbol{\theta})\mathbf{c})(\x)&=\Lambda(\x) b(\boldsymbol{\theta})\mathbf{c},  \quad\mathbf{c}\in\mathbb{C}^{n}=\mathfrak{N},\\
Z(\boldsymbol{\theta})\mathbf{u}&=0,  \quad\mathbf{u}\in\mathfrak{N}^{\perp}.
\end{align*}
Thus,
$$
Z(\boldsymbol{\theta})=\Lambda b(\boldsymbol{\theta})P.
\eqno(5.27)
$$

According to (1.12) and (1.13), the operator $R(\boldsymbol{\theta})$ takes the form
$$
(R(\boldsymbol{\theta})\mathbf{c})(\x)=h(\x) \left(b(\mathbf{D})\Lambda(\x)+\mathbf{1}_{m}\right)b(\boldsymbol{\theta})\mathbf{c},\quad\mathbf{c}\in\mathbb{C}^{n}=\mathfrak{N}.
$$
Then the spectral germ
$S(\boldsymbol{\theta})=R(\boldsymbol{\theta})^{*}R(\boldsymbol{\theta}):\mathfrak{N}\to\mathfrak{N}$
is given by
$$
S(\boldsymbol{\theta})=Pb(\boldsymbol{\theta})^{*} \left(b(\mathbf{D})\Lambda +\mathbf{1}_{m}\right)^{*} g \left(b(\mathbf{D})\Lambda +\mathbf{1}_{m}\right)b(\boldsymbol{\theta})|_{\mathfrak{N}}
$$
and acts as multiplication by the matrix
$b(\boldsymbol{\theta})^{*}g^{0}b(\boldsymbol{\theta})$, where
$$
g^{0}=\left|\Omega\right|^{-1}
\intop_{\Omega}\left(b(\mathbf{D})\Lambda(\mathbf{x})+\mathbf{1}_{m}\right)^{*}g(\x)
\left(b(\mathbf{D})\Lambda(\mathbf{x})+\mathbf{1}_{m}\right)\,d\mathbf{x}.
\eqno(5.28)
$$
Taking (5.26) into account, we rewrite $g^{0}$  as
$$
g^{0}=\left|\Omega\right|^{-1}\intop_{\Omega}\widetilde{g}(\mathbf{x})\,d\mathbf{x},
\eqno(5.29)
$$
where
$$
\widetilde{g}(\mathbf{x})=g(\mathbf{x}) \left(b(\mathbf{D})\Lambda(\mathbf{x})+\mathbf{1}_{m}\right).
\eqno(5.30)
$$
The constant matrix (5.29) is called the \emph{effective matrix}.
Automatically, $g^0$ is positive definite which can be easily seen from the representation (5.28).
Thus, we have proved that the spectral germ of the operator family $A(t,\boldsymbol{\theta})$ is represented as
$$
S(\boldsymbol{\theta})=b(\boldsymbol{\theta})^{*}g^{0}b(\boldsymbol{\theta}).
$$

As has been mentioned in Subsection~2.3, Condition 2.2 implies that
$$
S(\boldsymbol{\theta}) \geqslant c_{*}I_{\mathfrak{N}},
\eqno(5.31)
$$
where the constant $c_{*}$ is defined by (5.24) and does not depend on $\boldsymbol{\theta}$.
So, the spectral germ  $S(\boldsymbol{\theta})$ is nondegenerate.

\subsection{The effective operator. The properties of the effective matrix}

We put
$$
S(\mathbf{k}) = t^{2p}S(\boldsymbol{\theta})= b(\mathbf{k})^{*}g^{0}b(\mathbf{k}),
\quad \mathbf{k}=t\boldsymbol{\theta}\in\mathbb{R}^{d};
\eqno(5.32)
$$
this is the symbol of the differential operator
$$
A^{0}=b(\mathbf{D})^{*}g^{0}b(\mathbf{D})
\eqno(5.33)
$$
with constant coefficients called the \emph{effective operator} for $A$.
Relations (5.31) and (5.32) imply the following estimate for the symbol of the effective operator:
$$
b(\k)^* g^0 b(\k) \geqslant c_* |\k |^{2p} \1_n,\quad \k \in \R^d.
\eqno(5.34)
$$

Now we discuss some properties of the effective matrix $g^{0}$.

\smallskip\noindent\textbf{Proposition 5.3.}
\textit{Denote
$$
\overline{g}:  =\left|\Omega\right|^{-1} \int_{\Omega}g(\mathbf{x})d\mathbf{x},\quad
\underline{g}:  =\left(\left|\Omega\right|^{-1}\int_{\Omega}g(\mathbf{x})^{-1}d\mathbf{x}\right)^{-1}.
$$
Then the effective matrix  $g^{0}$ satisfies the following inequalities:
$$
\underline{g}\leqslant g^{0}\leqslant\overline{g}.
\eqno(5.35)
$$
If  $m=n$, then $g^0=\underline{g}$.}

\smallskip\noindent\textbf{Proof.}  The proof is similar to that of Theorem~1.5 from [BSu1, Chapter 3], where the second order DO's were
studied. Let ${\mathbf C}\in \C^m$, and let $\mathbf{v}_{\mathbf{C}}$ be the periodic solution of the problem (5.25).
Obviously,
$$
h {\mathbf C} = h(b(\D)\mathbf{v}_{\mathbf{C}} + {\mathbf C}) - h b(\D)\mathbf{v}_{\mathbf{C}}.
\eqno(5.36)
$$
The summands in the right-hand side of (5.36) are orthogonal to each other in $L_2(\Omega;\C^m)$,
since the first one belongs to ${\rm Ker}\,X_0^*$, while the second one belongs to ${\rm Ran}\, X_0$.
Hence,
$$
\|h(b(\D)\mathbf{v}_{\mathbf{C}} + {\mathbf C})\|^2_{L_2(\Omega)} \leqslant
\|h {\mathbf C}\|^2_{L_2(\Omega)},\quad {\mathbf C}\in \C^m.
$$
By (5.28), this implies
$$
\langle g^0 {\mathbf C}, {\mathbf C}\rangle \leqslant \langle \overline{g} {\mathbf C}, {\mathbf C}\rangle,
\quad {\mathbf C}\in \C^m,
$$
which is equivalent to the upper estimate (5.35).

To prove the lower estimate, note that ${\mathfrak P}:= (h^*)^{-1}{\mathfrak M} \subset \NN_*={\rm Ker}\,X_0^*$.
We put
$$
Q \w = |\Omega|^{-1} (h^*)^{-1} \underline{g} \intop_\Omega h^{-1} \w \,d\x,\quad \w \in \H_*=L_2(\Omega;\C^m).
$$
It is easily checked that $Q \w \in {\mathfrak P}$ for $\w\in \H_*$;
$Q \w=\w$ for $\w\in {\mathfrak P}$; and
$$
(Q \w,\w)_{\H_*}= |\Omega|^{-1} \langle \underline{g} {\mathbf C}_{\w},{\mathbf C}_{\w}\rangle,
\quad {\mathbf C}_{\w}= \intop_{\Omega} h^{-1}\w\,d\x, \quad \w \in \H_*.
\eqno(5.37)
$$
It follows that  $Q$ is the orthogonal projection of $\H_*$ onto the subspace ${\mathfrak P}$.
We apply $Q$ to both sides of (5.36). Since $h b(\D)\mathbf{v}_{\mathbf{C}} \in {\rm Ran}\,X_0 = \NN_*^\perp$,
then $Q h {\mathbf C}= Q h (b(\D)\mathbf{v}_{\mathbf{C}} + {\mathbf C})$.
Hence, by (5.28),
$$
\begin{aligned}
&\| Q h {\mathbf C}\|^2_{\H_*}= \| Q h (b(\D)\mathbf{v}_{\mathbf{C}} + {\mathbf C})\|^2_{\H_*}
\\
&\leqslant \| h (b(\D)\mathbf{v}_{\mathbf{C}} + {\mathbf C})\|^2_{\H_*} =
|\Omega| \langle g^0 {\mathbf C}, {\mathbf C}\rangle,
\quad {\mathbf C}\in \C^m.
\end{aligned}
\eqno(5.38)
$$
From (5.37) with $\w=h {\mathbf C}$ it follows that
$$
\| Q h {\mathbf C}\|^2_{\H_*} = ( Q h {\mathbf C},h {\mathbf C})_{\H_*}=
|\Omega| \langle \underline{g} {\mathbf C},{\mathbf C}\rangle.
\eqno(5.39)
$$
Together with (5.38) this implies the lower estimate (5.35).

If $m=n$, we have $n_*=m=n$. Then ${\mathfrak P} \subset \NN_*$, ${\rm dim}\,{\mathfrak P}=m=n$,
and ${\rm dim}\,{\mathfrak N}_*=n_*=n$. Hence, ${\mathfrak P}= \NN_*$.
Since $h(b(\D)\v_{\mathbf C} + {\mathbf C}) \in \NN_*$, then the inequality in
(5.38) becomes an identity. This means that $g^0 = \underline{g}$. $\square$

\smallskip
Estimates of the form (5.35) are known in homogenization theory for specific DOs as the Voigt-Reuss bracketing.
The following estimates for the norms of the effective matrix and its inverse follow from (5.35):
$$
\left| g^{0} \right| \leqslant \Vert g \Vert _{L_{\infty}},
\quad \bigl| \left(g^{0}\right)^{-1} \bigr| \leqslant \Vert g^{-1}\Vert _{L_{\infty}}.
\eqno(5.40)
$$

Now we distinguish the cases where one of the inequalities in  (5.35) becomes an identity.
The following two statements are similar to Propositions 1.6 and 1.7 from [BSu1, Chapter~3].

\smallskip\noindent\textbf{Proposition 5.4.}
\textit{Let ${\mathbf g}_k(\x)$, $k=1,\dots,m,$ be the columns of the matrix $g(\x)$.
Then the identity $g^0 = \overline{g}$ is equivalent to the relations
$$
b(\D)^*  {\mathbf g}_k(\x)=0,\quad k=1,\dots,m.
\eqno(5.41)
$$
}

\smallskip\noindent\textbf{Proof.}
By (5.28), the identity $g^0 = \overline{g}$ is equivalent to the relation
$$
\| h (b(\D)\mathbf{v}_{\mathbf{C}} + {\mathbf C})\|^2_{\H_*} = \| h {\mathbf C}\|^2_{\H_*},\quad {\mathbf C}\in \C^m.
\eqno(5.42)
$$
As has already been mentioned, the terms in the right-hand side of (5.36) are orthogonal to each other. Therefore,
(5.42) is equivalent to the condition that $hb(\D)\mathbf{v}_{\mathbf{C}}=0$ for any ${\mathbf C}\in \C^m$.
By (5.25), this condition is satisfied if and only if $b(\D)^* g(\x) {\mathbf C}=0$ for any ${\mathbf C}\in \C^m$.
The latter is equivalent to (5.41). $\square$

\smallskip\noindent\textbf{Proposition 5.5.}
\textit{Let ${\mathbf l}_k(\x)$, $k=1,\dots,m,$ be the columns of the matrix $g(\x)^{-1}$.
The identity $g^0 = \underline{g}$ is equivalent to the representations
$$
{\mathbf l}_k(\x) = {\mathbf l}_k^0 + b(\D) \v_k(\x),\quad {\mathbf l}_k^0 \in \C^m,\quad \v_k \in \wt{H}^p(\Omega;\C^n);\quad k=1,\dots,m.
\eqno(5.43)
$$
}

\smallskip\noindent\textbf{Proof.}
 According to (5.38) and (5.39), the identity $g^0 = \underline{g}$ is equivalent to the condition that
$ h (b(\D)\mathbf{v}_{\mathbf{C}} + {\mathbf C}) \in {\mathfrak P}$ for any ${\mathbf C}\in \C^m$.
In other words, for any  ${\mathbf C}\in \C^m$ there exists a vector ${\mathbf C}_*\in \C^m$ such that
$ h (b(\D)\mathbf{v}_{\mathbf{C}} + {\mathbf C})= (h^*)^{-1}{\mathbf C}_*$.
The latter relation means that
$$
g(\x)^{-1} {\mathbf C}_* = b(\D)\mathbf{v}_{\mathbf{C}}(\x) + {\mathbf C}, \quad {\mathbf C}\in \C^m.
\eqno(5.44)
$$
Integrating over $\Omega$, we obtain $\underline{g}^{-1} {\mathbf C}_* = {\mathbf C}$.

Obviously, (5.44) holds for any ${\mathbf C} \in \C^m$ if and only if it holds for
${\mathbf C}=\underline{g}^{-1} \e_k$ (i.~e., ${\mathbf C}_* = \e_k$), $k=1,\dots,m$.
The latter condition is equivalent to representations (5.43) for the columns ${\mathbf l}_k(\x)$, $k =1,\dots,m$. $\square$

\smallskip\noindent\textbf{Remark 5.6.} From the proof of  Proposition 5.5 it is seen that, if $g^0 = \underline{g}$, then
 the matrix (5.30) is constant: $\wt{g}(\x)=g^0 = \underline{g}$.

\subsection{Estimates for the matrix-valued function $\Lambda$}

In what follows, we need estimates for the norms of $\Lambda$.

\smallskip\noindent\textbf{Lemma 5.7.}
\textit{Let $\v_j(\x)$, $j=1,\dots,m,$ be the columns of the matrix-valued function $\Lambda(\x)$ which is a $\Gamma$-periodic
solution of problem} (5.26). \textit{Then we have}
\begin{align}
&\| b(\D)\v_j \|_{L_2(\Omega)} \leqslant |\Omega|^{1/2} \|g\|_{L_\infty}^{1/2}\|g^{-1}\|_{L_\infty}^{1/2},\quad j=1,\dots,m,
\tag{5.45}
\\
&\|\v_j \|_{L_2(\Omega)} \leqslant \alpha_0^{-1/2} (2r_0)^{-p}  |\Omega|^{1/2} \|g\|_{L_\infty}^{1/2}\|g^{-1}\|_{L_\infty}^{1/2},
\quad j=1,\dots,m.
\tag{5.46}
\end{align}

\smallskip\noindent\textbf{Proof.}
Recall that the function $\v_j \in \wt{H}^p(\Omega;\C^n)$ satisfies the identity
$$
(g (b(\D)\v_j + \e_j), b(\D) \w)_{L_2(\Omega)} =0,\quad \w \in \wt{H}^p(\Omega;\C^n),
\eqno(5.47)
$$
and also the condition $\int_\Omega \v_j \,d\x=0$. From (5.47) it follows that
$$
\| h b(\D)\v_j \|_{L_2(\Omega)} \leqslant \| h \e_j \|_{L_2(\Omega)} \leqslant |\Omega|^{1/2} \|g\|_{L_\infty}^{1/2},\quad j=1,\dots,m,
$$
which implies (5.45).

To estimate $\| \v_j \|_{L_2(\Omega)}$, we use the Fourier series, (4.4), (4.11), and the condition $\int_\Omega \v_j \,d\x=0$.
Then we obtain
$$
\| b(\D)\v_j \|_{L_2(\Omega)}^2 \geqslant \alpha_0
\sum_{\mathbf{0} \ne {\mathbf s} \in \wt{\Gamma}} |\s|^{2p} |\wh{\v}_{j,{\mathbf s}}|^2 \geqslant \alpha_0 (2r_0)^{2p} \|\v_j \|^2_{L_2(\Omega)},\quad j=1,\dots,m,
\eqno(5.48)
$$
where $\wh{\v}_{j,{\mathbf s}}$, ${\mathbf s} \in \wt{\Gamma}$, are the Fourier coefficients of the function $\v_j$.
Relations (5.45) and (5.48) imply (5.46).
$\square$

\smallskip\noindent\textbf{Corollary 5.8.}
\textit{Suppose that the matrix-valued function $\Lambda(\x)$ is the $\Gamma$-periodic solution of the problem} (5.26).
\textit{Then we have}
\begin{align}
\|\Lambda \|_{L_2(\Omega)} &\leqslant |\Omega|^{1/2} C^{(1)}_\Lambda,
\tag{5.49}
\\
\| b(\D)\Lambda \|_{L_2(\Omega)} &\leqslant |\Omega|^{1/2} C^{(2)}_\Lambda,
\tag{5.50}
\\
\|\Lambda \|_{H^p(\Omega)} &\leqslant |\Omega|^{1/2} C_\Lambda,
\tag{5.51}
\end{align}
where
$$
\begin{aligned}
C^{(1)}_\Lambda &:= m^{1/2}\alpha_0^{-1/2} (2r_0)^{-p} \|g\|_{L_\infty}^{1/2}\|g^{-1}\|_{L_\infty}^{1/2},
\\
C^{(2)}_\Lambda &:= m^{1/2} \|g\|_{L_\infty}^{1/2}\|g^{-1}\|_{L_\infty}^{1/2},
\\
C_{\Lambda}&:=  C_\Lambda^{(2)}\alpha_0^{-1/2} \biggl(\sum_{|\beta|\leqslant p} (2r_0)^{-2(p-|\beta|)} \biggr)^{1/2}.
\end{aligned}
$$

\smallskip\noindent\textbf{Proof.}
Inequalities (5.49) and (5.50) follow directly from (5.46) and (5.45), respectively.

To check (5.51), we apply the Fourier series expansion.
Similarly to  (5.14), taking (4.4) into account, we have
$$
\begin{aligned}
&\| \D^\beta \Lambda \|^2_{L_2(\Omega)} \leqslant
(2r_0)^{-2(p-|\beta|)}  \sum_{\mathbf{0} \ne {\mathbf s}\in \wt{\Gamma}} |{\mathbf s}|^{2p} |\wh{\Lambda}_{\mathbf s}|^2
 \\
&\leqslant (2r_0)^{-2(p-|\beta|)}  \alpha_0^{-1}\sum_{{\mathbf 0} \ne {\mathbf s}\in \wt{\Gamma}} |b({\mathbf s}) \wh{\Lambda}_{\mathbf s}|^2
 \\
 &= (2r_0)^{-2(p-|\beta|)}  \alpha_0^{-1} \| b(\D)\Lambda \|^2_{L_2(\Omega)},\quad |\beta|\leqslant p.
\end{aligned}
$$
Hence,
$$
\| \Lambda \|^2_{H^p(\Omega)} \leqslant
 \alpha_0^{-1} \| b(\D)\Lambda \|^2_{L_2(\Omega)}
 \biggl(\sum_{|\beta|\leqslant p} (2r_0)^{-2(p-|\beta|)} \biggr).
$$
Together with (5.50) this implies (5.51).
 $\square$

\section{Approximation of the resolvent $\left(A(\k)+\varepsilon^{2p}I\right)^{-1}$}

\subsection{Approximation of the resolvent  $\left(A(\k)+\varepsilon^{2p}I\right)^{-1}$ in the operator norm in $L_2(\Omega;\C^n)$}

We apply Theorem 2.4 to the operator family $A(t,\bt)$.
The number $t^0$ is defined by (5.20) and does not depend on $\bt$.
Due to the presence of the projection  $P$ onto the subspace (5.7)
of constant vector-valued functions, (5.32) implies that
$$
t^{2p}S(\bt)P = S(\k)P = b(\k)^*g^0 b(\k)P = b(\D + \k)^*g^0 b(\D+ \k)P = A^0(\k)P.
\eqno(6.1)
$$
Hence, the operator under the norm sign in (2.27) takes the form
\hbox{$(A(\k) + \eps^{2p}I)^{-1} - (A^0(\k) + \eps^{2p}I)^{-1}P$}.
Now the constant $C_A$ depends on $\bt$.
According to Remark 2.5, this constant is a polynomial of the variables $\wt{C}$, $\|X_p(\bt)\|$, $\delta^{-1/2p}$,
and $c_*^{-1/2p}$ with positive coefficients depending only on $p$.
Relations (5.16), (5.18), (5.19), and (5.24) show that the constants
$\delta$, $\wt{C}$, and $c_*$ do not depend on $\bt$;
by (5.6), the norm $\|X_p(\bt)\|$ can be replaced by $\alpha_1^{1/2}\|g\|^{1/2}_{L_\infty}$.
Thus, after possible overstating, we can assume that the constant $C_A$ does not depend on $\bt$;
it depends only on
$d$, $p$, $\alpha_0$, $\alpha_1$, $\|g\|_{L_\infty}$, $\|g^{-1}\|_{L_\infty}$, and $r_0$.

Applying Theorem 2.4, we arrive at the inequality
$$
\begin{aligned}
\eps^{2p-1} \| (A(\k) + \eps^{2p}I)^{-1} - (A^0(\k) + \eps^{2p}I)^{-1}P\|_{L_2(\Omega)\to L_2(\Omega)}
\leqslant C_A,
\\
\eps >0,\quad |\k|\leqslant t^0.
\end{aligned}
\eqno(6.2)
$$

Let us show that the projection $P$ under the norm sign in (6.2)
can be replaced by the identity; 
one should only change the constant on the right.
Using the discrete Fourier transformation and taking (4.11) and (5.34) into account, we see that
$$
\begin{aligned}
&\eps^{2p-1} \|(A^0(\k) + \eps^{2p}I)^{-1}P^\perp \|_{L_2(\Omega) \to L_2(\Omega)}
\\
&=\eps^{2p-1} \sup_{{\mathbf 0} \ne {\mathbf s}\in \wt{\Gamma}} \left|(b({\mathbf s}+\k)^* g^0 b({\mathbf s}+\k) + \eps^{2p}\1_n)^{-1}\right|
\\
&\leqslant \eps^{2p-1} \sup_{{\mathbf 0} \ne {\mathbf s}\in \wt{\Gamma}} \left( c_*|{\mathbf s}+\k|^{2p} + \eps^{2p} \right)^{-1}
\leqslant c_*^{-1/2p} r_0^{-1}, \quad \eps>0,\quad |\k| \leqslant t^0.
\end{aligned}
$$
Combining this with (6.2), we obtain
$$
\begin{aligned}
\eps^{2p-1} \| (A(\k) + \eps^{2p}I)^{-1} - (A^0(\k) + \eps^{2p}I)^{-1}\|_{L_2(\Omega)\to L_2(\Omega)}
\\
\leqslant C_A + c_*^{-1/2p} r_0^{-1},
\quad
\eps >0,\quad |\k|\leqslant t^0.
\end{aligned}
\eqno(6.3)
$$

For $\k \in \mathrm{clos}\,\wt{\Omega} \setminus {\mathcal B}(t^0)$ estimates are trivial.
Each summand under the norm sign in (6.3) is estimated separately by using estimate (5.23)
for the first eigenvalue of $A(\k)$ and the similar estimate for the effective operator.
We have
$$
\begin{aligned}
\eps^{2p-1} \| (A(\k) + \eps^{2p}I)^{-1} \|_{L_2(\Omega)\to L_2(\Omega)}
\leqslant c_*^{-1/2p} (t^0)^{-1},
\\
 \eps^{2p-1} \| (A^0(\k) + \eps^{2p}I)^{-1} \|_{L_2(\Omega)\to L_2(\Omega)}
\leqslant c_*^{-1/2p} (t^0)^{-1},
\\
\eps>0,\quad \mathbf{k}\in\mathrm{clos}\,\widetilde{\Omega}\setminus  {\mathcal B}(t^0).
\end{aligned}
\eqno(6.4)
$$

Combining (6.3) and (6.4), and denoting
$$
\wt{C}_A = \max\{ C_A + c_*^{-1/2p} r_0^{-1}, 2 c_*^{-1/2p} (t^0)^{-1}\},
$$
we arrive at the following result.

\smallskip\noindent\textbf{Theorem 6.1.} \textit{Suppose that the operator $A(\k)=b(\D+\k)^* g b(\D+\k)$ is defined in Subsection}~4.3.
\textit{Let $A^0(\k)=b(\D+\k)^* g^0 b(\D+\k)$, where $g^0$ is the effective matrix defined in Subsection}~5.3.
\textit{Then we have}
$$
\begin{aligned}
\| (A(\k) + \eps^{2p}I)^{-1} - (A^0(\k) + \eps^{2p}I)^{-1}\|_{L_2(\Omega)\to L_2(\Omega)}
\leqslant \wt{C}_A \eps^{1-2p},
\\
\eps >0,\quad \k \in \mathrm{clos}\,\widetilde{\Omega}.
\end{aligned}
$$
\textit{The constant $\wt{C}_A$ depends only on $d$, $p$, $\alpha_0$, $\alpha_1$, $\|g\|_{L_\infty}$, $\|g^{-1}\|_{L_\infty}$, and
the parameters of the lattice~$\Gamma$.}

\subsection{Approximation of the resolvent  $\left(A(\k)+\varepsilon^{2p}I\right)^{-1}$ in the energy norm}

Now we apply Theorem 3.1 to the operator family $A(t,\bt)$. Similarly to (6.1), by (5.27),
$$
t^p Z(\bt)= \Lambda b(\k) P = \Lambda b(\D+\k)P.
\eqno(6.5)
$$
By (6.1) and (6.5), the operator under the norm sign in (3.1) takes the form
$$
{\mathcal E}(\k,\eps):= A(\k)^{1/2} \left( (A(\k)+ \eps^{2p}I)^{-1} - (I + \Lambda b(\D+\k) ) (A^0(\k)+ \eps^{2p}I)^{-1} P \right).
\eqno(6.6)
$$
According to Remark 3.4 and relations (5.6), (5.16), (5.18), (5.19), (5.24),
after possible overstating, we can assume that $\check{C}_A$ does not depend on $\bt$; it depends only on
$d$, $p$, $\alpha_0$, $\alpha_1$, $\|g\|_{L_\infty}$, $\|g^{-1}\|_{L_\infty}$, and $r_0$.

Applying Theorem~3.1, we arrive at the inequality
$$
\eps^{p-1} \left\| {\mathcal E}(\k,\eps) \right\|_{L_2(\Omega)\to L_2(\Omega)}
\leqslant \check{C}_A,
\quad \eps >0,\quad |\k|\leqslant t^0.
\eqno(6.7)
$$

Now we estimate the norm of the operator (6.6) for $\k \in \mathrm{clos}\,\wt{\Omega} \setminus {\mathcal B}(t^0)$.
The operator (6.6) can be represented as the sum of three terms:
$$
\begin{aligned}
{\mathcal E}(\k,\eps)&= A(\k)^{1/2}  (A(\k)+ \eps^{2p}I)^{-1}
 - A(\k)^{1/2}  (A^0(\k)+ \eps^{2p}I)^{-1}P
\\
&- A(\k)^{1/2}\Lambda P_m b(\D+\k)  (A^0(\k)+ \eps^{2p}I)^{-1}P,
\end{aligned}
\eqno(6.8)
$$
where $P_m$ is the orthogonal projection onto the subspace $\mathfrak M$
of constant vector-valued functions in $L_2(\Omega;\C^m)$.

The first term is estimated with the help of (5.23):
$$
\begin{aligned}
\eps^{p-1} \| A(\k)^{1/2}  (A(\k)+ \eps^{2p}I)^{-1}\|_{L_2(\Omega) \to L_2(\Omega)}
\leqslant c_*^{-1/2p} (t^0)^{-1},
\\ \k \in \mathrm{clos}\,\wt{\Omega} \setminus {\mathcal B}(t^0).
\end{aligned}
\eqno(6.9)
$$
Using the presence of the projection $P$ and relations (4.4), (5.34),
we obtain the following estimate for the second term:
$$
\begin{aligned}
&\| A(\k)^{1/2}  (A^0(\k)+ \eps^{2p}I)^{-1} P \|_{L_2(\Omega) \to L_2(\Omega)}
\\
&= \| h b(\D+\k) (A^0(\k)+ \eps^{2p}I)^{-1} P \|_{L_2(\Omega) \to L_2(\Omega)}
\\
&\leqslant \|g\|^{1/2}_{L_\infty} \left| b(\k) \left( b(\k)^* g^0 b(\k) + \eps^{2p} \1_n \right)^{-1} \right|
\\
&\leqslant \alpha_1^{1/2} \|g\|^{1/2}_{L_\infty} |\k|^p \left( c_* |\k|^{2p}+\eps^{2p}\right)^{-1}.
\end{aligned}
\eqno(6.10)
$$
Hence,
$$
\begin{aligned}
&\eps^{p-1} \| A(\k)^{1/2}  (A^0(\k)+ \eps^{2p}I)^{-1} P \|_{L_2(\Omega) \to L_2(\Omega)}
\\
&\leqslant \alpha_1^{1/2} \|g\|^{1/2}_{L_\infty}  c_*^{-1/2-1/2p} (t^0)^{-1},
  \quad \k \in \mathrm{clos}\,\wt{\Omega} \setminus {\mathcal B}(t^0).
\end{aligned}
\eqno(6.11)
$$

Similarly to (6.10), for the third term in the right-hand side of (6.8) we have:
$$
\begin{aligned}
&\| A(\k)^{1/2}\Lambda P_m b(\D+\k)  (A^0(\k)+ \eps^{2p}I)^{-1}P\|_{L_2(\Omega)\to L_2(\Omega)}
\\
&\leqslant
\| A(\k)^{1/2}\Lambda P_m\|_{L_2(\Omega)\to L_2(\Omega)}  \|b(\D+\k)  (A^0(\k)+ \eps^{2p}I)^{-1}P\|_{L_2(\Omega)\to L_2(\Omega)}
\\
&\leqslant
\alpha_1^{1/2} |\k|^{p} \left( c_* |\k|^{2p}+\eps^{2p} \right)^{-1}  \| A(\k)^{1/2}\Lambda P_m\|_{L_2(\Omega)\to L_2(\Omega)}.
\end{aligned}
\eqno(6.12)
$$
Obviously,
$$
\begin{aligned}
&\| A(\k)^{1/2}\Lambda P_m\|_{L_2(\Omega)\to L_2(\Omega)}=
\| h b(\D+\k) \Lambda P_m\|_{L_2(\Omega)\to L_2(\Omega)}
\\
&\leqslant |\Omega|^{-1/2} \|g\|^{1/2}_{L_\infty} \| b(\D+\k) \Lambda \|_{L_2(\Omega)}.
\end{aligned}
\eqno(6.13)
$$
Recall that the norm $\| b(\D) \Lambda \|_{L_2(\Omega)}$ satisfies (5.50).
By (4.4) and (5.49),
$$
\| b(\k)\Lambda \|_{L_2(\Omega)} \leqslant |\k|^p \alpha_1^{1/2}
|\Omega|^{1/2} C_\Lambda^{(1)} \leqslant r_1^p \alpha_1^{1/2} |\Omega|^{1/2} C_\Lambda^{(1)},
\quad \k \in \mathrm{clos}\,\wt{\Omega},
\eqno(6.14)
$$
where $2r_1= \hbox{\rm{diam}}\,\wt{\Omega}$.
Relations (5.50), (6.13), and (6.14) imply that
$$
\| A(\k)^{1/2}\Lambda P_m\|_{L_2(\Omega)\to L_2(\Omega)}
\leqslant \|g\|^{1/2}_{L_\infty} \left( C_\Lambda^{(2)} + r_1^p \alpha_1^{1/2} C_\Lambda^{(1)} \right),
\quad \k \in \mathrm{clos}\,\wt{\Omega}.
\eqno(6.15)
$$
By (6.12) and (6.15), we obtain the following estimate for the third term in (6.8) for $|\k|> t^0$:
$$
\begin{aligned}
& \eps^{p-1}\|A(\k)^{1/2}\Lambda  P_m b(\D+ \k) (A^0(\k)+ \eps^{2p}I)^{-1}P \|_{L_2(\Omega)\to L_2(\Omega)}
\leqslant C_{9},
\\
&\k \in \mathrm{clos}\,\wt{\Omega} \setminus {\mathcal B}(t^0),
\quad
C_{9}:=  \alpha_1^{1/2} \|g\|^{1/2}_{L_\infty} \left( C_\Lambda^{(2)} + r_1^p \alpha_1^{1/2} C_\Lambda^{(1)} \right)
c_*^{-1/2-1/2p} (t^0)^{-1}.
 \end{aligned}
\eqno(6.16)
$$

As a result, relations (6.8), (6.9), (6.11), and (6.16) imply the following estimate for the operator
(6.6) for $|\k|>t^0$:
$$
\begin{aligned}
&\eps^{p-1}\| {\mathcal E}(\k,\eps)\|_{L_2(\Omega) \to L_2(\Omega)}
\leqslant \wh{C}_{A},\quad \eps>0,\quad \k \in \mathrm{clos}\,\wt{\Omega} \setminus {\mathcal B}(t^0),
\\
&\wh{C}_{A}:=  c_*^{-1/2p}(t^0)^{-1} \left(1+ c_*^{-1/2} \alpha_1^{1/2} \|g\|^{1/2}_{L_\infty}\right) + C_{9}.
\end{aligned}
\eqno(6.17)
$$

Inequalities (6.7) and (6.17) lead to the following result.

\smallskip\noindent\textbf{Theorem 6.2.} \textit{Suppose that the operator $A(\k)=b(\D+\k)^* g b(\D+\k)$ is defined in Subsection} 4.3.
\textit{Let $A^0(\k)=b(\D+\k)^* g^0 b(\D+\k)$, where $g^0$ is the effective matrix defined in Subsection}~5.3.
\textit{Let $P$ be the orthogonal projection of $L_2(\Omega;\C^n)$ onto the subspace} (5.7).
\textit{Suppose that $\Lambda(\x)$ is a $\Gamma$-periodic solution of the problem}~(5.26). \textit{Then for
$\eps>0$ and $\k\in \mathrm{clos}\,\wt{\Omega}$ we have}
$$
\begin{aligned}
\bigl\|A(\k)^{1/2} &\left( (A(\k)+ \eps^{2p}I)^{-1} \right.
\\
&\left.- (I + \Lambda b(\D+\k)) (A^0(\k)+ \eps^{2p}I)^{-1}P  \right)
\bigr\|_{L_2(\Omega)\to L_2(\Omega)}
\leqslant {C}_A^\circ \eps^{1-p}.
\end{aligned}
$$
\textit{The constant $C_A^\circ = \max\{\check{C}_A, \wh{C}_A\}$ depends only on $m$, $d$, $p$, $\alpha_0$, $\alpha_1$, $\|g\|_{L_\infty}$, $\|g^{-1}\|_{L_\infty}$, and the parameters of the lattice $\Gamma$.}

\section{Approximation of the resolvent of the operator $A$}

\subsection{Approximation of the resolvent $(A+\eps^{2p}I)^{-1}$ in the operator norm in $L_2(\R^d;\C^n)$}

We return to the operator $A$ acting in $L_2(\R^d;\C^n)$. By (4.17),
$$
\left(A+\varepsilon^{2p}I\right)^{-1} = {\mathcal U}^{-1} \left(\int_{\wt{\Omega}}\oplus (A(\k)+\varepsilon^{2p}I)^{-1} \,d\k \right) {\mathcal U}.
\eqno(7.1)
$$
The operator $\left(A^0 +\varepsilon^{2p}I\right)^{-1}$ admits a similar expansion. Hence,
$$
\begin{aligned}
&\| \left(A+\varepsilon^{2p}I\right)^{-1} - \left(A^0+\varepsilon^{2p}I\right)^{-1}\|_{L_2(\R^d)\to L_2(\R^d)}
\\
&= \textrm{ess-}\!\sup_{\k \in \wt{\Omega}}  \| (A(\k)+\varepsilon^{2p}I)^{-1}- (A^0(\k)+\varepsilon^{2p}I)^{-1} \|_{L_2(\Omega) \to L_2(\Omega)}.
\end{aligned}
$$
Together with Theorem 6.1 this implies the following result.

\smallskip\noindent\textbf{Theorem 7.1.} \textit{Suppose that the operator $A=b(\D)^* g b(\D)$ is defined in Subsection} 4.1.
\textit{Let $A^0=b(\D)^* g^0 b(\D)$ be the effective operator. Then}
$$
\| (A + \eps^{2p}I)^{-1} - (A^0 + \eps^{2p}I)^{-1} \|_{L_2(\R^d)\to L_2(\R^d)}
\leqslant \wt{C}_A \eps^{1-2p},\quad \eps >0.
\eqno(7.2)
$$
\textit{The constant $\wt{C}_A$ depends only on $d$, $p$, $\alpha_0$, $\alpha_1$, $\|g\|_{L_\infty}$, $\|g^{-1}\|_{L_\infty}$, and the
parameters of the lattice~$\Gamma$.}

\subsection{Approximation of the resolvent $\left(A+\varepsilon^{2p}I\right)^{-1}$ in the energy norm}

Now we will obtain approximation of the
resolvent $\left(A+\varepsilon^{2p}I\right)^{-1}$ with corrector taken into account, using Theorem 6.2 and
expansion~(7.1).
Recall that, under the Gelfand transformation, the operator $b(\D)$ expands in the direct integral of the operators
$b(\D+\k)$; and the operator of multiplication by the periodic matrix-valued function $\Lambda(\x)$
turns into multiplication by the same function on the fibers of the direct integral $\mathcal K$ (see (4.16)).
We also need the operator
$
\Pi := \mathcal{U}^{-1} [P] \mathcal{U}
$
acting in $L_2(\R^d;\C^n)$. Here $[P]$  is the operator in $\mathcal K$ that acts on the fibers as $P$.
It is easily seen (see [BSu3, (6.8)]) that $\Pi$ is the pseudodifferential operator with the symbol $\chi_{\wt{\Omega}}(\bxi)$,
where  $\chi_{\wt{\Omega}}$ is the characteristic function of the set $\wt{\Omega}$, i.~e.,
$$
(\Pi \u)(\x) = (2\pi)^{-d/2} \intop_{\wt{\Omega}} e^{i \langle\x ,\bxi \rangle} \wh{\u}(\bxi)\,d\bxi.
\eqno(7.3)
$$
(Note that the operator $\Pi$ is smoothing.)

It follows that the operator
$$
{\mathcal E}(\eps):= A^{1/2} \left( (A+ \eps^{2p}I)^{-1} - (I + \Lambda b(\D)) (A^0 + \eps^{2p}I)^{-1}\Pi \right)
$$
expands in the direct integral of the operators ${\mathcal E}(\k,\eps)$ (see (6.6)).
Therefore, Theorem 6.2 implies that
$$
\eps^{p-1}\|{\mathcal E}(\eps) \|_{L_2(\R^d)\to L_2(\R^d)}
\leqslant C_A^\circ,\quad \eps>0.
\eqno(7.4)
$$

From (7.4) we deduce the following result.

\smallskip\noindent\textbf{Theorem 7.2.} \textit{Suppose that the operator $A=b(\D)^* g b(\D)$ is defined in Subsection} 4.1.
\textit{Let $A^0=b(\D)^* g^0 b(\D)$ be the effective operator. Suppose that $\Lambda(\x)$ is the $\Gamma$-periodic
solution of the problem} (5.26), \textit{and let $\wt{g}(\x)$ be the matrix-valued function} (5.30).
\textit{Let $\Pi$ be the operator} (7.3). \textit{Then for $\eps>0$ we have}
$$
\| (A + \eps^{2p}I)^{-1} - (I+ \Lambda b(\D)\Pi)(A^0 + \eps^{2p}I)^{-1} \|_{L_2(\R^d)\to L_2(\R^d)}
\leqslant {C}_A^{(1)} \eps^{1-2p},
\eqno(7.5)
$$
$$
\begin{aligned}
\| A^{1/2} &\left( (A + \eps^{2p}I)^{-1} - (I+ \Lambda b(\D)\Pi)(A^0 + \eps^{2p}I)^{-1} \right) \|_{L_2(\R^d)\to L_2(\R^d)}
\\
& \leqslant {C}_A^{(2)} \eps^{1-p},
\end{aligned}
\eqno(7.6)
$$
$$
\| g b(\D) (A + \eps^{2p}I)^{-1} - \wt{g} b(\D) (A^0 + \eps^{2p}I)^{-1}\Pi \|_{L_2(\R^d)\to L_2(\R^d)}
\leqslant {C}_A^{(3)} \eps^{1-p}.
\eqno(7.7)
$$
\textit{The constants ${C}^{(1)}_A$, ${C}^{(2)}_A$, ${C}^{(3)}_A$ depend only on $m$, $d$, $p$, $\alpha_0$, $\alpha_1$, $\|g\|_{L_\infty}$, $\|g^{-1}\|_{L_\infty}$, and the parameters of the lattice $\Gamma$.}

\smallskip\noindent\textbf{Proof.}
To check (7.5), we use (7.2) and estimate the operator
$$
\Lambda b(\D)\Pi (A^0 + \eps^{2p}I)^{-1}=
\Lambda\Pi_m b(\D)(A^0 + \eps^{2p}I)^{-1}\Pi.
$$
Here $\Pi_m$ is the pseudodifferential operator in
$L_2(\R^d;\C^m)$ with the symbol $\chi_{\wt{\Omega}}(\bxi)$.
The operator $[\Lambda] \Pi_m$ is unitarily equivalent to the direct integral
of the operators $[\Lambda] P_m$, whence
$$
\| [\Lambda] \Pi_m \|_{L_2(\R^d)\to L_2(\R^d)} =
\| [\Lambda] P_m \|_{L_2(\Omega)\to L_2(\Omega)}
\leqslant |\Omega|^{-1/2} \| \Lambda\|_{L_2(\Omega)} \leqslant   C^{(1)}_\Lambda.
\eqno(7.8)
$$
We took (5.49) into account. Next, using (4.4), (5.34), and (7.3), we obtain
$$
\begin{aligned}
& \eps^{2p-1} \| b(\D) (A^0 + \eps^{2p}I)^{-1} \Pi\|_{L_2(\R^d) \to L_2(\R^d)}
\\
&= \eps^{2p-1} \sup_{\bxi \in \wt{\Omega}} \left| b(\bxi) (b(\bxi)^* g^0 b(\bxi) + \eps^{2p}\1_n)^{-1}\right|
\\
&\leqslant \eps^{2p-1} \alpha_1^{1/2} \sup_{\bxi \in \wt{\Omega}} |\bxi|^{p} \left( c_* | \bxi|^{2p} + \eps^{2p} \right)^{-1}
\\
&\leqslant \alpha_1^{1/2} c_*^{-1/2p} \sup_{\bxi \in \wt{\Omega}} |\bxi|^{p-1}
\leqslant  \alpha_1^{1/2} c_*^{-1/2p} r_1^{p-1},\quad \eps >0.
\end{aligned}
\eqno(7.9)
$$
As a result, relations (7.2), (7.8), and (7.9) imply estimate (7.5) with the constant $C_A^{(1)}= \wt{C}_A +
C_\Lambda^{(1)}\alpha_1^{1/2} c_*^{-1/2p} r_1^{p-1}$.

Now we prove (7.6) with the help of (7.4). By (4.4), (5.34), and (7.3),
$$
\begin{aligned}
&\| A^{1/2}  (A^0+ \eps^{2p}I)^{-1} (I-\Pi) \|_{L_2(\R^d) \to L_2(\R^d)}
\\
&=
\| h b(\D)  (A^0+ \eps^{2p}I)^{-1} (I-\Pi) \|_{L_2(\R^d) \to L_2(\R^d)}
\\
&\leqslant \|g\|^{1/2}_{L_\infty}
\sup_{|\bxi| \geqslant r_0} \left|b(\bxi) \left( b(\bxi)^* g^0 b(\bxi) + \eps^{2p}\1_n \right)^{-1}\right|
\\
&\leqslant \alpha_1^{1/2} \|g\|^{1/2}_{L_\infty}
\sup_{|\bxi| \geqslant r_0} |\bxi|^p \left( c_* |\bxi|^{2p} + \eps^{2p} \right)^{-1}.
\end{aligned}
\eqno(7.10)
$$
Hence,
$$
\eps^{p-1} \| A^{1/2}  (A^0+ \eps^{2p}I)^{-1} (I-\Pi) \|_{L_2(\R^d) \to L_2(\R^d)}
\leqslant
\alpha_1^{1/2} \|g\|^{1/2}_{L_\infty} c_*^{-1/2-1/2p} r_0^{-1}.
$$
Together with (7.4) this implies (7.6) with the constant
${C}_A^{(2)} = C_A^\circ +\alpha_1^{1/2} \|g\|^{1/2}_{L_\infty} c_*^{-1/2-1/2p} r_0^{-1}$.

We proceed to the proof of (7.7). Denote
$$
{\mathcal G}(\eps):= g b(\D) \left( (A + \eps^{2p}I)^{-1} - (I+ \Lambda b(\D)) (A^0 + \eps^{2p}I)^{-1} \Pi \right).
\eqno(7.11)
$$
From (7.4) it follows that
$$
\eps^{p-1} \| {\mathcal G}(\eps) \|_{L_2(\R^d) \to L_2(\R^d)}
\leqslant C_A^\circ \|g\|^{1/2}_{L_\infty},\quad \eps >0.
\eqno(7.12)
$$
By (4.3), for any sufficiently smooth function $\u(\x)$ in $\R^d$ one has
$$
b(\D) (\Lambda \u) = (b(\D) \Lambda) \u + \sum_{|\alpha|=p} b_\alpha \sum_{\beta \leqslant \alpha: \, |\beta| \geqslant 1}
C_\alpha^\beta (\D^{\alpha-\beta} \Lambda) \D^\beta \u.
$$
Then, recalling the notation (5.30), we represent the operator (7.11) as
$$
{\mathcal G}(\eps)= g b(\D) (A + \eps^{2p}I)^{-1} - \wt{g} b(\D)(A^0 + \eps^{2p}I)^{-1} \Pi -
 \wt{\mathcal G}(\eps),
\eqno(7.13)
$$
where
$$
\wt{\mathcal G}(\eps):= g
\sum_{|\alpha|=p} b_\alpha \sum_{\beta \leqslant \alpha: \, |\beta| \geqslant 1}
C_\alpha^\beta (\D^{\alpha-\beta} \Lambda) \Pi_m \D^\beta b(\D) (A^0 + \eps^{2p}I)^{-1} \Pi.
$$
Similarly to (7.8),
$$
\| [\D^{\alpha-\beta} \Lambda] \Pi_m \|_{L_2(\R^d)\to L_2(\R^d)}
\leqslant |\Omega|^{-1/2} \| \D^{\alpha-\beta} \Lambda\|_{L_2(\Omega)} \leqslant   C_\Lambda.
\eqno(7.14)
$$
In the last passage, (5.51) was used.
By (4.5) and (7.14),
$$
\begin{aligned}
\| \wt{\mathcal G}(\eps) \|_{L_2\to L_2}
&\leqslant \|g\|_{L_\infty} \alpha_1^{1/2} C_\Lambda
\\
&\times
\sum_{|\alpha|=p} \sum_{\beta \leqslant \alpha: \, |\beta| \geqslant 1}
C_\alpha^\beta \| \D^\beta b(\D) (A^0 + \eps^{2p}I)^{-1} \Pi\|_{L_2\to L_2}.
\end{aligned}
\eqno(7.15)
 $$
Similarly to (7.9), from (4.4), (5.34), and (7.3) it follows that
$$
 \| \D^\beta b(\D) (A^0 + \eps^{2p}I)^{-1} \Pi\|_{L_2(\R^d) \to L_2(\R^d)}
\leqslant \alpha_1^{1/2} \sup_{\bxi \in \wt{\Omega}} |\bxi|^{p+|\beta|} \left( c_* | \bxi|^{2p} + \eps^{2p} \right)^{-1},
$$
whence
$$
\eps^{p-1} \| \D^\beta b(\D) (A^0 + \eps^{2p}I)^{-1} \Pi\|_{L_2(\R^d) \to L_2(\R^d)}
\leqslant  \alpha_1^{1/2} c_*^{-1/2-1/2p} r_1^{|\beta|-1}
\eqno(7.16)
$$
for $|\beta| \geqslant 1$. By (7.15) and (7.16),
$$
\eps^{p-1} \| \wt{\mathcal G}(\eps) \|_{L_2(\R^d) \to L_2(\R^d)}
\leqslant C_{10},\quad \eps>0,
\eqno(7.17)
$$
where
$$
C_{10}:=  \|g\|_{L_\infty} \alpha_1 C_\Lambda c_*^{-1/2-1/2p} \biggl( \sum_{|\alpha|=p} \sum_{\beta \leqslant \alpha: \, |\beta| \geqslant 1}
C_\alpha^\beta r_1^{|\beta|-1}\biggr).
$$
As a result, relations (7.12), (7.13), and (7.17) imply the required inequality (7.7) with the constant
$C_A^{(3)} = C_A^\circ \|g\|_{L_\infty}^{1/2} + C_{10}$. $\square$

\smallskip
Now we distinguish the special cases. If $g^0=\overline{g}$, then conditions (5.41) are satisfied,
whence the solution $\Lambda$ of the problem (5.26) is equal to zero.
Therefore, (7.6) simplifies, and we arrive at the following statement.

\smallskip\noindent\textbf{Proposition 7.3.} \textit{If $g^0=\overline{g}$} (\textit{i.~e., conditions} (5.41)
\textit{are satisfied}), \textit{then}
$$
\| A^{1/2} \left( (A + \eps^{2p}I)^{-1} - (A^0 + \eps^{2p}I)^{-1} \right) \|_{L_2(\R^d)\to L_2(\R^d)}
\leqslant {C}_A^{(2)} \eps^{1-p},\quad \eps>0.
\eqno(7.18)
$$

\smallskip
Now we consider the case where $g^0=\underline{g}$.

\smallskip\noindent\textbf{Proposition 7.4.} \textit{If $g^0=\underline{g}$} (\textit{i.~e., representations} (5.43)
\textit{are satisfied}), \textit{then}
$$
\| gb(\D) (A + \eps^{2p}I)^{-1} - g^0 b(\D)(A^0 + \eps^{2p}I)^{-1}  \|_{L_2(\R^d)\to L_2(\R^d)}
\leqslant \wt{C}_A^{(3)} \eps^{1-p},\quad \eps>0.
\eqno(7.19)
$$
\textit{The constant $\wt{C}_A^{(3)}$ depends only on $m$, $d$, $p$, $\alpha_0$, $\alpha_1$, $\|g\|_{L_\infty}$, $\|g^{-1}\|_{L_\infty}$, and
 the parameters of the lattice $\Gamma$.}

\smallskip\noindent\textbf{Proof.}
By Remark 5.6, in the case under consideration we have $\wt{g}(\x)=g^0=\underline{g}$.
Then the inequality (7.7) takes the form
$$
\| g b(\D) (A + \eps^{2p}I)^{-1} - {g}^0 b(\D) (A^0 + \eps^{2p}I)^{-1}\Pi \|_{L_2(\R^d)\to L_2(\R^d)}
\leqslant {C}_A^{(3)} \eps^{1-p}
\eqno(7.20)
$$
for $\eps >0$. Similarly to (7.10), by (5.40), we have
$$
\begin{aligned}
&\| {g}^0 b(\D) (A^0 + \eps^{2p}I)^{-1}(I-\Pi) \|_{L_2(\R^d)\to L_2(\R^d)}
\\
&\leqslant  \|g\|_{L_\infty} \alpha_1^{1/2} \sup_{|\bxi| \geqslant r_0} |\bxi|^p \left( c_* |\bxi|^{2p} + \eps^{2p} \right)^{-1}.
\end{aligned}
$$
Hence,
$$
\eps^{p-1} \| {g}^0 b(\D) (A^0 + \eps^{2p}I)^{-1}(I-\Pi) \|_{L_2(\R^d)\to L_2(\R^d)}
\leqslant  \alpha_1^{1/2} \|g\|_{L_\infty}  c_*^{-1/2 -1/2p} r_0^{-1}.
\eqno(7.21)
$$
Relations (7.20) and (7.21) imply (7.19) with the constant
$\wt{C}_A^{(3)} = {C}_A^{(3)} + \alpha_1^{1/2} \|g\|_{L_\infty}  c_*^{-1/2 -1/2p} r_0^{-1}$. $\square$

\subsection{Approximation of the resolvent  $\left(A-\zeta\varepsilon^{2p}I\right)^{-1}$ for $\zeta \in \C \setminus \R_+$}

Now we carry the results of Theorems 7.1 and 7.2 over to the case
of the resolvent at arbitrary regular point from $\C \setminus \R_+$.
For this, we apply appropriate identities for the resolvents (cf. [Su]).

Let $\zeta \in \C \setminus \R_+$. We put $\zeta = |\zeta|e^{i\varphi}$, $0 <\varphi < 2\pi$, and denote
$$
c(\varphi) =
\begin{cases}
|\sin \varphi|^{-1}, & \varphi \in (0,\pi/2) \cup (3\pi/2, 2\pi) \\
1, & \varphi \in [\pi/2, 3\pi/2]
\end{cases}.
\eqno(7.22)
$$

\smallskip\noindent\textbf{Theorem 7.5.} \textit{Suppose that the assumptions of Theorem} 7.1
\textit{are satisfied. Let $\zeta = |\zeta|e^{i\varphi} \in \C \setminus \R_+$, and let $c(\varphi)$ be defined by} (7.22).
\textit{Then for $\eps>0$ we have}
$$
\| (A - \zeta \eps^{2p}I)^{-1} - (A^0 -\zeta \eps^{2p}I)^{-1} \|_{L_2(\R^d)\to L_2(\R^d)}
\leqslant {\mathcal C}_1 c(\varphi)^2 \eps^{1-2p} |\zeta|^{1/2p -1}.
\eqno(7.23)
$$
\textit{The constant ${\mathcal C}_1=4 \wt{C}_A$ depends only on $d$, $p$, $\alpha_0$, $\alpha_1$, $\|g\|_{L_\infty}$, $\|g^{-1}\|_{L_\infty}$,
and the parameters of the lattice~$\Gamma$.}

\smallskip\noindent\textbf{Proof.}
Let $\wh{\zeta} = e^{i\varphi}$, $0 <\varphi < 2\pi$. We rely on the identity
$$
\begin{aligned}
 &(A - \wh{\zeta} \eps^{2p}I)^{-1} - (A^0 - \wh{\zeta} \eps^{2p}I)^{-1}
 =
 (A + \eps^{2p}I) (A - \wh{\zeta} \eps^{2p}I)^{-1}
 \\
 &\times \left((A + \eps^{2p}I)^{-1} - (A^0 + \eps^{2p}I)^{-1}\right)
 (A^0 + \eps^{2p}I) (A^0 - \wh{\zeta} \eps^{2p}I)^{-1}.
\end{aligned}
\eqno(7.24)
$$
Since the spectrum of $A$ is contained in $\R_+$, then
$$
\begin{aligned}
&\|  (A + \eps^{2p}I) (A - \wh{\zeta} \eps^{2p}I)^{-1}\|_{L_2(\R^d)\to L_2(\R^d)}
\leqslant \sup_{x \geqslant 0} (x+\eps^{2p}) \bigl|x - \wh{\zeta} \eps^{2p} \bigr|^{-1}
\\
&= \sup_{y \geqslant 0} (y + 1) \bigl|y - \wh{\zeta}\bigr|^{-1} \leqslant 2 c(\varphi).
\end{aligned}
\eqno(7.25)
$$
Similarly,
$$
\|  (A^0 + \eps^{2p}I) (A^0 - \wh{\zeta} \eps^{2p}I)^{-1}\|_{L_2(\R^d)\to L_2(\R^d)}
 \leqslant 2 c(\varphi).
\eqno(7.26)
$$
From (7.2), (7.24)--(7.26) it follows that
$$
\|(A - \wh{\zeta} \eps^{2p}I)^{-1} - (A^0 - \wh{\zeta} \eps^{2p} I)^{-1} \|_{L_2(\R^d)\to L_2(\R^d)}
\leqslant 4 c(\varphi)^2 \wt{C}_A \eps^{1-2p},\quad \eps>0.
\eqno(7.27)
$$
Replacing $\eps$ by $\eps|\zeta|^{1/2p}$ in (7.27), we arrive at the required inequality (7.23). $\square$

\smallskip\noindent\textbf{Theorem 7.6.} \textit{Suppose that the assumptions of Theorem} 7.2
\textit{are satisfied. Let $\zeta = |\zeta|e^{i\varphi} \in \C \setminus \R_+$, and let $c(\varphi)$ be given by} (7.22).
\textit{Then for $\eps>0$ we have}
$$
\begin{aligned}
&\| (A - \zeta \eps^{2p}I)^{-1} - (I + \Lambda b(\D)\Pi)(A^0 -\zeta \eps^{2p}I)^{-1} \|_{L_2(\R^d)\to L_2(\R^d)}
\\
&\leqslant {\mathcal C}_2 c(\varphi)^2 \eps^{1-2p} |\zeta|^{1/2p -1},
\end{aligned}
\eqno(7.28)
$$
$$
\begin{aligned}
&\| A^{1/2}\left((A - \zeta \eps^{2p}I)^{-1} - (I + \Lambda b(\D)\Pi)(A^0 -\zeta \eps^{2p}I)^{-1}\right) \|_{L_2(\R^d)\to L_2(\R^d)}
\\
&\leqslant {\mathcal C}_3 c(\varphi)^2 \eps^{1-p} |\zeta|^{1/2p -1/2},
\end{aligned}
\eqno(7.29)
$$
$$
\begin{aligned}
&\| g b(\D)(A - \zeta \eps^{2p}I)^{-1} - \wt{g} b(\D) (A^0 -\zeta \eps^{2p}I)^{-1} \Pi \|_{L_2(\R^d)\to L_2(\R^d)}
\\
&\leqslant {\mathcal C}_4 c(\varphi)^2 \eps^{1-p} |\zeta|^{1/2p -1/2}.
\end{aligned}
\eqno(7.30)
$$
\textit{The constants ${\mathcal C}_2$, ${\mathcal C}_3$, ${\mathcal C}_4$
depend only on $m$, $d$, $p$, $\alpha_0$, $\alpha_1$, $\|g\|_{L_\infty}$, $\|g^{-1}\|_{L_\infty}$, and the parameters of the lattice $\Gamma$.}

\smallskip\noindent\textbf{Proof.}
Let $\wh{\zeta} = e^{i\varphi}$, $0 <\varphi < 2\pi$. We use the identity
$$
\begin{aligned}
 &(A - \wh{\zeta} \eps^{2p}I)^{-1} - (I + \Lambda b(\D)\Pi)(A^0 - \wh{\zeta} \eps^{2p}I)^{-1}
 =
 (A + \eps^{2p}I) (A - \wh{\zeta} \eps^{2p}I)^{-1}
 \\
 &\times \left((A + \eps^{2p}I)^{-1} - (I + \Lambda b(\D)\Pi)(A^0 + \eps^{2p}I)^{-1}\right)
 (A^0 + \eps^{2p}I) (A^0 - \wh{\zeta} \eps^{2p}I)^{-1}
 \\
 &+  \eps^{2p} (\wh{\zeta}+1) (A - \wh{\zeta} \eps^{2p}I)^{-1} \Lambda b(\D)\Pi (A^0 - \wh{\zeta} \eps^{2p}I)^{-1}.
\end{aligned}
\eqno(7.31)
$$
Denote the consecutive terms in the right-hand side of (7.31) by ${\mathcal J}_1(\wh{\zeta},\eps)$ and ${\mathcal J}_2(\wh{\zeta},\eps)$.

From (7.5), (7.25), and (7.26) it follows that
$$
\eps^{2p-1} \| {\mathcal J}_1(\wh{\zeta},\eps)\|_{L_2(\R^d) \to L_2(\R^d)}
\leqslant 4 c(\varphi)^2 C_A^{(1)}.
\eqno(7.32)
$$
The second term satisfies
$$
\begin{aligned}
&\| {\mathcal J}_2(\wh{\zeta},\eps)\|_{L_2(\R^d) \to L_2(\R^d)}
\leqslant 2 \eps^{2p} \|(A - \wh{\zeta} \eps^{2p}I)^{-1} \|_{L_2\to L_2}
\\
&\times \|\Lambda \Pi_m b(\D) (A^0 + \eps^{2p}I)^{-1}\Pi\|_{L_2\to L_2}
\|(A^0 + \eps^{2p}I)(A^0 - \wh{\zeta} \eps^{2p}I)^{-1}\|_{L_2\to L_2}.
\end{aligned}
\eqno(7.33)
$$
Note that the norm of the resolvent $(A - \wh{\zeta} \eps^{2p}I)^{-1}$
does not exceed the inverse distance from the point $\wh{\zeta} \eps^{2p}$ to $\R_+$:
$$
\|(A - \wh{\zeta} \eps^{2p}I)^{-1} \|_{L_2(\R^d) \to L_2(\R^d)} \leqslant \eps^{-2p} c(\varphi).
\eqno(7.34)
$$
Combining (7.8), (7.9), (7.26), (7.33), and (7.34), we obtain
$$
\eps^{2p-1}\| {\mathcal J}_2(\wh{\zeta},\eps)\|_{L_2(\R^d) \to L_2(\R^d)}
\leqslant 4 c(\varphi)^2 C_\Lambda^{(1)} \alpha_1^{1/2} c_*^{-1/2p} r_1^{p-1}.
\eqno(7.35)
$$
Relations (7.31), (7.32), and (7.35) imply that
$$
\begin{aligned}
&\| (A - \wh{\zeta} \eps^{2p}I)^{-1} - (I + \Lambda b(\D)\Pi)(A^0 -\wh{\zeta} \eps^{2p}I)^{-1} \|_{L_2(\R^d)\to L_2(\R^d)}
\\
&\leqslant {\mathcal C}_2  c(\varphi)^2 \eps^{1-2p},\quad \eps>0,
\end{aligned}
\eqno(7.36)
$$
where ${\mathcal C}_2 = 4 {C}_A^{(1)} + 4 C_\Lambda^{(1)} \alpha_1^{1/2} c_*^{-1/2p} r_1^{p-1}$.
Replacing $\eps$ by $\eps |\zeta|^{1/2p}$ in (7.36), we arrive at (7.28).

To check (7.29), we apply the operator $A^{1/2}$ to both sides of (7.31):
$$
\begin{aligned}
A^{1/2} \left((A - \wh{\zeta} \eps^{2p}I)^{-1} - (I + \Lambda b(\D)\Pi)(A^0 - \wh{\zeta} \eps^{2p}I)^{-1}\right)
= {\mathcal T}_1(\wh{\zeta},\eps) + {\mathcal T}_2(\wh{\zeta},\eps),
\end{aligned}
\eqno(7.37)
$$
where
$$
\begin{aligned}
{\mathcal T}_1(\wh{\zeta},\eps) &= (A + \eps^{2p}I)(A - \wh{\zeta} \eps^{2p}I)^{-1}
\\
&\times
A^{1/2}\left((A + \eps^{2p}I)^{-1} - (I + \Lambda b(\D)\Pi)(A^0 + \eps^{2p}I)^{-1}\right)
\\
&\times
(A^0 + \eps^{2p}I)(A^0 - \wh{\zeta} \eps^{2p}I)^{-1},
\\
{\mathcal T}_2(\wh{\zeta},\eps)&=
\eps^{2p} (\wh{\zeta}+1) A^{1/2} (A - \wh{\zeta} \eps^{2p}I)^{-1}
\Lambda b(\D)\Pi (A^0 - \wh{\zeta} \eps^{2p}I)^{-1}.
\end{aligned}
$$
Estimate for the first term is deduced from (7.6), (7.25), and (7.26):
$$
\eps^{p-1} \|{\mathcal T}_1(\wh{\zeta},\eps) \|_{L_2(\R^d)\to L_2(\R^d)}
\leqslant 4 C_A^{(2)} c(\varphi)^2,\quad \eps >0.
\eqno(7.38)
$$
The second term satisfies
$$
\begin{aligned}
 & \|{\mathcal T}_2(\wh{\zeta},\eps) \|_{L_2(\R^d)\to L_2(\R^d)}
\leqslant 2 \eps^{2p}  \|A^{1/2} (A - \wh{\zeta} \eps^{2p}I)^{-1}\|_{L_2\to L_2}
\\
&\times \|\Lambda \Pi_m b(\D) (A^0 + \eps^{2p}I)^{-1}\Pi \|_{L_2\to L_2}
\|(A^0 + \eps^{2p}I) (A^0 - \wh{\zeta} \eps^{2p}I)^{-1} \|_{L_2\to L_2}.
\end{aligned}
\eqno(7.39)
$$
From (4.4), (5.34), (7.3), and (7.8) it follows that
$$
\|\Lambda \Pi_m b(\D) (A^0 + \eps^{2p}I)^{-1}\Pi \|_{L_2\to L_2}
\leqslant C_\Lambda^{(1)} \alpha_1^{1/2}
\sup_{\bxi \in \wt{\Omega}} |\bxi|^{p} \left( c_* |\bxi|^{2p} + \eps^{2p} \right)^{-1}.
\eqno(7.40)
$$
By (7.25),
$$
\|A^{1/2} (A - \wh{\zeta} \eps^{2p}I)^{-1}\|_{L_2\to L_2}
\leqslant 2 c(\varphi) \|A^{1/2} (A + \eps^{2p}I)^{-1}\|_{L_2\to L_2}
\leqslant 2 c(\varphi) \eps^{-p}.
\eqno(7.41)
$$
Now, relations (7.26) and (7.39)--(7.41) imply that
$$
\begin{aligned}
 & \eps^{p-1} \|{\mathcal T}_2(\wh{\zeta},\eps) \|_{L_2(\R^d)\to L_2(\R^d)}
\leqslant 8  c(\varphi)^2 C_\Lambda^{(1)} \alpha_1^{1/2} \eps^{2p-1}
\sup_{\bxi \in \wt{\Omega}} |\bxi|^{p} \left( c_* |\bxi|^{2p} + \eps^{2p} \right)^{-1}
\\
&\leqslant
8  c(\varphi)^2 C_\Lambda^{(1)} \alpha_1^{1/2} c_*^{-1/2p} r_1^{p-1},\quad \eps>0.
\end{aligned}
\eqno(7.42)
$$
As a result, relations (7.37), (7.38), and (7.42) yield
$$
\begin{aligned}
&\| A^{1/2} \left((A - \wh{\zeta} \eps^{2p}I)^{-1} - (I + \Lambda b(\D)\Pi)(A^0 - \wh{\zeta} \eps^{2p}I)^{-1}\right)\|_{L_2(\R^d)\to L_2(\R^d)}
\\
&\leqslant {\mathcal C}_3 c(\varphi)^2 \eps^{1-p},\quad \eps >0,
\end{aligned}
\eqno(7.43)
$$
with the constant
${\mathcal C}_3  = 4{C}_A^{(2)} + 8  C_\Lambda^{(1)} \alpha_1^{1/2} c_*^{-1/2p} r_1^{p-1}$.
Replacing $\eps$ by $\eps|\zeta|^{1/2p}$ in (7.43), we arrive at the required inequality (7.29).

It remains to check (7.30). We use the identity
$$
\begin{aligned}
& g b(\D)(A - \wh{\zeta} \eps^{2p}I)^{-1} - \wt{g} b(\D) (A^0 - \wh{\zeta} \eps^{2p}I)^{-1} \Pi
\\
&= \left(g b(\D)(A + \eps^{2p}I)^{-1} - \wt{g} b(\D) (A^0 + \eps^{2p}I)^{-1} \Pi\right)
(A^0 + \eps^{2p}I) (A^0 - \wh{\zeta} \eps^{2p}I)^{-1}
\\
& +
(\wh{\zeta}+1) \eps^{2p} g b(\D) (A + \eps^{2p}I)^{-1}\left((A - \wh{\zeta} \eps^{2p}I)^{-1}- (A^0 - \wh{\zeta} \eps^{2p}I)^{-1}\right).
\end{aligned}
\eqno(7.44)
$$
Denote the consecutive summands in the right-hand side of (7.44) by ${\mathcal L}_1(\wh{\zeta},\eps)$ and
${\mathcal L}_2(\wh{\zeta},\eps)$. From (7.7) and (7.26) it follows that
$$
\|{\mathcal L}_1(\wh{\zeta},\eps)\|_{L_2(\R^d) \to L_2(\R^d)}
\leqslant 2 C_A^{(3)} c(\varphi) \eps^{1-p}.
\eqno(7.45)
$$
The second term is estimated with the help of (7.27):
$$
\begin{aligned}
&\|{\mathcal L}_2(\wh{\zeta},\eps)\|_{L_2(\R^d) \to L_2(\R^d)}
\leqslant 8\eps c(\varphi)^2 \wt{C}_A \|g\|^{1/2}_{L_\infty} \|A^{1/2} (A + \eps^{2p}I)^{-1}\|_{L_2\to L_2}
\\
&\leqslant 8\eps^{1-p} c(\varphi)^2 \wt{C}_A \|g\|^{1/2}_{L_\infty}.
\end{aligned}
\eqno(7.46)
$$
Now relations (7.44)--(7.46) imply that
$$
\| g b(\D)(A - \wh{\zeta} \eps^{2p}I)^{-1} - \wt{g} b(\D) (A^0 - \wh{\zeta} \eps^{2p}I)^{-1} \Pi \|_{L_2(\R^d)\to L_2(\R^d)}
\leqslant  {\mathcal C}_4 c(\varphi)^2 \eps^{1-p}
\eqno(7.47)
$$
for $\eps>0$, where ${\mathcal C}_4 = 2 C_A^{(3)} + 8 \wt{C}_A \|g\|^{1/2}_{L_\infty}$.
Replacing $\eps$ by $\eps|\zeta|^{1/2p}$ in (7.47), we arrive at (7.30). $\square$

\subsection{Special cases}
Now we prove analogs of Propositions 7.3 and 7.4 for the resolvent $(A-\zeta \eps^{2p}I)^{-1}$.
If $g^0=\overline{g}$, then $\Lambda=0$, and (7.29) leads to the following statement.

\smallskip\noindent\textbf{Proposition 7.7.} \textit{Suppose that the assumptions of Theorem} 7.5
\textit{are satisfied. If $g^0=\overline{g}$} (\textit{i.~e., conditions} (5.41)
\textit{are satisfied}), \textit{then}
$$
\begin{aligned}
&\| A^{1/2} \left( (A - \zeta \eps^{2p}I)^{-1} - (A^0 - \zeta \eps^{2p}I)^{-1} \right) \|_{L_2(\R^d)\to L_2(\R^d)}
\\
&\leqslant {\mathcal C}_3 c(\varphi)^2 \eps^{1-p} |\zeta|^{1/2p - 1/2},\quad \eps>0.
\end{aligned}
$$

\smallskip
Now we consider the case where $g^0=\underline{g}$.

\smallskip\noindent\textbf{Proposition 7.8.} \textit{Suppose that the assumptions of Theorem} 7.5
\textit{are satisfied. If $g^0=\underline{g}$} (\textit{i.~e., representations} (5.43)
\textit{are satisfied}), \textit{then}
$$
\begin{aligned}
&\| gb(\D) (A -\zeta \eps^{2p}I)^{-1} - g^0 b(\D)(A^0 -\zeta \eps^{2p}I)^{-1}  \|_{L_2(\R^d)\to L_2(\R^d)}
\\
&\leqslant {\mathcal C}_4^\circ c(\varphi)^2 \eps^{1-p} |\zeta|^{1/2p - 1/2},\quad \eps>0.
\end{aligned}
\eqno(7.48)
$$
\textit{The constant ${\mathcal C}_4^\circ$ depends only on $m$, $d$, $p$, $\alpha_0$, $\alpha_1$, $\|g\|_{L_\infty}$,
$\|g^{-1}\|_{L_\infty}$, and the parameters of the lattice $\Gamma$.}

\smallskip\noindent\textbf{Proof.}
By Remark 5.6, in the case under consideration we have $\wt{g}(\x)=g^0=\underline{g}$.
First, we consider the resolvent at the point $\wh{\zeta}\eps^{2p}$, where $\wh{\zeta}= e^{i\varphi}$.
The following analog of the identity (7.44) is true:
$$
gb(\D) (A -\wh{\zeta} \eps^{2p}I)^{-1} - g^0 b(\D)(A^0 - \wh{\zeta} \eps^{2p}I)^{-1} =
{\mathcal L}^\circ_1(\wh{\zeta},\eps) + {\mathcal L}_2(\wh{\zeta},\eps)
$$
where
$$
\begin{aligned}
{\mathcal L}^\circ_1(\wh{\zeta},\eps) &=
\left(gb(\D) (A + \eps^{2p}I)^{-1} - g^0 b(\D)(A^0 + \eps^{2p}I)^{-1}\right)
\\
&\times (A^0+\eps^{2p}I) (A^0 - \wh{\zeta}\eps^{2p}I )^{-1},
\end{aligned}
$$
and the second term ${\mathcal L}_2(\wh{\zeta},\eps)$ is the same as in (7.44).

From (7.19) and (7.26) it follows that
$$
\| {\mathcal L}^\circ_1(\wh{\zeta},\eps)   \|_{L_2(\R^d)\to L_2(\R^d)} \leqslant 2 c(\varphi) \wt{C}_A^{(3)}\eps^{1-p},\quad \eps >0.
$$
Combining this with (7.46), we obtain
$$
\| gb(\D) (A -\wh{\zeta} \eps^{2p}I)^{-1} - g^0 b(\D)(A^0 -\wh{\zeta} \eps^{2p}I)^{-1}  \|_{L_2(\R^d)\to L_2(\R^d)}
\leqslant {\mathcal C}_4^\circ c(\varphi)^2 \eps^{1-p}
\eqno(7.49)
$$
for $\eps >0$, where ${\mathcal C}_4^\circ = 2 \wt{C}_A^{(3)} + 8 \wt{C}_A \|g\|^{1/2}_{L_\infty}$.
Replacing $\eps$ by $\eps|\zeta|^{1/2p}$ in (7.49), we arrive at (7.48). $\square$

\subsection{Removal of the smoothing operator}

It turns out that, under some additional assumptions on the matrix-valued function $\Lambda(\x)$,
it is possible to remove the smoothing operator $\Pi$ in approximations (7.28)--(7.30).
However, the order of estimates for the terms containing $I-\Pi$
is different from the order of estimates (7.28)--(7.30);
see Proposition 7.12 below.

\smallskip\noindent\textbf{Condition 7.9.} \textit{Suppose that the $\Gamma$-periodic
solution $\Lambda\in \wt{H}^p(\Omega)$ of the problem} (5.26)
\textit{is bounded and is a multiplier from $H^p(\R^d;\C^m)$ to $H^p(\R^d;\C^n)$}:
$$
\Lambda \in L_\infty(\R^d) \cap M(H^p(\R^d;\C^m)\to H^p(\R^d;\C^n)).
$$

\smallskip
Since  $\Lambda$ is periodic, Condition 7.9 is equivalent to the relation
$\Lambda \in L_\infty(\Omega) \cap M(H^p(\Omega;\C^m)\to H^p(\Omega;\C^n)).$
The norm of the operator  $[\Lambda]$ of multiplication by the matrix-valued function $\Lambda(\x)$ is denoted by
$$
M_\Lambda := \|[\Lambda]\|_{H^p(\R^d) \to H^p(\R^d)}.
\eqno(7.50)
$$

Description of the spaces of multipliers in the Sobolev classes can be found in the book [MSh].
The following statement gives some sufficient conditions ensuring that Condition 7.9 is satisfied.

\smallskip\noindent\textbf{Proposition 7.10.}
\textit{Suppose that at least one of the following two assumptions is satisfied}:

\noindent
$1^\circ$. $2p>d$;

\noindent
$2^\circ$. $g^0=\underline{g}$, \textit{i.~e., representations} (5.43) \textit{are satisfied}.

\noindent
\textit{Then Condition} 7.9 \textit{is satisfied. Moreover, $\|\Lambda\|_{L_\infty}$ and the multiplier norm} (7.50)
\textit{are controlled in terms of $m$, $n$, $d$, $p$,
$\alpha_0$, $\alpha_1$, $\|g\|_{L_\infty}$, $\|g^{-1}\|_{L_\infty}$, and the parameters of the lattice $\Gamma$.}

\smallskip\noindent\textbf{Proof.}
Since $\Lambda\in \wt{H}^p(\Omega)$, in the case where $2p >d$
it follows from the Sobolev embedding theorem and from theorem of [MSh, Subsection 1.3.3]
that Condition 7.9 is satisfied.
Herewith, $\|\Lambda\|_{L_\infty}$ and $M_\Lambda$ are estimated by
$C \|\Lambda\|_{H^p(\Omega)}$, where $C$ depends on $m$, $n$, $d$, $p$, and the domain $\Omega$.
Taking estimate (5.51) into account, we prove the first statement.

Let us prove the second statement. We assume that $2p \leqslant d$ (otherwise, the first statement can be applied).
Suppose that $g^0=\underline{g}$. By Remark 5.6, we have
$\wt{g}= g(b(\D)\Lambda+\1_m)= g^0$.
Then $\Lambda \in \wt{H}^p(\Omega)$ is the $\Gamma$-periodic solution of the problem
$$
b(\D)^* b(\D)\Lambda(\x)= b(\D)^* g(\x)^{-1} g^0,\quad \intop_\Omega \Lambda(\x)\,d\x =0.
\eqno(7.51)
$$
The operator $b(\D)^* b(\D)$ is a matrix elliptic operator with constant coefficients.
Therefore, the solution of the problem (7.51) can be described in terms of the Fourier coefficients:
$$
\wh{\Lambda}_{\mathbf 0}=0;\quad
\wh{\Lambda}_{\mathbf s} = \left(b({\mathbf s})^* b({\mathbf s})\right)^{-1} b({\mathbf s})^* \wh{(g^{-1})}_{\mathbf s} g^0,
\quad {\mathbf 0} \ne {\mathbf s} \in \wt{\Gamma}.
$$
Since $g^{-1} g^0 \in L_\infty \subset L_q(\Omega)$ for any $q<\infty$, from the well known Marcinkiewicz theorem about
the multipliers for the Fourier series (see [Ma]) it follows that $\Lambda \in \wt{W}^p_q(\Omega)$ for any $q<\infty$.
Let us fix $q$ such that $pq >d$ (for instance, $q=p^{-1}(d+1)$).
By the Marcinkiewicz theorem, the norm $\|\Lambda\|_{W_q^p(\Omega)}$ is controlled in terms of
$m$, $n$, $d$, $p$, $\alpha_0$, $\alpha_1$, $\|g\|_{L_\infty}$, and $\|g^{-1}\|_{L_\infty}$.
Next, by the Sobolev embedding theorem and Corollary 1 of [MSh, Subsection 1.3.4],
relation $\Lambda \in \wt{W}^p_q(\Omega)$ ensures that Condition 7.9 is satisfied.
Herewith, $\|\Lambda\|_{L_\infty}$ and $M_\Lambda$ are estimated by $C \|\Lambda\|_{W^p_q(\Omega)}$,
where $C$ depends on $m$, $n$, $d$, $p$, and the domain $\Omega$.
This completes the proof of the second statement.
$\square$

\smallskip
Now we estimate the operator $b(\D) (I-\Pi)(A^0 - \zeta \eps^{2p}I)^{-1}$ in the \hbox{$(L_2 \to H^p)$}-norm.

\smallskip\noindent\textbf{Lemma 7.11.}
\textit{For $\eps>0$ and $\zeta \in \C \setminus \R_+$ we have
$$
\| b(\D) (I-\Pi)(A^0 - \zeta \eps^{2p}I)^{-1} \|_{L_2(\R^d)\to H^p(\R^d)}
\leqslant C_{11} c(\varphi),
\eqno(7.52)
$$
where $C_{11}=2 \alpha_1^{1/2} c_*^{-1} \left(1+ r_0^{-2}\right)^{p/2}$. }

\smallskip\noindent\textbf{Proof.}
Using (4.4), (5.34), and (7.3), we obtain
$$
\begin{aligned}
&\| b(\D) (I-\Pi)(A^0 + |\zeta| \eps^{2p}I)^{-1} \|_{L_2(\R^d)\to H^p(\R^d)}
\\
&= \sup_{\bxi \in \R^d} (1- \chi_{\wt{\Omega}}(\bxi)) (1+ |\bxi|^2)^{p/2} \left|b(\bxi)\left(b(\bxi)^* g^0 b(\bxi) + |\zeta| \eps^{2p} \1_n\right)^{-1}\right|
\\
&\leqslant  \alpha_1^{1/2} \sup_{|\bxi| \geqslant r_0}  (1+ |\bxi|^2)^{p/2} |\bxi|^p \left(c_* |\bxi|^{2p} + |\zeta| \eps^{2p} \right)^{-1}
\leqslant  \alpha_1^{1/2} c_*^{-1} \left(1+ r_0^{-2}\right)^{p/2}.
\end{aligned}
\eqno(7.53)
$$
Obviously,
$$
\begin{aligned}
&\| (A^0 + |\zeta| \eps^{2p}I) (A^0 - \zeta \eps^{2p}I)^{-1}\|_{L_2(\R^d)\to L_2(\R^d)}
\leqslant
\sup_{x\geqslant 0} (x+ |\zeta| \eps^{2p}) |x - \zeta \eps^{2p}|^{-1}
\\
&=
\sup_{y\geqslant 0} (y+  1) |y - \wh{\zeta} |^{-1} \leqslant 2 c(\varphi).
\end{aligned}
\eqno(7.54)
$$
Relations (7.53) and (7.54) imply (7.52).
$\square$

\smallskip\noindent\textbf{Proposition 7.12.} \textit{Suppose that the assumptions of Theorem} 7.6
\textit{and Condition} 7.9 \textit{are satisfied. Then for $\eps>0$ we have}
$$
\|  \Lambda b(\D)(I -\Pi)(A^0 -\zeta \eps^{2p}I)^{-1} \|_{L_2(\R^d)\to L_2(\R^d)} \leqslant \mathcal{C}_5 c(\varphi),
\eqno(7.55)
$$
$$
\| A^{1/2} \left(\Lambda b(\D)(I -\Pi)(A^0 -\zeta \eps^{2p}I)^{-1} \right) \|_{L_2(\R^d)\to L_2(\R^d)} \leqslant \mathcal{C}_6 c(\varphi),
\eqno(7.56)
$$
$$
\|  \wt{g} b(\D)(I -\Pi)(A^0 -\zeta \eps^{2p} I)^{-1}  \|_{L_2(\R^d)\to L_2(\R^d)} \leqslant \mathcal{C}_7 c(\varphi).
\eqno(7.57)
$$
\textit{The constants $\mathcal{C}_5$, $\mathcal{C}_6$, and $\mathcal{C}_7$
depend only on $m$, $d$, $p$, $\alpha_0$, $\alpha_1$, $\|g\|_{L_\infty}$, $\|g^{-1}\|_{L_\infty}$,
the parameters of the lattice $\Gamma$, and also on $M_\Lambda$ and $\|\Lambda\|_{L_\infty}$.}

\smallskip\noindent\textbf{Proof.}
Estimate (7.55) with ${\mathcal C}_5= \|\Lambda\|_{L_\infty} C_{11}$ follows from Condition 7.9 and estimate (7.52).

To prove (7.56), note that
$$
\begin{aligned}
&\| A^{1/2} \left(\Lambda b(\D)(I -\Pi)(A^0 -\zeta \eps^{2p}I)^{-1} \right) \|_{L_2(\R^d)\to L_2(\R^d)}
\\
 &\leqslant \|g\|_{L_\infty}^{1/2} \alpha_1^{1/2}
\| \Lambda b(\D)(I -\Pi)(A^0 -\zeta \eps^{2p}I)^{-1} \|_{L_2(\R^d)\to H^p(\R^d)}.
\end{aligned}
\eqno(7.58)
$$
Combining Condition~7.9 and inequalities (7.52), (7.58), we obtain (7.56) with the constant
${\mathcal C}_6 = \alpha_1^{1/2} \|g\|_{L_\infty}^{1/2} M_\Lambda C_{11}$.

To prove (7.57), note that,  by Lemma 1 of [MSh, Subsection 1.3.2],
from Condition 7.9 it follows that $b(\D)\Lambda$ is a multiplier from $H^p(\R^d;\C^m)$ to $L_2(\R^d;\C^m)$, and
its multiplier norm is controlled in terms of $\alpha_1$, $\|\Lambda\|_{L_\infty}$, and $M_\Lambda$:
$$
\| [b(\D)\Lambda] \|_{H^p(\R^d)\to L_2(\R^d)} \leqslant {\mathfrak C}_\Lambda = {\mathfrak C}_\Lambda(\alpha_1,\|\Lambda\|_{L_\infty}, M_\Lambda).
$$
Then the matrix-valued function $\wt{g}=g(b(\D)\Lambda + \1_m)$ is a multiplier from
$H^p(\R^d;\C^m)$ to $L_2(\R^d;\C^m)$, and
$\| [\wt{g}] \|_{H^p(\R^d)\to L_2(\R^d)} \leqslant \|g\|_{L_\infty}({\mathfrak C}_\Lambda +1)$.
Together with (7.52) this implies (7.57) with $\mathcal{C}_7 = \|g\|_{L_\infty}({\mathfrak C}_\Lambda +1) C_{11}$.
$\square$

Now, Theorem 7.6 and Proposition 7.12 imply the following result.

\smallskip\noindent\textbf{Theorem 7.13.} \textit{Suppose that the assumptions of Theorem} 7.6
\textit{and Condition} 7.9 \textit{are satisfied. Then for $\eps>0$ we have}
$$
\begin{aligned}
&\| (A - \zeta \eps^{2p}I)^{-1} - (I + \Lambda b(\D))(A^0 -\zeta \eps^{2p}I)^{-1} \|_{L_2(\R^d)\to L_2(\R^d)}
\\
&\leqslant {\mathcal C}_2 c(\varphi)^2 \eps^{1-2p} |\zeta|^{1/2p -1} + {\mathcal C}_5 c(\varphi),
\end{aligned}
$$
$$
\begin{aligned}
&\| A^{1/2}\left((A - \zeta \eps^{2p}I)^{-1} - (I + \Lambda b(\D))(A^0 -\zeta \eps^{2p}I)^{-1}\right) \|_{L_2(\R^d)\to L_2(\R^d)}
\\
&\leqslant {\mathcal C}_3 c(\varphi)^2 \eps^{1-p} |\zeta|^{1/2p -1/2}+ {\mathcal C}_6 c(\varphi),
\end{aligned}
$$
$$
\begin{aligned}
&\| g b(\D)(A - \zeta \eps^{2p}I)^{-1} - \wt{g} b(\D) (A^0 -\zeta \eps^{2p}I)^{-1} \|_{L_2(\R^d)\to L_2(\R^d)}
\\
&\leqslant {\mathcal C}_4 c(\varphi)^2 \eps^{1-p} |\zeta|^{1/2p -1/2} + {\mathcal C}_7 c(\varphi).
\end{aligned}
$$
\textit{The constants ${\mathcal C}_2$, ${\mathcal C}_3$, and ${\mathcal C}_4$ depend only on
$m$, $d$, $p$, $\alpha_0$, $\alpha_1$, $\|g\|_{L_\infty}$, $\|g^{-1}\|_{L_\infty}$, and the parameters of the lattice $\Gamma$.
The constants ${\mathcal C}_5$, ${\mathcal C}_6$, and ${\mathcal C}_7$ depend on the same parameters and also on $\|\Lambda\|_{L_\infty}$ and $M_\Lambda$.}

\section{Homogenization of the operator $A_\eps$}

\subsection{Approximation of the resolvent of the operator $A_\eps$ in the operator norm in $L_2(\R^d;\C^n)$}
For any $\Gamma$-periodic function $\varphi(\x)$ in $\R^d$ we denote
$$
\varphi^\eps(\x):= \varphi(\eps^{-1}\x),\quad \eps>0.
$$

In $L_2(\R^d;\C^n)$, we consider the operator $A_\eps$, $\eps>0$, given formally by the differential expression
$$
A_\eps = b(\D)^* g^\eps(\x) b(\D),\quad \eps>0.
\eqno(8.1)
$$
As usual, the precise definition of the operator $A_\eps$ is given in terms of the corresponding
closed quadratic form
$$
a_\eps [\u,\u] = \intop_{\R^d} \langle g^\eps(\x) b(\D)\u, b(\D)\u \rangle \,d\x,
\quad \u \in H^p(\R^d;\C^n).
$$
The form $a_\eps$ is subject to the following estimates similar to (4.8):
$$
\alpha_0 \|g^{-1}\|^{-1}_{L_\infty} \intop_{\R^d} |\bxi|^{2p} |\wh{\u}(\bxi)|^2\,d\bxi
\leqslant a_\eps [\u,\u] \leqslant \alpha_1 \|g\|_{L_\infty} \intop_{\R^d} |\bxi|^{2p} |\wh{\u}(\bxi)|^2\,d\bxi.
\eqno(8.2)
$$

For small $\eps$ the coefficients of the operator (8.1) oscillate rapidly.
A typical homogenization problem as applied to the operator (8.1)
is to approximate its resolvent for small $\eps$. Using the results of \S 7 and the scaling transformation, we deduce
theorems about approximation of the resolvent $(A_\eps - \zeta I)^{-1}$ for $\zeta \in \C \setminus \R_+$.

Let $T_\eps$ be the unitary scaling transformation in $L_2(\R^d;\C^n)$ given by
$$
(T_\eps \u)(\x) := \eps^{d/2} \u(\eps \x).
$$
It is easily seen that
$$
A_\eps = \eps^{-2p} T_\eps^* A T_\eps,
$$
where $A$ is the operator (4.1). Hence,
$$
(A_\eps - \zeta I)^{-1} = \eps^{2p} T_\eps^* (A - \zeta \eps^{2p} I)^{-1} T_\eps.
\eqno(8.3)
$$
A similar identity is true for the operator $A^0$:
$$
(A^0 - \zeta I)^{-1} = \eps^{2p} T_\eps^* (A^0 - \zeta \eps^{2p}I)^{-1} T_\eps.
\eqno(8.4)
$$
Subtracting (8.4) from (8.3) and using that the operator $T_\eps$ is unitary, we obtain
$$
\begin{aligned}
& \| (A_\eps - \zeta I)^{-1} - (A^0 - \zeta I)^{-1} \|_{L_2(\R^d) \to L_2(\R^d)}
 \\
 & = \eps^{2p}\| (A - \zeta \eps^{2p}I)^{-1} - (A^0 - \zeta \eps^{2p}I)^{-1}\|_{L_2(\R^d) \to L_2(\R^d)}.
\end{aligned}
\eqno(8.5)
$$
Theorem 7.5 together with (8.5) imply the following result.

\smallskip\noindent\textbf{Theorem 8.1.} \textit{Let $A_\eps$ be the operator} (8.1),
\textit{and let $A^0$ be the effective operator} (5.33). \textit{Let $\zeta = |\zeta|e^{i\varphi} \in \C \setminus \R_+$,
 and let $c(\varphi)$ be defined by} (7.22). \textit{For $\eps>0$ we have}
$$
\| (A_\eps - \zeta I)^{-1} - (A^0 -\zeta I)^{-1} \|_{L_2(\R^d)\to L_2(\R^d)}
\leqslant {\mathcal C}_1  c(\varphi)^2 \eps |\zeta|^{1/2p -1}.
\eqno(8.6)
$$
\textit{The constant ${\mathcal C}_1$ depends only on $d$, $p$, $\alpha_0$, $\alpha_1$, $\|g\|_{L_\infty}$, $\|g^{-1}\|_{L_\infty}$, and
the parameters of the lattice $\Gamma$.}

\subsection{Approximation of the resolvent of the operator $A_\eps$ in the energy norm}

Now, using Theorem 7.6, we obtain approximation of the resolvent $(A_\eps - \zeta I)^{-1}$
in the norm of operators acting from $L_2(\R^d;\C^n)$ to the Sobolev space $H^p(\R^d;\C^n)$,
and also approximation of the operator $g^\eps b(\D) (A_\eps - \zeta I)^{-1}$
(corresponding to the "flux") in the norm of operators acting from $L_2(\R^d;\C^n)$ to $L_2(\R^d;\C^m)$.

Let $\Pi_\eps$ be the pseudodifferential operator in $L_2(\R^d;\C^n)$ with the symbol $\chi_{\wt{\Omega}/\eps}(\bxi)$, i.~e.,
$$
(\Pi_\eps \u)(\x) = (2\pi)^{-d/2} \intop_{\wt{\Omega}/\eps} e^{i \langle\x ,\bxi \rangle} \wh{\u}(\bxi)\,d\bxi.
\eqno(8.7)
$$
The operators (7.3) and (8.7) satisfy the following identity
$$
\Pi_\eps = T_\eps^* \Pi T_\eps.
\eqno(8.8)
$$

We put
$$
K(\zeta;\eps):= \Lambda^\eps b(\D) (A^0 - \zeta I)^{-1} \Pi_\eps.
\eqno(8.9)
$$
The operator (8.9) is called a \textit{corrector}; this operator is a continuous mapping of
$L_2(\R^d;\C^n)$ to $H^p(\R^d;\C^n)$.

\smallskip\noindent\textbf{Theorem 8.2.} \textit{Suppose that the assumptions of Theorem} 8.1
\textit{are satisfied. Let $\Pi_\eps$ be the operator} (8.8), \textit{let $K(\zeta;\eps)$ be given by} (8.9),
\textit{and let $\wt{g}(\x)$ be the matrix-valued function} (5.30).
\textit{Then for $\eps>0$ we have}
$$
\begin{aligned}
&\| (A_\eps - \zeta I)^{-1} - (A^0 -\zeta I)^{-1} - \eps^{2p} K(\zeta;\eps)\|_{L_2(\R^d)\to H^p(\R^d)}
\\
&\leqslant  \eps c(\varphi)^2 |\zeta|^{1/2p} \left({\mathcal C}'  |\zeta|^{-1} + {\mathcal C}''  |\zeta|^{-1/2}\right),
\end{aligned}
\eqno(8.10)
$$
$$
\begin{aligned}
&\| g^\eps b(\D)(A_\eps - \zeta I)^{-1} - \wt{g}^\eps b(\D)(A^0 -\zeta I)^{-1} \Pi_\eps \|_{L_2(\R^d)\to L_2(\R^d)}
\\
&\leqslant  \eps c(\varphi)^2 {\mathcal C}_4  |\zeta|^{1/2p -1/2}.
\end{aligned}
\eqno(8.11)
$$
\textit{The constants ${\mathcal C}'$, ${\mathcal C}''$, and ${\mathcal C}_4$ depend only on
$m$, $d$, $p$, $\alpha_0$, $\alpha_1$, $\|g\|_{L_\infty}$, $\|g^{-1}\|_{L_\infty}$, and the parameters of the lattice $\Gamma$.}

\smallskip\noindent\textbf{Proof.}
Similarly to (8.3), by (8.8),
$$
K(\zeta;\eps)  =
\eps^{p} T_\eps^* \Lambda b(\D)(A^0 - \zeta \eps^{2p} I)^{-1} \Pi T_\eps.
\eqno(8.12)
$$

From (8.3), (8.4), and (8.12) it follows that
$$
\begin{aligned}
& \| (A_\eps - \zeta I)^{-1} - (A^0 - \zeta I)^{-1} - \eps^p K(\zeta;\eps)\|_{L_2(\R^d) \to L_2(\R^d)}
 \\
 & = \eps^{2p}\| (A - \zeta \eps^{2p}I)^{-1} - (I + \Lambda b(\D)\Pi)(A^0 - \zeta \eps^{2p}I)^{-1}  \|_{L_2(\R^d) \to L_2(\R^d)}.
\end{aligned}
\eqno(8.13)
$$
By (7.28) and (8.13), for $\eps>0$ we have
$$
\| (A_\eps - \zeta I)^{-1} - (A^0 - \zeta I)^{-1} - \eps^p K(\zeta;\eps)\|_{L_2(\R^d) \to L_2(\R^d)}
\leqslant {\mathcal C}_2 c(\varphi)^2 \eps |\zeta|^{1/2p-1}.
\eqno(8.14)
$$

Similarly,
$$
\begin{aligned}
& \| A_\eps^{1/2}\left((A_\eps - \zeta I)^{-1} - (A^0 - \zeta I)^{-1} - \eps^p K(\zeta;\eps)\right)\|_{L_2(\R^d) \to L_2(\R^d)}
 \\
 & = \eps^{p}\|A^{1/2}\left( (A - \zeta \eps^{2p}I)^{-1} - (I + \Lambda b(\D)\Pi)(A^0 - \zeta \eps^{2p}I)^{-1}\right)  \|_{L_2(\R^d) \to L_2(\R^d)}.
\end{aligned}
$$
Together with (7.29) this yields
$$
\begin{aligned}
&\| A_\eps^{1/2}\left((A_\eps - \zeta I)^{-1} - (A^0 - \zeta I)^{-1} - \eps^p K(\zeta;\eps)\right)\|_{L_2(\R^d) \to L_2(\R^d)}
\\
&\leqslant {\mathcal C}_3  c(\varphi)^2 \eps |\zeta|^{1/2p-1/2}.
\end{aligned}
\eqno(8.15)
$$

Since $(1+|\bxi|^2)^p\leqslant 2^{p-1}(1+|\bxi|^{2p})$, taking the lower estimate (8.2) into account, for any
$\u\in H^p(\R^d;\C^n)$ we have
$$
\begin{aligned}
&\| \u \|^2_{H^p(\R^d)}=
\intop_{\R^d} (1+|\bxi|^{2})^p |\wh{\u}(\bxi)|^2\,d\bxi
\leqslant
2^{p-1} \intop_{\R^d} (1+|\bxi|^{2p}) |\wh{\u}(\bxi)|^2\,d\bxi
\\
&\leqslant 2^{p-1} \left(\|\u\|_{L_2(\R^d)}^2 + \alpha_0^{-1}\|g^{-1}\|_{L_\infty}\|A_\eps^{1/2} \u\|^2_{L_2(\R^d)}\right).
\end{aligned}
$$
Combining this with  (8.14) and (8.15), we deduce
(8.10) with the constants
$$
{\mathcal C}'= 2^{(p-1)/2} {\mathcal C}_2,\quad {\mathcal C}'' = 2^{(p-1)/2} {\mathcal C}_3 \alpha_0^{-1/2}\|g^{-1}\|^{1/2}_{L_\infty}.
$$

Inequality (8.11) follows from (7.30) with the help of the scaling transformation. $\square$

\smallskip\noindent\textbf{Remark 8.3.}
1) For a fixed $\zeta \in \C \setminus \R_+$, estimates of Theorems 8.1 and 8.2 are of sharp order $O(\eps)$.
For large $|\zeta|$ the order improves due to the presence of the factors $|\zeta|^{-s}$ (with $s>0$)
in the right-hand sides. 2) Estimates (8.6), (8.10), and (8.11) are uniform with respect to $\varphi$ in any sector of the form
$\{\zeta=|\zeta| e^{i\varphi} \in \C:\ \varphi_0 \leqslant \varphi \leqslant 2\pi - \varphi_0\}$
with arbitrarily small $\varphi_0$.

\subsection{Special cases}

If $g^0=\overline{g}$, then $\Lambda =0$ and the corrector (8.9) is equal to zero.
In this case (8.10) simplifies.

\smallskip\noindent\textbf{Proposition 8.4.} \textit{Suppose that the assumptions of Theorem} 8.1 \textit{are satisfied.
Suppose that $g^0=\overline{g}$} (\textit{i.~e., conditions} (5.41)
\textit{are satisfied}). \textit{Then for $\eps>0$ we have}
$$
\begin{aligned}
&\|  (A_\eps - \zeta I)^{-1} - (A^0 - \zeta I)^{-1}  \|_{L_2(\R^d)\to H^p(\R^d)}
\\
&\leqslant \eps  c(\varphi)^2  |\zeta|^{1/2p}\left( {\mathcal C}'|\zeta|^{-1} + {\mathcal C}''|\zeta|^{-1/2}\right).
\end{aligned}
$$

\smallskip
The following statement is deduced from Proposition 7.8 by the scaling transformation.

\smallskip\noindent\textbf{Proposition 8.5.} \textit{Suppose that the assumptions of Theorem} 8.1
\textit{are satisfied. If $g^0=\underline{g}$} (\textit{i.~e., representations} (5.43)
\textit{are true}), \textit{then for $\eps>0$ we have}
$$
\begin{aligned}
&\| g^\eps b(\D) (A -\zeta I)^{-1} - g^0 b(\D)(A^0 -\zeta I)^{-1}  \|_{L_2(\R^d)\to L_2(\R^d)}
\\
&\leqslant {\mathcal C}_4^\circ c(\varphi)^2 \eps |\zeta|^{1/2p - 1/2}.
\end{aligned}
$$

\subsection{Removal of the smoothing operator}

Now, we suppose that Condition 7.9 is satisfied.
Then instead of the corrector (8.9) one can use the operator
$$
K^0(\zeta;\eps):= \Lambda^\eps b(\D) (A^0 - \zeta I)^{-1},
\eqno(8.16)
$$
which in this case is a continuous mapping of $L_2(\R^d;\C^n)$ to $H^p(\R^d;\C^n)$.
Note that (8.16) is the traditional corrector used in the homogenization theory.

The following result is deduced from Theorem 7.13 by the scaling transformation (cf. the proof of Theorem 8.2).

\smallskip\noindent\textbf{Theorem 8.6.} \textit{Suppose that the assumptions of Theorem} 8.1
\textit{and Condition} 7.9 \textit{are satisfied. Let $K^0(\zeta;\eps)$ be given by} (8.16), \textit{and let
$\wt{g}(\x)$ be the matrix-valued function} (5.30).
\textit{Then for $0< \eps \leqslant 1$ we have}
$$
\begin{aligned}
&\| (A_\eps - \zeta I)^{-1} - (A^0 -\zeta I)^{-1} - \eps^p K^0(\zeta;\eps) \|_{L_2(\R^d)\to H^p(\R^d)}
\\
&\leqslant
\eps  c(\varphi)^2  |\zeta|^{1/2p}\left( {\mathcal C}'|\zeta|^{-1} + {\mathcal C}''|\zeta|^{-1/2}\right)
+ {\mathcal C}_8 \eps^p c(\varphi),
\end{aligned}
\eqno(8.17)
$$
$$
\begin{aligned}
&\| g^\eps b(\D)(A - \zeta I)^{-1} - \wt{g}^\eps b(\D) (A^0 -\zeta I)^{-1} \|_{L_2(\R^d)\to L_2(\R^d)}
\\
&\leqslant \eps {\mathcal C}_4 c(\varphi)^2 |\zeta|^{1/2p -1/2} + {\mathcal C}_7 \eps^p c(\varphi).
\end{aligned}
\eqno(8.18)
$$
\textit{The constants ${\mathcal C}'$, ${\mathcal C}''$, and ${\mathcal C}_4$ depend only on
$m$, $d$, $p$, $\alpha_0$, $\alpha_1$, $\|g\|_{L_\infty}$, $\|g^{-1}\|_{L_\infty}$, and the parameters of the lattice $\Gamma$.
The constants ${\mathcal C}_7$ and ${\mathcal C}_8$ depend on the same parameters and also on
$\|\Lambda\|_{L_\infty}$ and $M_\Lambda$.}

\smallskip

Note that in Theorem 8.6 it is assumed that $0< \eps \leqslant 1$. This is
because in the proof of (8.17) it is used that $\eps^{2p}\leqslant \eps^p$.
Besides, estimates (8.17) and (8.18) are interesting for small $\eps$.
The constant ${\mathcal C}_8$ is given by
${\mathcal C}_8 = 2^{(p-1)/2}\left( {\mathcal C}_5 + \alpha_0^{-1/2} \|g^{-1}\|^{1/2}_{L_\infty} {\mathcal C}_6 \right)$.

Comparing Proposition 7.10 and Theorem 8.6, we arrive at the following statement.

\smallskip\noindent\textbf{Corollary 8.7.}
 \textit{Suppose that the assumptions of Theorem} 8.1
 \textit{are satisfied. Let $K^0(\zeta;\eps)$ be the operator} (8.16), \textit{and let $\wt{g}(\x)$ be the matrix-valued function} (5.30).
\textit{Moreover, suppose that at least one of the following two assumptions is satisfied}:

\noindent
$1^\circ$. $2p>d$;

\noindent
$2^\circ$. $g^0=\underline{g}$, \textit{i.~e., representations} (5.43) \textit{are true}.

\noindent
\textit{Then estimates} (8.17) \textit{and} (8.18)
\textit{are true for $0< \eps \leqslant 1$. All the constants in these estimates depend only on
 $m$, $n$, $d$, $p$, $\alpha_0$, $\alpha_1$, $\|g\|_{L_\infty}$, $\|g^{-1}\|_{L_\infty}$, and the parameters of the lattice} $\Gamma$.

\smallskip\noindent\textbf{Remark 8.8.}
1) For fixed $\zeta \in \C \setminus \R_+$ estimates of Theorem 8.6 are of sharp order $O(\eps)$.
2) Estimates (8.17) and (8.18) are uniform with respect to $\varphi$ in any sector of the form
$\{\zeta=|\zeta| e^{i\varphi} \in \C:\ \varphi_0 \leqslant \varphi \leqslant 2\pi - \varphi_0\}$
with arbitrarily small $\varphi_0$.
3) The assumptions of Corollary 8.7 are satisfied in the following cases that are interesting for applications:
a) if $p=2$ and $d=2$ or $d=3$, then $2p>d$; b) if $m=n$, then $g^0=\underline{g}$. For instance, this is the case
for the operator $A_\eps = \Delta g^\eps(\x)\Delta$.

\end{document}